\documentclass{amsart}
\usepackage{amsxtra}
\usepackage{color}
\theoremstyle{plain}

\newcommand{\cleqn}{\setcounter{equation}{0}}
\newcommand{\clth}{\setcounter{theorem}{0}}
\newcommand {\sectionnew}[1]{\section{#1}\cleqn\clth}
\newtheorem{theorem}{Theorem}[section]
\newtheorem{lemma}[theorem]{Lemma}
\newtheorem{definition-theorem}[theorem]{Definition-Theorem}
\newtheorem{proposition}[theorem]{Proposition}
\newtheorem{corollary}[theorem]{Corollary}
\newtheorem{definition}[theorem]{Definition}
\newtheorem{example}[theorem]{Example}
\newtheorem{remark}[theorem]{Remark}
\newtheorem{notation}[theorem]{Notation}
\newtheorem{assumption}[theorem]{Assumption}
\newtheorem{lemma-definition}[theorem]{Lemma-Definition}
\newtheorem{lemma-notation}[theorem]{Lemma-Notation}
\newtheorem{question}[theorem]{Question}
\newtheorem{remark-definition}[theorem]{Remarks-Definition}
\newtheorem{notation-remark}[theorem]{Notation-Remarks}
\newcommand \bth[1] { \begin{theorem}\label{t#1} }
\newcommand \ble[1] { \begin{lemma}\label{l#1} }

\newcommand \bpr[1] { \begin{proposition}\label{p#1} }
\newcommand \bco[1] { \begin{corollary}\label{c#1} }
\newcommand \bde[1] { \begin{definition}\label{d#1}\rm }
\newcommand \bex[1] { \begin{example}\label{e#1}\rm }
\newcommand \bre[1] { \begin{remark}\label{r#1}\rm }

\newcommand \bnota[1] {\begin{notation}\label{n#1}\rm }
\newcommand \bas[1] { \begin{assumption}\label{a#1}\rm }
\newcommand \bld[1] { \begin{lemma-definition}\label{ld#1} }

\newcommand \bqu[1] { \begin{question}\label{q#1}\rm }

\newcommand {\eth} { \end{theorem} }
\newcommand {\ele} { \end{lemma} }

\newcommand {\epr} { \end{proposition} }
\newcommand {\eco} { \end{corollary} }
\newcommand {\ede} { \end{definition} }
\newcommand {\eex} { \end{example} }
\newcommand {\ere} { \end{remark} }
\newcommand {\enota} { \end{notation} }
\newcommand {\eas} {\end{assumption}}
\newcommand {\eld}{ \end{lemma-definition} }

\newcommand {\equ} {\end{question}}
\newcommand \thref[1]{Theorem \ref{t#1}}
\newcommand \leref[1]{Lemma \ref{l#1}}
\newcommand \prref[1]{Proposition \ref{p#1}}

\newcommand \deref[1]{Definition \ref{d#1}}
\newcommand \exref[1]{Example \ref{e#1}}
\newcommand \reref[1]{Remark \ref{r#1}}
\newcommand \lb[1]{\label{#1}}

\newcommand \ldref[1]{Lemma-Definition \ref{ld#1}}

\def \d {{\partial}}   %differentials and partials

\def \Rset {{\mathbb R}}         %mathsets
\def \Cset {{\mathbb C}}

\def \F  {{\mathcal{F}}}

\def \la {\lambda}

\def \ra  {\rightarrow}           %maps

\def \la {\langle}
\def \ra {\rangle}
\def \Ad { {\mathrm{Ad}} }

\def \g  {\mathfrak{g}}   % Lie algebra letters
\def \h  {\mathfrak{h}}
\def \f  {\mathfrak{f}}

\def \b  {\mathfrak{b}}
\def \p  {\mathfrak{p}}

\def \a  {\mathfrak{a}}

\def \q  {\mathfrak{q}}
\def \d  {\mathfrak{d}}

\def \c  {\mathfrak{c}}

\def \k {\mathfrak{k}}
\DeclareMathOperator \ad { {\mathrm{ad}} }

\newcommand{\beqa}{\begin{eqnarray*}}                     %added by Lu
\newcommand{\eeqa}{\end{eqnarray*}}
\def \hs {\hspace{.2in}}
\def \lara {\la \, , \, \ra}

\def \lara {\la \, , \, \ra}

\def \piG {\pi_G}

\def \lam {\lambda}
\def \V {\mathcal V}
\def \vr {\varrho}

\def \pist {\pi_{\rm st}}

\def \sQ {{\scriptscriptstyle Q}}
\def \sD {{\scriptscriptstyle D}}
\def \sF {{\scriptscriptstyle F}}
\def \sG {{\scriptscriptstyle G}}
\def \sX {{\scriptscriptstyle X}}
\def \sY {{\scriptscriptstyle Y}}
\def \sZ {{\scriptscriptstyle Z}}

\def \sM {{\scriptscriptstyle M}}

\def \piG {{\pi_{\scriptscriptstyle G}}}
\def \piD {{\pi_{\scriptscriptstyle D}}}
\def \piGs {{\pi_{{\scriptscriptstyle G}^*}}}
\def \piX {{\pi_{\scriptscriptstyle X}}}
\def \piY {{\pi_{\scriptscriptstyle Y}}}
\def \piZ {{\pi_{\scriptscriptstyle Z}}}

\def \Cset {{\mathbb C}}
\def \lrw {\longrightarrow}
\def \Pist {\Pi_{\rm st}}
\def \Lam {{\Lambda}}

\def \sV {{\scriptscriptstyle V}}
\def \sGQ {{\scriptscriptstyle G/Q}}
\def \sGQp {{\scriptscriptstyle{G/Q_+}}}

\def \sH {{\scriptscriptstyle H}}
\def \sJ {{\scriptscriptstyle J}}

\def \ot {\otimes}

\def \Fbb {\mathbb F}
\def \Xbb {{\mathbb{X}}}

\def \sFbb {\scriptscriptstyle{\mathbb F}}

\def \sXbb {\scriptscriptstyle{\mathbb X}}
\def \sYbb {\scriptscriptstyle{\mathbb Y}}

\def \swF {\scriptscriptstyle{\widetilde{F}}}
\def \wF {\widetilde{F}}
\def \swFbb {\scriptscriptstyle{\widetilde{{\mathbb{F}}}}}
\def \wFbb {\widetilde{{\mathbb{F}}}}
\def \sW {\scriptscriptstyle{W}}
\begin{document}

\setlength{\baselineskip}{1.2\baselineskip}
%%%%%%%%%%%%%%%%%%%%%%%%%%%%%%%%%%%%%%%%%%%%%%%%%%%%%%%%%%%%%%%%%%%%%%%%%%%
%%%%%%%%%%%%%%%%%%%%%%    Title    %%%%%%%%%%%%%%%%%%%%%%%%%%%%%%%%%%%%%%%%
\title[Mixed product Poisson structures]
{Mixed product Poisson structures associated to Poisson Lie groups and Lie bialgebras}
\author{Jiang-Hua Lu}
\address{
Department of Mathematics   \\
The University of Hong Kong \\
Pokfulam Road               \\
Hong Kong}
\email{jhlu@maths.hku.hk}
\author{Victor Mouquin}
\address{
Department of Mathematics   \\
University of Toronto \\
Toronto, Canada}               
\email{mouquinv@math.toronto.edu}
\date{}
\begin{abstract} We introduce  and study some 
mixed product Poisson structures on product manifolds associated to Poisson Lie groups and Lie bialgebras. For
quasitriangular Lie bialgebras, our construction is equivalent to that of fusion products of quasi-Poisson $G$-manifolds introduced by Alekseev, Kosmann-Schwarzbach, and Meinrenken. Our primary examples include four series of 
holomorphic Poisson structures on products of flag varieties and related spaces of complex semi-simple Lie groups.
\end{abstract}
\maketitle
%\tableofcontents
%%%%%%%%%%%%%%%%%%%%   Introduction   %%%%%%%%%%%%%%%%%%%%%%%%%%%%%%%%%%%%%%%%
\vspace{-.35in}
\sectionnew{Introduction and outlines of main results}\lb{sec-intro}

\subsection{Introduction}\lb{subsec-intro}
This paper is the first of a series of papers devoted to a detailed study of some naturally defined
holomorphic Poisson structures on products of flag varieties and related spaces of complex semi-simple Lie groups
(see $\S$\ref{subsec-moti-intro} and $\S$\ref{subsec-flags} for the definitions of these Poisson structures).
Aspects of the Poisson structures to be investigated include calculations in coordinates, symplectic groupoids over (extended) Bruhat cells, orbits of symplectic leaves under the action of a maximal torus $T$ of $G$, also called the $T$-leaves, log-moment maps for the $T$-action, cluster structures on the $T$-leaves, and quantization. This paper sets up the general framework, namely that of {\it mixed product Poisson structures} and Poisson structures {\it defined by quasitriangular $r$-matrices}, in which the above mentioned aspects of the Poisson structures will be studied in forthcoming papers. Although we are motivated by  particular examples, the theory developed in the paper is applicable to much wider classes of Poisson structures and is of independent interest.

\bde{de-mixed-general}
Given two manifolds $Y_1$ and $Y_2$, by a {\it mixed product Poisson structure} on the product manifold $Y_1 \times Y_2$ we mean a Poisson bi-vector field $\pi$ on $Y_1 \times Y_2$ that projects to well-defined Poisson structures on $Y_1$ and $Y_2$. Given $Y_i$, $1 \leq i \leq n$, where $n \geq 2$, a Poisson structure $\pi$ on the product manifold $Y =Y_1 \times \cdots \times Y_n$ is said to be a {\it mixed product} if the projection of $\pi$ to $Y_i \times Y_j$ is a well-defined mixed product Poisson structure on $Y_i \times Y_j$ for any $1 \leq i < j \leq n$, and in this case, we also call the pair
$(Y_1 \times \cdots \times Y_n, \pi)$ a {\it mixed product Poisson manifold}.
\ede

Mixed product Poisson manifolds are semiclassical analogs of {\it locally factored algebras} introduced by P. Etingof and D. Kazhdan in \cite{EK:III}. When each $Y_j$ is a finite dimensional vector space, {\it linear} mixed product Poisson structures on the direct sum vector space $Y = Y_1 \oplus \cdots \oplus Y_n$ are in one-to-one correspondence with Lie algebra structures on 
its dual space $Y^* \cong Y_1^* \oplus \cdots \oplus Y_n^*$ such that each $Y_i^* \hookrightarrow Y^*$ is a Lie subalgebra and $[Y_i^*, Y_j^*]\subset Y_i^* \oplus Y_j^*$ for each pair $i \neq j$. In this paper, we study a class of mixed product Poisson
structures in the category of Poisson manifolds with Poisson actions by Lie bialgebras.

The constructions in this paper are all special cases of {\it twists} of Lie bialgebras and Poisson actions thereof as introduced by Drinfeld \cite{dr:quasi} in the context of quasi-Hopf algebras. More specifically, given
Lie bialgebras 
$(\g_j, \delta_{\g_j})$, $1 \leq j \leq n$, we introduce ($\S$\ref{subsec-mixed-twists}) the notion of 
{\it mixed twists} of their direct product Lie bialgebra, which are twists by  
{\it mixed twisting elements} of the direct product Lie bialgebra structure on 
$\g = \g_1 \oplus \cdots \oplus \g_n$. 
If $(Y_j, \pi_j, \lam_j)$ is a $(\g_j, \delta_{\g_j})$-Poisson space for each $1 \leq j \leq n$
(see \deref{de-Poisson-space}), a mixed twisting element $t$ is used to obtain 
a mixed product Poisson structure $\piY$ on the product manifold $Y = Y_1 \times \cdots \times Y_n$
with the properties 1) the projection of $\piY$ to each $Y_j$ is  $\pi_j$,  and 2) the direct product action of $\g$ on $(Y, \piY)$ is now a 
Poisson action of the twist by $t$ of the direct product Lie bialgebra structure on $\g$ (see \prref{pr-mixed-mixed}).

As our first example of mixed twists of direct product Lie bialgebras, we introduce 
in $\S$\ref{subsec-polyuble} the notion of {\it polyubles} of an arbitrary Lie bialgebra, generalizing that of the Drinfeld double of a Lie bialgebra. We then construct mixed product Poisson structures ($\S$\ref{sec-mixed} and $\S$\ref{sec-polyuble}) 
that have Poisson actions by polyuble Lie bialgebras or their duals (\prref{pr-n-mixed-product} and \reref{re-n-extended}). Polyubles for factorizable Lie bialgebras have been previously constructed by V. Fock and A. Rosly \cite{fock-rosly:polyubles, fock-rosly:r-matrix}, from whom we have borrowed the term.

Polyubles of Lie bialgebras are special cases of a more general construction: if
$r \in \g \otimes \g$ is a  quasitriangular $r$-matrix 
on a Lie algebra $\g$,  we construct in $\S$\ref{subsec-mixed-powers}, for each $n \geq 2$, a quasitriangular $r$-matrix $r^{(n)}$ on the direct product Lie algebra $\g^n$, which is a certain alternating sum of $r$ plus a
mixed twisting element of the $n$-fold direct product 
of the quasitriangular Lie bialgebra $(\g, r)$ with itself.
The dual of the quasitriangular Lie bialgebra $(\g^n, r^{(n)})$ 
is precisely the {\it locally factored Lie bialgebra with equal components} constructed from $(\g, r)$ by P. Etingof and D. Kazhdan in \cite[Proposition 1.9]{EK:III}. If $(Y_j, \pi_j, \lam_j)$ is a $(\g, r)$-Poisson space for each $1 \leq j \leq n$
(see \deref{de-Poisson-space}),
one then obtains a mixed product Poisson structure on the product manifold $Y = Y_1 \times \cdots \times Y_n$ such that 
the direct product action of $\g^n$ on $Y$ is a Poisson action for the quasitriangular Lie bialgebra $(\g^n, r^{(n)})$
(see \thref{th-piY-Mix-r}).
In particular, we are naturally lead to the notion of {\it fusion products} of Poisson spaces of a quasitriangular Lie bialgebra (\deref{de-mixed-product-r}), which, under twisting, is equivalent to that of fusion products of quasi-Poisson spaces introduced by Alekseev, Kosmann-Schwarzbach, and Meinrenken \cite{AA} (see $\S$\ref{subsec-mixed-fusion}). 

Another emphasis of our paper is on {\it Poisson structures defined by quasitriangular $r$-matrices}, i.e., 
Poisson structures on a manifold $Y$ of the form $\lam(r)$, where $r$ is a quasitriangular $r$-matrix on a Lie algebra 
$\g$ and $\lam$ is a Lie algebra action of $\g$ on $Y$ (see $\S$\ref{subsec-intro-r} and $\S$\ref{subsec-Poi-r}).
$\S$\ref{sec-mixed-quotient} and $\S$\ref{sec-primary-examples} are devoted to our primary examples, 
which are mixed product Poisson structures arising from quotients of product Poisson Lie groups, which, at the
same time, are defined by quasitriangular $r$-matrices of the form $r^{(n)}$ 
(\thref{th-GQ-n-mixed} and \thref{th-main-GGG-QQQ}). This class of examples
includes the four series of holomorphic Poisson structures on products of flag varieties and related spaces of
complex semi-simple Lie groups (see
$\S$\ref{subsec-moti-intro}).

We also give in $\S$\ref{sec-review} a somewhat detailed review on Lie bialgebras. A partial purpose of the review is
to fix our conventions on signs and constants that are to be used in this and subsequent papers.

\subsection{Poisson structures defined by quasitriangular $r$-matrices}\lb{subsec-intro-r}
It is a simple observation, see \prref{pr-admi-Poi-r} (see also \cite[Section 2.1]{David-Severa:quasi-Hamiltonian-groupoids}), that given a quasitriangular $r$-matrix $r =\sum_i x_i \ot y_i \in \g \ot \g$ on a Lie algebra $\g$ and a left Lie algebra action $\lam$ of $\g$ on a manifold $Y$, the $2$-tensor field 
\[
\lam(r):=\sum_i \lam(x_i) \ot \lam(y_i)
\]
on $Y$ 
is Poisson as long as it is skew-symmetric. In this case,  $(Y, -\lam(r), \lam)$ is a left $(\g, r)$-Poisson space, and we say that the Poisson structure $-\lam(r)$ is {\it defined by the action $\lam$ and the quasitriangular $r$-matrix $r$}.  An equivalent condition for $\lam(r)$ to be skew-symmetric is that the stabilizer subalgebra of $\g$ at each point of $Y$ be coisotropic with respect to the symmetric part of $r$ (see $\S$\ref{subsec-Poi-r}). 

\subsection{Primary Examples}\lb{subsec-moti-intro} 
For a Lie group $G$ and an integer $n \geq 1$, let  $G^n$ act on itself from the right by 
\begin{equation} \lb{eq-Gn-Gn-1}
(g_1, g_2, \ldots, g_n) \cdot (h_1, h_2, \ldots, h_n) = (g_1h_1, h_1^{-1}g_2h_2, \ldots, h_{n-1}^{-1}g_nh_n),
\hs g_i, h_i \in G.
\end{equation}
If $Q_1, \ldots, Q_n$ are closed subgroups of $G$, let
\begin{equation} \lb{eq-Z}
Z = G \times_{Q_1} \cdots \times_{Q_{n-1}} G/Q_n
\end{equation}
be the quotient manifold of $G^n$ by the action of $Q_1 \times \cdots \times Q_n \subset G^n$ given in \eqref{eq-Gn-Gn-1}. Assume now that $(G, \piG)$ is a Poisson Lie group and that $Q_1, \ldots, Q_n$ are closed Poisson Lie subgroups of $(G, \piG)$. Then (see \ldref{ld-GGG-QQQ}) the direct product Poisson structure $\pi^n_{\sG}$ on $G^n$ projects to a well-defined Poisson structure $\pi_{\scriptscriptstyle Z} =\varpi_{\sZ}(\pi_{\sG}^n)$ on $Z$, where $\varpi_{\sZ}: G^n \to Z$ is the quotient map, 
and we refer to $(Z, \pi_\sZ)$ as a {\it quotient Poisson manifold} of the product Poisson Lie group $(G^n, \pi_\sG^n)$.

When $(G, \pist)$ is a {\it standard complex semi-simple Poisson Lie group}, defined by the choice of the triple $(B, B_-, \!\lara_\g)$,
 where $(B, B_-)$ is a pair of opposite Borel subgroups of the connected complex semi-simple Lie group $G$ and $\lara_\g$ a fixed symmetric non-degenerate 
invariant bilinear form on the Lie algebra of $\g$ (see $\S$\ref{subsec-flags}), both $B$ and $B_-$ are Poisson Lie subgroups of $(G, \pist)$. For each integer $n \geq 1$, one thus 
has the quotient Poisson manifolds $(F_n, \pi_n)$ and $(\wF_n, \, \tilde{\pi}_n)$, where
\[
F_n = G \times_B \cdots \times_B G/B  \hs \mbox{and} \hs \wF_n = G \times_B \cdots \times_B G
\]
are both quotients of $G^n$, and $\pi_n = \pi_{\sF_n}$ and $\tilde{\pi}_n = \pi_{\swF_n}$ are the respective projections of $\pi_{\rm st}^n$ to $F_n$ and $\wF_n$.
Associated to $(G, \pist)$ is its Drinfeld double Poisson Lie group $(G \times G, \Pist)$, where $G \times G$ has the
product Lie group structure, and $B \times B_-$ is a Poisson Lie subgroup of $(G \times G, \Pist)$. For each integer $n \geq 1$, one also has the quotient Poisson manifolds $(\Fbb_n, \Pi_n)$  and  $(\wFbb_n, \widetilde{\Pi}_n)$, where
\begin{align*}
\Fbb_n &= (G \times G)\times_{  (B\times B_-)} \cdots \times_{(B \times B_-)} (G \times G)/(B \times B_-),\\
\widetilde{{\mathbb F}}_n & = (G \times G)\times_{(B\times B_-)} \cdots \times_{(B \times B_-)} (G \times G)
\end{align*}
are both quotients of $(G \times G)^n$, and $\Pi_n = \pi_{\sFbb_n}$ and $\widetilde{\Pi}_n = \pi_{\swFbb_n}$ are the respective projections of $\Pi_{\rm st}^n$ to $\Fbb_n$ and $\wFbb_n$.
These four series of Poisson manifolds are the focus of study in this and forthcoming papers.

Returning to the quotient Poisson manifold $(Z, \pi_{\sZ})$ in \eqref{eq-Z} for an arbitrary Poisson Lie group $(G, \piG)$ and
closed Poisson Lie subgroups $Q_1, \ldots, Q_n$, note that the manifold $Z$ is diffeomorphic to the product manifold
$(G/Q_1) \times \cdots \times (G/Q_n)$ via the diffeomorphism
$J_\sZ: Z \to (G/Q_1) \times \cdots \times (G/Q_n)$ given by
\begin{equation}\lb{eq-IZ}
J_{\sZ}(\varpi_{\sZ}(g_1, g_2, \ldots, g_n)) = (g_1Q_1, \; g_1g_2Q_2, \; \ldots, \;g_1g_2\cdots g_nQ_n), \notag \hs 
g_1, \ldots, g_n  \in G.
\end{equation}
We study $\pi_\sZ$ via the Poisson structure $J_\sZ(\pi_\sZ)$  on the 
product manifold
$(G/Q_1) \times \cdots \times (G/Q_n)$. The main results 
(\thref{th-GQ-n-mixed} and \thref{th-main-GGG-QQQ}) on the Poisson structure $J_\sZ(\pi_\sZ)$
are summarized as follows.

\bth{th-A}
1) As Poisson structures on the product manifold $(G/Q_1) \times \cdots \times (G/Q_n)$, 
\[
J_\sZ(\pi_\sZ) =-\sigma \left(r_\d^{(n)}\right),
\]
where $r_\d$ is the quasitriangular $r$-matrix on the double Lie algebra of the Lie bialgebra $(\g, \delta_\g)$ 
of $(G, \piG)$, and  $\sigma$ is the direct product action  of $\d^n$ on
$(G/Q_1)  \times \cdots \times (G/Q_n)$ of naturally defined (see \eqref{eq-sigma-GQ}) actions of $\d$ on each factor $G/Q_j$;

2) If $(\g, \delta_\g)$ has a quasitriangular structure $r \in \g \otimes \g$ (see \deref{de-quasi-r}) and if
$Q_1 = \cdots = Q_n = Q$, where the Lie algebra $\q$ of $Q$ satisfies ${\rm Im}(r_+) \subset \q$ (see
\eqref{eq-r-pm-0}), then  
\[
J_{\sZ}(\pi_{\sZ}) = -\lambda(r^{(n)}),
\]
where $\lam$ is the direct product action of $\g^n$ on  $(G/Q)^n$ induced by left translation.
\eth

\thref{th-A} applies directly to the Poisson manifolds $(F_n, \pi_n)$ and $(\Fbb_n, \Pi_n)$. We can thus 
identify both $\pi_n$ and $\Pi_n$ as mixed product Poisson structures, respectively on 
$(G/B)^n$ and $((G \times G)/(B \times B_-))^n$, that are also defined by quasitriangular $r$-matrices.
Similar results holds for the Poisson structures $\tilde{\pi}_n$ and $\widetilde{\Pi}_n$ by adjusting the
quasitriangular $r$-matrices (see $\S$\ref{subsec-DD}). 

Let $T = B \cap B_-$, a maximal torus  of $G$. Then $T$ acts on $G$ by left translation and on 
$G \times G$ by left translation via $T \cong \{(t, t): t \in T\}$.  The actions of $T$ on
$Z \in \{F_n, \;\wF_n,\; \Fbb_n, \; \wFbb_n\}$, induced by the $T$-actions on the first
factors of $G^n$ and $(G \times G)^n$, preserve the Poisson structures $\pi_\sZ$.
In \cite{Lu-Victor:flags}, a sequel to the present paper, we describe the $T$-leaves and compute the
ranks of all the symplectic leaves in  $(Z, \pi_\sZ)$ for $Z \in \{F_n, \;\wF_n,\; \Fbb_n, \; \wFbb_n\}$ by
first developing a general theory on torus orbits of symplectic leaves for Poisson structures 
defined by quasitriangular $r$-matrices. Our descriptions of 
the $T$-leaves 
of these Poisson manifolds are in terms of 
{\it extended Bruhat cells, extended Richardson varieties,} and {\it extended
Double Bruhat cells associated to conjugacy classes} (see \cite{Lu-Victor:flags} for detail).
In \cite{Balazs-Lu:BS}, the induced Poisson structures on extended
Bruhat cells are computed explicitly in the so-called Bott-Samelson coordinates in terms of root strings and the
structure constants of the Lie algebra $\g$, and they are shown to give rise to 
polynomial Poisson algebras that are {\it symmetric nilpotent 
semi-quadratic Poisson-Ore extensions} of ${\mathbb C}$ in the sense of \cite[Definition 4]{Goodearl-Yakimov:PNAS}.
Symplectic groupoids of extended Bruhat cells in terms of extended double Bruhat cells will be studied in 
another paper.

We would like to point out that some examples of Poisson structures defined by quasitriangular $r$-matrices and Lie algebra actions, together with their quantizations, have been studied by several authors (see, for example, \cite{D-G-Majid:R, E-Yvette:q-quantum, EK:homog, Zwicknagl}). Our paper gives new classes of such Poisson structures, and it would be very interesting to
study their quantizations.

\subsection{Acknowledgments}
The authors wish to thank Anton Alekseev, David Li-Bland, and Xiaomeng Xu for helpful comments. This work was partially supported by a University of Hong Kong Post-graduate Studentship and by the Research Grants Council of the Hong Kong SAR, China (GRF HKU 704310 and 703712). 

\subsection{Notation}\lb{subsec-intro-notation}
Throughout the paper, vector spaces are understood to be over ${\bf k} = \Rset$ or $\Cset$. The non-degenerate bilinear pairing between a finite dimensional vector space $V$ and its dual space $V^*$ will always be denoted by $\lara$. The annihilator of a vector subspace $U$ of $V$ is, by definition, 
$U^0 =\{\xi \in V^*: \la \xi, U\ra =0\} \subset V^*$. For each integer $k \geq 1$, we identify the tensor product $V^{\otimes k}$ with the space of multi-linear maps from $V^* \times \cdots \times V^*$ ($k$-times) to ${\bf k}$, so that
\[
(v_1 \otimes v_2 \otimes \cdots \otimes v_k)(\xi_1, \xi_2, \ldots, \xi_k) = \la v_1, \xi_1\ra \la v_2, \xi\ra \cdots \la v_k, \xi_k\ra, \hs v_j \in V, \; \xi_j \in V^*.
\]
Let $\wedge^k V$ be the subspace of skew-symmetric elements in $V^{\otimes k}$, and define, for $v_1, \ldots, v_k \in V$,
\begin{equation}\lb{eq-wedge}
v_1 \wedge v_2 \wedge \cdots \wedge v_k = \sum_{\sigma \in S_k} {\rm sign}(\sigma) v_{\sigma(1)} \otimes v_{\sigma(2)} \otimes \cdots \otimes v_{\sigma(k)} \in \wedge^k V \subset V^{\otimes k}.
\end{equation}
If $\lara_{(\sV, \sW)}$ is a bilinear pairing between $V$ and another vector space $W$, for each integer
$k \geq 1$, we extend $\lara_{(\sV, \sW)}$ to a bilinear pairing between $\wedge^kV$ and $\wedge^kW$, also denoted as
$\lara_{(\sV, \sW)}$, via
\[
\la x_1 \wedge \cdots \wedge x_k, \; y_1 \wedge \cdots \wedge y_k \ra_{(\sV, \sW)} 
=\det ((\la x_i, y_j\ra_{(\sV, \sW)})_{i, j = 1, \ldots, k}), \hs x_j\in V, \, y_j \in W.
\]
In particular, we identify $\wedge^k(V^*)$ with $(\wedge^k V)^*$ via the bilinear pairing
between $\wedge^k V$ and $\wedge^k (V^*)$ which is 
the extension of the pairing $\lara$ between $V$ and $V^*$.
For any element $r = \sum_i u_i \otimes v_i \in V \otimes V$, define 
\begin{equation}\lb{eq-r-V-sharp}
r^\#: \;\;\; V^* \lrw V, \;\;\; r^\#(\xi) = \sum_i \la \xi, u_i\ra v_i, \hs \xi \in V^*.
\end{equation}
If $\Lam =\sum_i u_i \wedge v_i \in \wedge^2 V$, then $\Lam^\#: V^* \to V$ is given by $\Lam^\#(\xi) = \sum_i \la \xi, u_i\ra v_i - \la \xi, v_i\ra u_i$, or
\begin{equation}\lb{eq-Lam-sharp}
\la \Lam^\#(\xi), \; \eta\ra = \Lam (\xi,  \, \eta) = \la \Lambda, \; \xi \wedge \eta\ra, \hs \xi, \,\eta \in V^*.
\end{equation}
If $\sigma: V \to W$ is a linear map, denote also by $\sigma$ the linear map $V^{\otimes k}\to W^{\otimes k}$
given by $\sigma(v_1 \otimes \cdots \otimes v_k) = \sigma(v_1) \otimes \cdots \otimes \sigma(v_k)$.
%Define the pairing between $\wedge^k V \subset V^{\otimes k}$ and $\wedge^k V^* \subset (V^*)^{\otimes k}$ by 
%\[
%\la x_1 \wedge \cdots \wedge x_k, \; \xi_1 \wedge \cdots \wedge \xi_k \ra_k =\det ((\la x_i, \xi_j\ra_\sV)_{i, j = 1, \ldots, %k}), \hs x_j \in V, \, \xi_j \in V^*.
%\]
%A symmetric bilinear form $\lara_\sV$ on $V$ will be extended bilinearly to one on $\wedge^k V$ via
%\[
%\la x_1 \wedge \cdots \wedge x_k, \; y_1 \wedge \cdots \wedge y_k \ra_\sV =\det ((\la x_i, y_j\ra_\sV)_{i, j = 1, \ldots, k}), %\hs x_j, \,y_j \in V.
%\] 
%A subspace $U$ of $V$ is said to be {\it coisotropic} (resp. {\it isotropic, Lagrangian}) with respect to $\lara_\sV$ if $U^\perp \subset U$ (resp. $U \subset U^\perp$, $U = U^\perp$), where $U^\perp = \{v \in V : \la v, U \ra_\sV = 0 \} \subset V$.
If $V = V_1 \oplus  \cdots \oplus V_n$ is a direct sum of vector spaces, and if $A_j \in V_j^{\ot k}$, $1 \leq j \leq n$, denote by  $(0, \ldots, A_j, \ldots, 0) \in V^{\otimes k}$  the image of $A_j$ under the embedding of $V_j$ into $V$, and let
\[
(A_1, \ldots, A_n) = \sum_{j=1}^n(0, \ldots, A_j, \ldots, 0) \in V^{\ot k}.
\]
 
We will work in the category of finite dimensional smooth real or complex manifolds. If $Y$ is a real or complex manifold, for an integer $k \geq 1$, $\V^k(Y)$ denotes the space of $k$-vector fields on $Y$, i.e., smooth or holomorphic sections of $\wedge^k TY$, where $TY$ is the smooth or holomorphic tangent bundle of $Y$. One has the Schouten bracket $[ \; , \; ]: \V^k(Y) \times \V^l(Y) \to \V^{k+l-1}(Y)$ which restricts to the Lie brackets on $\V^1(Y)$ and satisfies $[A, f] = A(f)$ for $A \in \V^1(Y)$, $f \in \V^0(Y)$, and
\begin{align*}
&[A, \, B \wedge C] = [A, B]\wedge C + (-1)^{(k-1)l} B \wedge [A, C], \hs A \in \V^k(Y), \, B \in \V^l(Y), \, C \in \V^m(Y),\\
&[A, B] = -(-1)^{(k-1)(l-1)} [B, A], \hs A \in \V^k(Y), \, B \in \V^l(Y).
\end{align*}
%If $F: X \to Y$, the same letter $F$ also denotes the induced map $TX \to TY$. For $V_\sX \in \V^k(X)$ and $V_\sY \in \V^k(Y)$, %by $F(V_\sX) = V_\sY$ we mean that $F(V_\sX(x)) = V_\sY(F(x))$ for all $x \in X$, and we also say that $V_\sX$ and $V_\sY$ are %$F$-related.
A real (resp. complex) {\it Poisson manifold} is a pair $(Y, \pi_\sY)$, where $Y$ is a real (resp. complex) manifold and $\pi_\sY$ is a smooth (resp. holomorphic) bi-vector field on $Y$ such that $[\pi_\sY, \pi_\sY]=0$.
%A map $F: (Y, \pi_\sY)\to (Z, \pi_\sZ)$ between Poisson manifolds is said to be Poisson if $F(\pi_\sY(y)) = \pi_\sZ(F(y))$ for all $y \in Y$. 
Given Poisson manifolds $(Y_j, \pi_{\sY_j})$, $1 \leq j \leq n$, the direct product Poisson structure on the product manifold $Y_1 \times \cdots \times Y_n$ is $(\pi_{\sY_1}, \ldots, \pi_{\sY_n})$, also denoted by $\pi_{\sY_1} \times \cdots \times \pi_{\sY_n}$.  The $n$-fold product of a Poisson manifold $(Y, \piY)$ with itself is denoted by $(Y^n, \pi_\sY^n)$.

Let $G$ be a Lie group with Lie algebra $\g$.  
%The left and right translations on $G$ by $g \in G$ will be denoted by $l_g$ and $r_g$ respectively. For an integer $k \geq 1$ and for $x \in \wedge^k\g$ and $\xi \in \wedge^k\g^*$, $x^l$ and $x^R$ (resp. $\xi^l$ and $\xi^R$)  will respectively denote the left and right invariant $k$-vector fields (resp. $k$-forms) on $G$ with value $x$ (resp. $\xi$) at $e$.
Identifying $\wedge \g$ with the space of left invariant multivector fields on $G$, one has the Schouten bracket on $\wedge \g$ which is also denoted by $[ \; , \; ]$. In particular, for $\Lam \in \wedge^2 \g$ and using our convention in \eqref{eq-wedge}, $[\Lam, \Lam] \in \wedge^3 \g$ is given by
\begin{equation}\lb{eq-wedge-Lam}
\frac{1}{2}[\Lam, \Lam](\xi, \eta, \zeta) = \la \xi, \; [\Lam^\#(\eta), \, \Lam^\#(\zeta)]\ra + \la \eta, \; [\Lam^\#(\zeta), \, \Lam^\#(\xi)]\ra + \la \zeta, \; [\Lam^\#(\xi), \, \Lam^\#(\eta)]\ra, 
\end{equation}
for $\xi, \eta, \zeta \in \g^*$. A {\it left action} of $\g$ on a manifold $Y$ is a Lie algebra anti-homomorphism $\lam: \g \to \V^1(Y)$, while a {\it right action} of $\g$ on $Y$ is a Lie algebra homomorphism $\rho: \g \to \V^1(Y)$. If $\lam: G \times Y \to Y, (g, y) \mapsto gy$ is a left action of $G$ on $Y$, one has the induced left action of $\g$ on $Y$, also denoted by $\lam$, given by 
\[
\lam:\; \g \longrightarrow \V^1(Y), \;\; \lam(x)(y) = \frac{d}{dt}|_{t=0} \exp(tx)y, \hs x \in \g, \, y \in Y.
\]
Similarly, a right Lie group action $\rho: Y \times G \to Y, (y, g) \to yg$ induces a right Lie algebra action 
\[
\rho:\; \g \longrightarrow \V^1(Y), \;\;\rho(x)(y) = \frac{d}{dt}|_{t=0} y\exp(tx), \hs x \in \g, \, y \in Y.
\]

\sectionnew{Review on Poisson Lie groups, Lie bialgebras,  and Poisson actions}\lb{sec-review}
In this section, we set up our conventions (on signs and constants) and review some facts on Lie bialgebras that will be used in this paper and in subsequent papers. Except for that in $\S$\ref{subsec-Poi-r}, all materials are standard and well-known, and we refer to \cite{chari-pressley, dr:Poi-Lie, dr:quantum, etingof-schiffmann,  Yvette:lectures, k-s:quantum,  lu-we:Poi, RT:factorizable, STS1, STS2, STS:RIMS} for more details. All Lie bialgebras in this paper are assumed to be finite dimensional.

\subsection{Lie bialgebras}\lb{subsec-Lie-bialgebras} 
A {\it Lie bialgebra} is a pair $(\g, \delta_\g)$ where $\g$ is a Lie algebra and $\delta_\g: \g \to \wedge^2 \g$ a linear map, called the {\it co-bracket}, such that 1)  the dual map $\delta_\g^*: \wedge^2\g^* \to \g^*$ of $\delta_\g$ is a Lie bracket on the dual space $\g^*$ of $\g$, and 2) $\delta_\g$ satisfies the $1$-cocycle condition
\begin{equation}\lb{eq-cocycle-def}
\delta_\g[x, y] = [x, \delta_\g(y)] + [\delta_\g(x), y], \hs x, y \in \g.
\end{equation}
A {\it Lie bialgebra homomorphism} from one Lie bialgebra $(\g_1, \delta_{\g_1})$ to another  $(\g_2, \delta_{\g_2})$ is a Lie algebra homomorphism $\phi: \g_1 \to \g_2$ such that $\phi(\delta_{\g_1}(x)) = \delta_{\g_2}(\phi(x))$ for all $x \in \g_1$.

If $(\g, \delta_\g)$ is a Lie bialgebra, so is $(\g^*, \delta_{\g^*})$, where the Lie bracket on $\g^*$ is given by $\delta_\g^*: \wedge^2\g^* \to \g^*$, the dual map of $\delta_\g$, and $\delta_{\g^*}: \g^* \to \wedge^2 \g^*$ is the dual map of the Lie bracket on $\g$. The Lie bialgebra $(\g^*, \delta_{\g^*})$ is called the  {\it dual Lie bialgebra} of $(\g, \delta_\g)$. If $(\g', \delta_{\g'})$ is a Lie bialgebra isomorphic to the dual Lie bialgebra $(\g^*, \delta_{\g^*})$, we will call $(\g, \delta_\g)$ and $(\g', \delta_{\g'})$ a {\it dual pair of Lie bialgebras}. 

A {\it sub-Lie bialgebra} of a Lie bialgebra $(\g, \delta_\g)$ is a Lie bialgebra $(\p, \delta_\p)$ with an injective Lie bialgebra homomorphism $(\p, \delta_\p) \to (\g, \delta_\g)$. A {\it quotient Lie bialgebra} of $(\g, \delta_\g)$ is a Lie bialgebra $(\q, \delta_\q)$ with a surjective Lie bialgebra homomorphism $(\g, \delta_\g) \to (\q, \delta_\q)$.

\subsection{The Double Lie (bi)algebra of a Lie bialgebra}\lb{subsec-double}
Recall that a quadratic Lie algebra is a pair $(\d, \lara_\d)$, where $\d$ is a Lie algebra and $\lara_\d$ is a non-degenerate symmetric ad-invariant bilinear form on $\d$. Let $(\g, \delta_\g)$ be a Lie bialgebra and let $(\g^*, \delta_{\g^*})$ be its dual Lie bialgebra. Consider the direct sum vector space $\d = \g \oplus \g^*$ and denote its elements by $x+\xi$, where $x \in \g$ and $\xi \in \g^*$. Then $\d$ carries the symmetric bilinear form 
\[
\la x + \xi, \; y + \eta \ra_\d = \la x, \eta\ra + \la \xi, y\ra, \hs x, y \in \g, \, \xi, \eta \in \g^*,
\]
and a unique Lie bracket $[\, , \, ]$ extending those on $\g$ and $\g^*$ such that $\lara_\d$ is ad-invariant. Namely
\begin{equation}\lb{eq-bra-d}
[x+\xi, \; y + \eta] = [x, \; y] + \ad_\xi^* y - \ad_\eta^* x + [\xi, \; \eta] + \ad_x^* \eta - \ad_y^* \xi, 
\hs x, y \in \g, \, \xi, \eta \in \g^*,
\end{equation}
where $\ad^*$ denotes both the coadjoint action of $\g$ on $\g^*$ and of $\g^*$ on $\g$, i.e.,
\[
\la \ad_x^*\xi, \; y\ra = \la \xi, \; [y, x]\ra, \hs \la \ad_\xi^* x, \; \eta\ra = \la x, \; [\eta, \xi]\ra,
\hs x, y \in \g, \; \xi \in \g^*.
\]
The quadratic Lie algebra $(\d, \lara_\d)$ is called the {\it (Drinfeld) double Lie algebra} of the Lie bialgebra $(\g, \delta_\g)$.

\bde{de-lag-splitting}
Assume that $(\d, \lara_\d)$ is an even dimensional quadratic Lie algebra. By a {\it Lagrangian subalgebra} of $(\d, \lara_\d)$ we mean a Lie subalgebra $\g$ of $\d$ which is also a Lagrangian, i.e., maximal isotropic, with respect to $\lara_\d$. A decomposition $\d = \g + \g'$ of vector spaces, where both $\g$ and $\g'$ are Lagrangian subalgebras of $(\d, \lara_\d)$, is called a {\it Lagrangian splitting} of $(\d, \lara_\d)$, and in this case  $((\d, \lara_\d), \g, \g')$ is also called a {\it Manin triple}. Given a Lagrangian splitting $\d = \g + \g'$, by identifying $\g$ and $\g'$ as the dual spaces of each other using $\lara_\d$, one obtains $\delta_\g: \g \to \wedge^2 \g$ and $\delta_{\g^\prime}: \g' \to \wedge^2 \g'$, the dual maps of the Lie brackets on $\g'$ and on $\g$ respectively, and one has the dual pair of Lie bialgebras $(\g, \delta_\g)$ and $(\g', \delta_{\g'})$. Their double Lie algebra is clearly isomorphic to $(\d, \lara_\d)$. The notion of Lie bialgebras is therefore equivalent to that of even-dimensional quadratic Lie algebras with Lagrangian splittings.
\ede

Let $(\g, \delta_\g)$ be a Lie bialgebra, with dual Lie bialgebra $(\g^*, \delta_{\g^*})$ and double Lie algebra $(\d, \lara_\d)$. Consider the direct product Lie algebra  $\d^2 = \d \oplus \d$ and denote its elements by $(a, b)$, where $a, b \in \d$. Let $\lara_{\d^2}$ be the bilinear form on $\d^2$ given by
\begin{equation}\lb{eq-lara-d2-0}
\la (a_1, a_2), \,(a_1^\prime, a_2^\prime)\ra_{\d^2} = \la a_1, \,a_1^\prime \ra_\d - \la a_2, \,a_2^\prime\ra_\d, \hs a_1, 
a_2, \, a_1^\prime, \, a_2^\prime \in \d.
\end{equation}
The quadratic Lie algebra $(\d^2, \lara_{\d^2})$ has the Lagrangian splitting $\d^2 = \d_{\rm diag} + \d'$, where
\[
\d_{\rm diag}  = \{(a,\,a): a \in \d\} \hs \mbox{and} \hs \d' = \g^* \oplus \g = \{(\xi, x): \xi \in \g^*, x \in \g\}.
\]
Identifying $\d \cong \d_{\rm diag}$ via $a \mapsto (a, a)$, one obtains the dual pair of Lie bialgebras $(\d, \delta_\d)$ and $(\d', \delta_{\d'})$. The Lie bialgebra $(\d, \delta_\d)$ is called the {\it (Drinfeld) double Lie bialgebra} of the Lie bialgebra $(\g, \delta_\g)$. Note that the non-degenerate pairing between $\d \cong \d_{\rm diag}$ and $\d'$ induced by $\lara_{\d^2}$ is given by
\begin{equation}\lb{eq-pairing-dd}
\la a, \; (\xi, x)\ra_{(\d, \d^\prime)} = \la \xi - x, \; a \ra_\d, \hs  a \in \d, \;  (\xi, x) \in \d'.
\end{equation}
It is straightforward to check (see also \cite[$\S$4.2]{etingof-schiffmann} and \cite[$\S$1.1]{Yakimov-leaves}) that
\begin{equation}\lb{eq-delta-d-g-g}
\delta_\d|_\g = \delta_\g: \;\; \g \lrw \wedge^2 \g \hs \mbox{and} \hs \delta_\d|_{\g^*} = -\delta_{\g^*}: \;\;\g^* \lrw \wedge^2 \g^*.
\end{equation}
In particular, both $(\g, \delta_\g)$ and $(\g^*, -\delta_{\g^*})$ are sub-Lie bialgebras of $(\d, \delta_\d)$.

\bex{ex-double-of-sub}
Let $(\g, \delta_\g)$ be a Lie bialgebra and $(\p, \delta_\p)$ a sub-Lie bialgebra of $(\g, \delta_\g)$. The double Lie bialgebra $(\d_\p, \, \delta_{\d_\p})$ of $(\p, \delta_\p)$ can be identified with a ``sub-quotient" of the double Lie bialgebra $(\d, \delta_\d)$ of $(\g, \delta_\g)$: indeed, $(\p + \g^*, \;\delta_\d|_{(\p +\g^*)})$ is a sub-Lie bialgebra of $(\d,\delta_\d)$, and the map $\p + \g^* \to \d_\p: x + \xi \mapsto x + \xi|_\p$, $x \in \p, \xi \in \g^*$, is a surjective Lie bialgebra homomorphism from $(\p + \g^*, \;\delta_\d|_{(\p +\g^*)})$ to $(\d_\p, \, \delta_{\d_\p})$. \hfill $\diamond$
\eex

\subsection{Quasitriangular $r$-matrices}\lb{subsec-r-matrices}
Let $\g$ be any Lie algebra. For $r =\sum_i x_i \ot y_i \in \g \ot \g$, define 
\begin{equation}\lb{eq-delta-r}
\delta_r: \;\; \g \lrw \g \otimes \g, \;\;\; \delta_r(x) = \ad_x r = \sum_i ([x, x_i] \ot y_i + x_i \ot [x, y_i]), \hs x \in \g.
\end{equation}
An element $r \in \g \ot \g$ is called a {\it classical $r$-matrix (or simply an $r$-matrix) on $\g$} if $(\g, \delta_r)$ is a Lie bialgebra. Writing $r$ as $r = \Lam + s$ with $\Lam \in \wedge^2 \g$ and $s \in S^2 \g$, it follows from the definition that $r$ is a classical $r$-matrix on $\g$ if and only if $s \in (S^2 \g)^\g$ and $[\Lam, \Lam] \in (\wedge^3 \g)^\g$, where $(S^2\g)^\g$ and $(\wedge^3 \g)^\g$ denote respectively the ad-invariant subspaces of $S^2 \g$ and $\wedge^3 \g$, and $[ \; , \; ]$ is the Schouten bracket on $\wedge \g$ (see $\S$\ref{subsec-intro-notation} and \eqref{eq-wedge-Lam}). In this paper, we will refer to classical $r$-matrices simply as $r$-matrices. Given a Lie bialgebra $(\g, \delta_\g)$, an element $r \in \g \otimes \g$ such that $\delta_\g = \delta_r$ is called an {\it $r$-matrix for $(\g, \delta_\g)$}. A Lie bialgebra $(\g, \delta_\g)$ is said to be {\it co-boundary} if it has an $r$-matrix.

Recall \cite{dr:quantum} the classical Yang-Baxter operator ${\rm CYB}: \g \otimes \g \to \g \ot \g \ot \g$ given by  
\begin{equation}\lb{eq-CYB-r}
{\rm CYB}(r) = [r_{12}, \; r_{13}] + [r_{12}, \; r_{23}] + [r_{13}, \; r_{23}] \in \g \otimes \g \otimes \g \subset
U(\g)^{\otimes 3}, \hs r \in \g \ot \g,
\end{equation}
where $U(\g)$ is the universal enveloping algebra of $\g$, $[ \; , \; ]$ on the right hand side denotes the commutator bracket in $U(\g)^{\otimes 3}$, and $r_{12}, r_{13}, r_{23} \in U(\g)^{\otimes 3}$ are respectively given by
\[
r_{12} = \sum_i x_i \otimes y_i \otimes 1, \hs r_{23} = \sum_i 1 \otimes x_i \otimes y_i, \hs 
r_{13} = \sum_i x_i \otimes 1 \otimes y_i
\]
if $r = \sum_i x_i \otimes y_i \in \g \otimes \g$. One checks directly \cite{dr:Poi-Lie, Yvette:lectures} (note our convention in $\S$\ref{subsec-intro-notation} and \eqref{eq-r-V-sharp}) that 
\begin{align}\lb{eq-CYB-Lam}
{\rm CYB}(\Lam) &= \frac{1}{2} [\Lam, \; \Lam] \in \wedge^3 \g\; \mbox{if} \; \Lam \in \wedge^2 \g;\\
\lb{eq-CYB-s}
{\rm CYB}(s) &\in (\wedge^3 \g)^\g \; \mbox{and} \; {\rm CYB}(s)(\xi, \eta, \zeta) = 
\la \xi, \; [s^\#(\eta), \, s^\#(\zeta)]\ra, \; \xi, \eta, \zeta \in \g^*, \; \mbox{if} \; s \in (S^2 \g)^\g;\\
\lb{eq-CYB-Lam-s}
{\rm CYB}(r) &= {\rm CYB}(\Lam) + {\rm CYB}(s) \in \wedge^3 \g  \; \mbox{if} \; r  = \Lam + s \; \mbox{with} \; 
\Lam\in \wedge^2 \g \; \mbox{and} \; s \in (S^2\g)^\g.
\end{align}
An element $r \in \g \otimes \g$ is called a {\it quasitriangular $r$-matrix} on $\g$ if the symmetric part of $r$ is ad-invariant and if $r$ satisfies the {\it Classical Yang Baxter Equation} ${\rm CYB}(r) = 0$. By \eqref{eq-CYB-Lam} - \eqref{eq-CYB-Lam-s}, a quasitriangular $r$-matrix on $\g$ is indeed an $r$-matrix on $\g$ (see also \cite{dr:quantum, Yvette:lectures, Majid:quantum}). 
Let
\[
r^{21} = \sum_i y_i \otimes x_i \hs \mbox{if} \hs r = \sum_i x_i \otimes y_i,
\]
It is clear that $r \in \g \ot \g$ is a quasitriangular $r$-matrix on $\g$ if and only if $r^{21}$ is. 

\bde{de-quasi-r}
1) A quasitriangular $r$-matrix for a  Lie bialgebra $(\g, \delta_\g)$ is also called a {\it quasitriangular structure}
of $(\g, \delta_\g)$. If $r$ is a quasitriangular $r$-matrix on $\g$, we call the pair $(\g, r)$ a 
{\it quasitriangular Lie bialgebra} (with $\delta_r: \g \to \wedge^2\g$ as its co-bracket).
A quasitriangular $r$-matrix $r$ on $\g$ is said to be {\it factorizable} if its symmetric part is non-degenerate, and in this case we call $(\g, r)$ a factorizable Lie bialgebra;

2) Given two quasitriangular Lie bialgebras $(\g_1, r_1)$ and $(\g_2, r_2)$, we say 
$\phi: (\g_1, r_1) \to (\g_2, r_2)$ is a {\it morphism of quasitriangular Lie bialgebras} if 
$\phi: \g_1 \to \g_2$ is a Lie algebra homomorphism and $\phi(r_1) = r_2$.  
\ede

\bre{re-phi-r} A quasitriangular Lie bialgebra is thus the Lie bialgebra together with a
quasitriangular quasitriangular structure.
A morphism $\phi: (\g_1, r_1) \to (\g_2, r_2)$ of quasitriangular Lie bialgebras induces a morphism
$\phi: (\g_1, \delta_{r_1}) \to (\g_2, \delta_{r_2})$ of Lie bialgebras, but the converse fails in general (for example,
when both $\g_1$ and $\g_2$ are abelian Lie algebras).
\hfill $\diamond$
\ere

Let $(\g, \delta_\g)$ be a Lie bialgebra and $(\d, \lara_\d)$ its double Lie algebra.  Let $\dim \g = m$ and let ${\rm Gr}_0(m, \d)$ be the set of all $m$-dimensional subspaces of $\d$ complementary to $\g$. Then
\begin{equation}\lb{eq-r-kr}
\g \ot \g \lrw {\rm Gr}_0(m, \d), \;\;\; r \longmapsto \k_r = \{-r^\#(\xi) + \xi: \;  \xi \in \g^*\}
\end{equation}
is a bijection, where $r^\#: \g^* \to \g$ is given in \eqref{eq-r-V-sharp}. Note that $\k_r^\perp = \k_{-r^{21}}$ for any $r \in \g \ot \g$, where $\q^\perp = \{a \in \d: \la a, \q\ra_\d = 0\}$ for $\q \subset \d$. In particular, $\k_r$ is a Lagrangian subspace of $\d$ with respect to $\lara_\d$ if and only if $r \in \wedge^2 \g$. A proof of the following result of Drinfeld can be found in \cite{A-M:linearization}.

\ble{le-drinfi-r-0} {\rm (Drinfeld \cite{dr:quasi})}
An element $r \in \g \ot \g$ is an $r$-matrix for $(\g, \delta_\g)$ if and only if $[\g, \k_r] \subset \k_r$, and $r$ is a  quasitriangular $r$-matrix for $(\g, \delta_\g)$ if and only if $\k_r$ is a Lie ideal of $\d$.
\ele

Applying \leref{le-drinfi-r-0} to the double Lie bialgebra $(\d, \delta_\d)$ of any Lie bialgebra $(\g, \delta_\g)$, one sees that $(\d, \delta_\d)$ is co-boundary: the two Lie ideals $0 \oplus \d$ and $\d \oplus 0$ of $\d^2 = \d \oplus \d$ give rise to the two factorizable quasitriangular matrices $r_\d$ and $-(r_\d)^{21}$ for $(\d, \delta_\d)$ respectively, with
\begin{equation}\lb{eq-double-r}
r_\d = \sum_{i=1}^m x_i \otimes \xi_i  \in \d \otimes \d ,
\end{equation}
where $\{x_i\}_{i=1}^m$ is any basis of $\g$ and $\{\xi_i\}_{i=1}^m$ its dual basis of $\g^*$. It is straightforward to check \cite[$\S$4.2]{etingof-schiffmann} that $\delta_\d = \delta_{r_\d}: \d \to \wedge^2 \d$.

\bde{de-r-on-d}
We call $r_\d \in \d \ot \d$ in \eqref{eq-double-r}  {\it the quasitriangular $r$-matrix on $\d$ associated to the Lagrangian splitting $\d = \g + \g^*$ of $(\d, \lara_\d)$}. Its skew-symmetric part
\begin{equation}\lb{eq-Lambda-d}
\Lambda_{\g, \g^*}= \frac{1}{2} \sum_{i=1}^m x_i \wedge \xi_i \in \wedge^2 \d 
\end{equation}
will be called the {\it skew-symmetric $r$-matrix} on $\d$ associated to the Lagrangian splitting $\d = \g + \g^*$.
\ede

Let $(\g, \delta_\g)$ be a Lie bialgebra with a quasitriangular structure $r \in \g \ot \g$, and let $(\g^*, \delta_{\g^*})$ 
and $(\d, \delta_\d)$ be respectively the dual and the double Lie bialgebra of $(\g, \delta_\g)$. Define (see \eqref{eq-r-V-sharp})
\begin{align}\lb{eq-r-pm-0}
&r_+ = r^\#: \; \g^* \lrw \g, \hs   r_- = -(r^{21})^\# =-r_+^*: \; \g^* \lrw \g,\\
\lb{eq-ppm}
&p_\pm: \;\d \lrw \g, \;\;\; p_+(x + \xi) = x + r_+(\xi), \hs p_-(x+\xi) = x + r_-(\xi), \hs x \in \g, \; \xi \in \g^*.
\end{align}
It follows from \leref{le-drinfi-r-0} (see also \cite{dr:almost, etingof-schiffmann, Hodges-Yakimov, Yvette:lectures, RT:factorizable}) that both $p_+$ and $p_-$ 
are Lie algebra homomorphisms.
With $r_\d \in \d \otimes \d$ given in \eqref{eq-double-r}, one has 
\begin{equation}\lb{eq-r-ppm}
p_+(r_\d) = \sum_{i=1}^m x_i \otimes r_+(\xi_i) = r, \hs p_-(r_\d) = \sum_{i=1}^m x_i \otimes r_-(\xi_i) = -r^{21}.
\end{equation}
In particular, $p_+$ and $p_-$ are also Lie bialgebra homomorphisms from $(\d, \delta_\d)$ to $(\g, \delta_\g)$,
and it follows that $r_+, r_-: (\g^*, -\delta_{\g^*}) \to (\g, \delta_\g)$ are also Lie bialgebra homomorphisms. Let 
\begin{equation}\lb{eq-de-l-pm}
\f_- = {\rm Im}(r_-) \subset \g, \hs \f_+ = {\rm Im} (r_+) \subset \g.
\end{equation}
Then $(\f_-, \delta_\g|_{\f_-})$ and $(\f_+, \delta_\g|_{\f_+})$ are sub-Lie bialgebras of $(\g, \delta_\g)$. Consider the well-defined non-degenerate bilinear pairing $\lara_{(\f_-, \f_+)}$ between $\f_-$ and $\f_+$ given by
\begin{equation}\lb{eq-pairing-ll}
\la r_-(\xi), \; r_+(\eta)\ra_{(\f_-, \f_+)} = -\la \xi, \; r_+(\eta)\ra =\la r_-(\xi), \; \eta\ra, \hs \xi,\, \eta \in \g^*.
\end{equation}

\ble{le-lpm-dual} (See also \cite[Proposition 4.1]{etingof-schiffmann})
Under the non-degenerate pairing $\lara_{(\f_-, \f_+)}$ between $\f_-$ and $\f_+$, $(\f_-, \delta_\g|_{\f_-})$ and $(\f_+, -\delta_\g|_{\f_+})$ form a dual pair of Lie bialgebras. 
\ele

\begin{proof}
Let $x_-, y_- \in \f_-$ and $x_+, y_+ \in \f_+$. Let $x_- = r_-(\xi)$ and $x_+ = r_+(\eta)$, where $\xi, \eta \in \g^*$. As $r_{\pm}: (\g^*, -\delta_{\g^*}) \to (\g, \delta_\g)$ are Lie bialgebra homomorphisms, one has
\begin{align*}
\la x_- \wedge y_-, \; \delta_\g(x_+)\ra_{(\f_-, \f_+)} 
& = \la x_- \wedge y_-, \; \delta_\g(r_+(\eta))\ra_{(\f_-, \f_+)} = 
-\la x_- \wedge y_-, \; r_+(\delta_{\g^*}(\eta))\ra_{(\f_-, \f_+)}\\
& = -\la x_- \wedge y_-, \; \delta_{\g^*}(\eta)\ra = -\la [x_-, \, y_-], \; \eta\ra 
= -\la [x_-, \, y_-], \; x_+\ra_{(\f_-, \f_+)},\\
\la \delta_\g(x_-), \; x_+ \wedge y_+\ra_{(\f_-, \f_+)} & = \la \delta_\g(r_-(\xi)), \; x_+ \wedge y_+\ra_{(\f_-, \f_+)}
= -\la r_-(\delta_{\g^*}(\xi)), \; x_+ \wedge y_+ \ra_{(\f_-, \f_+)} \\
& = -\la \delta_{\g^*}(\xi), \; x_+ \wedge y_+ \ra = -\la \xi, \; [x_+, y_+]\ra = \la x_-, \; [x_+, y_+]\ra_{(\f_-, \f_+)}.
\end{align*}
\end{proof}

\bre{re-r-fpm}
Let $\{x_i\}_{i=1}^m$ be a basis of $\g$ such that $\{x_i\}_{i=1}^l$ is a basis of $\f_-$, and let $\{\xi_i\}_{i=1}^m$ be the dual basis of $\g^*$. As ${\rm Span}\{\xi_{l+1}, \ldots, \xi_m\} = \f_-^0 = \ker r_+$, $\{r_+(\xi_i)\}_{i=1}^l$ is a basis of $\f_+$, dual to the basis $\{x_i\}_{i=1}^l$ of $\f_-$ under the pairing $\lara_{(\f_-, \f_+)}$. Consequently,
\begin{equation}\lb{eq-r-fpm}
r = \sum_{i=1}^m x_i \otimes r_+(\xi_i) = \sum_{i=1}^l x_i \otimes r_+(\xi_i) \in \f_- \otimes \f_+ \subset \g \otimes \g.
\end{equation}
\hfill$\diamond$
\ere

Let $(\d_{\f_-}, \delta_{\d_{\f_-}})$ be the double Lie bialgebra of $(\f_-, \delta_\g|_{\f_-})$, and let $r_{\d_{\f_-}}$ be the quasitriangular $r$-matrix on $\d_{\f_-}$ associated to the Lagrangian splitting $\d_{\f_-} = \f_- + \f_+$. Consider the linear map
\begin{equation}\lb{eq-q}
q: \;\; \d_{\f_-} \lrw \g, \;\;\; q(x_-, x_+) = x_- + x_+, \hs  x_- \in \f_-, \; x_+ \in \f_+.
\end{equation}

\ble{le-rr-quasi} {\rm (\cite[Proposition 4.1]{etingof-schiffmann})} The map
$q: (\d_{\f_-}, \delta_{\d_{\f_-}}) \to (\g, \delta_\g)$ is a Lie bialgebra homomorphism. In fact, $q(r_{\d_{\f_-}}) = r$. 
\ele

\begin{proof}
Identify $(\d_{\f_-}, \delta_{\d_{\f_-}})$ with the quotient of the sub-Lie bialgebra $(\f_- + \g^*, \, \delta_\d|_{\f_- + \g^*})$ of $(\d, \delta_\d)$ as in \exref{ex-double-of-sub}. Then $q$, being the map induced by the Lie algebra homomorphism $p_+: \d \to \g$, is a Lie algebra homomorphism. It follows from the definition of $r_{\d_{\f_-}}$ and \eqref{eq-r-fpm} that
$q(r_{\d_{\f_-}}) = r$.
\end{proof}

\subsection{Twists of Lie bialgebras and of quasitriangular $r$-matrices}\lb{subsec-twists}
Let $(\g, \delta_\g)$ be a Lie bialgebra and $(\d, \lara_\d)$ its double Lie algebra. An element $t \in \wedge^2 \g$ is called a {\it twisting element} for $(\g, \delta_\g)$ if the Lagrangian subspace  $\k_t = \{-t^\#(\xi) + \xi: \; \xi \in \g^*\}$ of $\d$ (see \eqref{eq-r-kr}) is a Lie subalgebra of $\d$. Using \eqref{eq-wedge-Lam}, one readily checks \cite{dr:quasi}
that $t \in \wedge^2 \g$ is a twisting element for $(\g, \delta_\g)$ if and only if
\begin{equation}\lb{eq-delta-12-C}
\delta_{\g}(t) +\frac{1}{2} [t, \,t] = 0,
\end{equation}
where the co-bracket $\delta_\g: \g \to \wedge^2 \g$ is linearly extended to $\delta_\g: \wedge^2 \g \to \wedge^3 \g$ via
\begin{equation}\lb{eq-delta-g-wedge-3}
\delta_{\g}(x \wedge y) = \delta_{\g}(x) \wedge y - x \wedge \delta_{\g}(y), \hs x, \, y \in \g.
\end{equation}
A twisting element $t$ for $(\g, \delta_\g)$ gives rise to the Lagrangian splitting $\d = \g + \k_{t}$ and thus another co-bracket $\delta_{\g, t}: \g \to \wedge^2 \g$ on $\g$ given by 
\begin{equation}\lb{eq-g1g1}
\delta_{\g, t}(x) = \delta_\g(x) - [x, \, t] = \delta_\g(x) + [t, \, x], \hs x \in \g.
\end{equation}
The Lie bialgebra $(\g, \delta_{\g, t})$ is called the {\it twist} of $(\g, \delta_\g)$ by the twisting element $t$, and we also say that the Lagrangian splitting $\d = \g + \k_{t}$ is the {\it twist} of the Lagrangian splitting $\d = \g + \g^*$ by $t$. Note that the two quasitriangular $r$-matrices $r_\d$ and $r_{\d, t}$ on $\d$ associated to the two Lagrangian splittings $\d = \g + \g^*$ and $\d = \g + \k_{t}$ are now related by $r_{\d, t} = r_\d - t$, where $t \in \wedge^2 \g$ is regarded as an element in $\wedge^2 \d$ via the embedding $\g \hookrightarrow \d$. 

\ble{le-t-r}
Assume that the Lie bialgebra $(\g, \delta_\g)$ has a quasitriangular structure $r \in \g \otimes \g$, 
Then
\[
\wedge^2 \g \lrw \g \ot \g, \;\;\; t \longmapsto r - t,
\]
is a one-to-one correspondence between twisting elements for the Lie bialgebra $(\g, \delta_\g)$ and quasitriangular $r$-matrices on $\g$ that have the same symmetric part as $r$. 
\ele

\begin{proof}
Let $r = \Lam + s$, where $\Lam \in \wedge^2 \g$. Let $t \in \wedge^2 \g$. By \eqref{eq-CYB-Lam} - \eqref{eq-CYB-Lam-s},
\[
{\rm CYB}(r-t) =  \frac{1}{2}[\Lam - t, \; \Lam - t] + {\rm CYB}(s) =-[t, \Lam] + \frac{1}{2}[t, t] = \delta_\g(t) +\frac{1}{2}[t,t].
\]
Thus $r - t$ is a quasitriangular matrix on $\g$ if and only if $t$ is a twisting element for $(\g, \delta_\g)$.
\end{proof}

\bde{de-twist-r}
Given a quasitriangular $r$-matrix $r$ on a Lie algebra $\g$, any quasitriangular $r$-matrix $r'$ on $\g$ with the same 
symmetric part as $r$ will be called a {\it twist} of $r$.
\ede

\subsection{Poisson Lie groups and Poisson actions} \lb{subsec-poisson-lie-gps}
A Poisson bi-vector field $\piG$ on a Lie group $G$ is said to be {\it multiplicative} if the group multiplication map $(G \times G, \, \piG \times \piG) \to (G, \piG), (g_1, g_2) \mapsto g_1g_2$, is Poisson. A {\it Poisson Lie group} is a pair $(G, \piG)$, where $G$ is a Lie group and $\piG$ is a multiplicative Poisson bi-vector field on $G$. A Lie subgroup of $G$ that is also a Poisson submanifold with respect to $\piG$ is called a {\it Poisson Lie subgroup} of $(G, \piG)$.

Let $(G, \piG)$ be a Poisson Lie group and $\g$ the Lie algebra of $G$. Then $\piG(e) = 0$, where $e \in G$ is the identity element. Let $d_e\piG: \g \to \wedge^2 \g$ be the linearization of $\piG$ at $e$, given by $(d_e\piG)(x) = [\tilde{x}, \piG](e)$, where for $x \in \g$, $\tilde{x}$ is any local vector field on $G$ with $\tilde{x}(e) = x$. Then the pair $(\g, \delta_\g = d_e\piG)$ is a Lie bialgebra, called the {\it  Lie bialgebra of the Poisson Lie group $(G, \piG)$}. 

Let $(G, \piG)$ be a Poisson Lie group with Lie bialgebra $(\g, \delta_\g)$. Any Poisson Lie group $(G^*, \piGs)$ whose Lie bialgebra is isomorphic to the dual Lie bialgebra of $(\g, \delta_\g)$ is called a {\it dual Poisson Lie group} of $(G, \piG)$, and we also call $(G, \piG)$ and $(G^*, \piGs)$  a {\it dual pair of Poisson Lie groups}. If {\bf CHECK} $(D, \piD)$ is a 
connected Poisson Lie group whose Lie bialgebra $(\d, \delta_\d)$ is isomorphic to the double Lie bialgebra of
$(\g, \delta_\g)$, and if the Lie algebra embedding $\g \to \d$ integrates to a Lie group homomorphism $G \to D$, we
call $(D, \piD)$ a {\it Drinfeld double} of $(G, \piG)$.

A {\it left Poisson action of a Poisson Lie group} $(G, \piG)$ on a Poisson manifold $(Y, \pi_\sY)$ is  a left group action $\sigma: G \times Y \to Y$ such that $\sigma: (G \times Y, \, \piG \times \piY) \to (Y, \piY)$ is a Poisson map. Similarly, a right group action of $G$ on $Y$ given by $\sigma: Y \times G \to Y$ is a {\it right Poisson action} of $(G, \piG)$ on $(Y, \pi_\sY)$ if $\sigma: (Y\times G, \pi_\sY \times \piG) \to (Y, \pi_\sY)$ is a Poisson map.

A {\it left (resp. right) Poisson action of a Lie bialgebra} $(\g, \delta_\g)$ on a Poisson manifold $(Y, \pi_\sY)$ is a Lie algebra anti-homomorphism (resp. Lie algebra homomorphism) $\sigma: \g \to \V^1(Y)$ such that
\[
[\sigma(x), \pi_\sY] = \sigma(\delta_\g(x)), \hs \forall x \in \g.
\]
Here recall our convention from $\S$\ref{subsec-intro-notation} that 
for any integer $k \geq 2$, $\sigma$ also denotes the linear map 
\begin{equation}\lb{eq-sigma-ext}
\wedge^k \g \lrw \V^k(Y), \;\;\;
\sigma(x_1 \wedge \cdots \wedge x_k) = \sigma(x_1) \wedge \cdots \wedge \sigma(x_k), \hs x_1, \ldots, x_k \in \g.
\end{equation}

\bde{de-Poisson-space} 1) For a Poisson Lie group $(G, \piG)$, a {\it left (or right) $(G, \piG)$-Poisson space} 
is a triple $(Y, \piY, \sigma)$, where $(Y, \piY)$ is a Poisson manifold and $\sigma$ a left (or right) Poisson action of $(G, \piG)$ on $(Y, \piY)$; 

2) For a Lie bialgebra $(\g, \delta_\g)$, a {\it left (or right) $(\g, \delta_\g)$-Poisson space} is a triple $(Y, \piY, \sigma)$, where $(Y, \piY)$ is a Poisson manifold and $\sigma: \g \to \V^1(Y)$ a left (or right) Poisson action of $(\g, \delta_\g)$ on $(Y, \piY)$. 

3) For a quasitriangular Lie algebra $(\g, r)$ (see \deref{de-quasi-r}),
a left (or right) $(\g, r)$-Poisson space is a left (or right)-Poisson space of the Lie bialgebra
$(\g, \delta_r)$, where $\delta_r: \g \to \wedge^2 \g$ is given in \eqref{eq-delta-r}.
\ede

Let $(G, \pi)$ be a connected Poisson Lie group with  Lie bialgebra $(\g, \delta_\g)$, let $(Y, \pi_\sY)$ be a Poisson manifold, and let $\sigma: G \times Y \to Y$ (resp. $\sigma: Y \times G \to Y$) be a left (resp. right) action of $G$ on $Y$. Recall from $\S$\ref{subsec-intro-notation} our convention on Lie algebra actions induced from Lie group actions.

\ble{le-poi-action-infini}\cite{Weinstein-dressing} 
$(Y, \piY, \sigma)$ is a $(G, \piG)$-Poisson space if and only if it is a $(\g, \delta_\g)$-Poisson space.
\ele

\bre{re-quotient-action-0}
By \leref{le-poi-action-infini}, if $(Y, \piY, \lam)$ is a left $(\g, \delta_\g)$-Poisson space, then $(Y, \piY, -\lam)$ is a right $(\g, -\delta_\g)$-Poisson space, and $(Y, -\piY, -\lam)$ is a right $(\g, \delta_\g)$-Poisson space. \hfill $\diamond$
\ere

%\bre{re-reduction-1} {\bf (Reduction to sub-Lie bialgebras, I.)} Let $(\g, \delta_\g)$ be a Lie bialgebra and let
%$(Y, \piY, \sigma)$ be a (left or right) $(\g, \delta_\g)$-Poisson space. 
%Let $\p$ be a Lie subalgebra of $\g$ and let $\sigma|_\p: \p \to \V^1(Y)$ be the restriction of $\sigma$ to $\p$.
%If $\p$ is a co-ideal of $(\g, \delta_\g)$, then $(Y, \piY, \sigma|_\p)$ is a $(\p, \delta_\p:=\delta_\g|_\p)$-Poisson space.
%If $\p$ is an ideal of $(\g, \delta_\g)$ and if $\sigma(x) = 0$ for all $x \in \p$, then $\sigma$ induces a
%well-defined Lie algebra action $\sigma': \g/\h \to \V^1(Y)$ of the quotient Lie algebra $\g/\h$ on $Y$ by
%$\sigma'(x+\p) = \sigma(x)$ for $x \in \g$, and $(Y, \piY, \sigma')$ is a $(\g/\h, \delta_{\g/\h})$-Poisson
%space.
%\hfill $\diamond$
%\ere

\bex{ex-dressing} ({\bf Dressing actions})
Let $(G, \piG)$ be a Poisson Lie group with Lie bialgebra $(\g, \delta_\g)$, and let $(\g^*, \delta_{\g^*})$ and $(\d, \delta_\d)$ be respectively the dual and double Lie bialgebra of $(\g, \delta_\g)$. A Poisson map $\Phi: (Y, \piY) \to (G, \piG)$ {\it generates} a right Poisson action of $(\g^*, \delta_{\g^*})$ on $(Y, \piY)$ by
\[
\g^* \longrightarrow \V^1(Y), \;\; \xi \longmapsto \pi_{\sY}^\#(\Phi^*(\xi^R)), \hs \xi \in \g^*,
\]
where for $\xi \in \g^*$, $\xi^R$ denotes the right invariant $1$-form on $G$ with value $\xi$ at $e$. In particular, the identity map $G \to G$ generates \cite{STS2} the {\it right dressing action} $\vr$ of $(\g^*, \delta_{\g^*})$ on $(G, \piG)$ given by
\begin{equation}\lb{eq-vr-p-0}
\varrho: \; \g^* \lrw \V^1(G), \;\; \varrho(\xi) = \pi_\sG^\#(\xi^R), \hs \xi \in \g^*.
\end{equation}
The vector fields $\vr(\xi)$, $\xi \in \g^*$, are called the {\it right dressing vector fields} on $G$. For $g \in G$, denote by $\Ad_g$ the Adjoint action of $g$ on $\g$ and by $l_g$ and $r_g$ the left and right translations on $G$ by $g$. Then \cite{dr:homog} the adjoint representation of $\g$ on $\d$ integrates to an action of $G$ on $\d$, again denoted by $\Ad_g$ for $g \in G$, via
\begin{equation}\lb{eq-G-on-D-adjoint}
\Ad_{g^{-1}} (x+\xi) = \Ad_{g^{-1}} x  -l_{g^{-1}} \pi_\sG^\#(\xi^R)(g)+ \Ad_g^* \xi, \hs x \in \g, \, \xi \in \g^*, \, g \in G.
\end{equation}
Let $p_\g: \d \to \g$ be the projection with respect to the decomposition $\d = \g + \g^*$. By
\eqref{eq-G-on-D-adjoint},
\begin{equation}\lb{eq-vr-p}
\vr(\xi)(g)= -l_g p_\g \Ad_{g^{-1}}\xi, \hs \xi \in \g^*.
\end{equation}
The following ``multiplicativity" of the dressing vector fields $\varrho(\xi)$, $\xi \in \g^*$, can be proved directly using \eqref{eq-vr-p} (see also \cite[Lemma 3.4.8]{k-s:quantum} and \cite[Corollary 3.10]{lu:thesis}):
\begin{equation}\lb{eq-multi-vr}
\vr(\xi)(g_1g_2) = r_{g_2} \vr(\xi)(g_1) + l_{g_1} \vr(\Ad_{g_1}^* \xi)(g_2), \hs g_1, g_2 \in G.
\end{equation}
Let $\lam$ be the left action of $\g$ on $G$ given by
\[
\lam: \;\;\; \g \lrw \V^1(G), \;\;\; \lam(x) = x^R, \hs x \in \g,
\]
where for $x \in \g$, $x^R$ denotes the right invariant vector field on $G$ with value $x$ at $e$. It follows from  \eqref{eq-multi-vr} that one has the left Lie algebra action
\begin{equation}\lb{eq-varsig-G}
\varsigma: \;\; \d \lrw \V^1(G), \;\; \varsigma(x+\xi) =\lam(x) - \vr(\xi) = x^R -\pi_\sG^\#(\xi^R), \hs x \in \g, \, \xi \in \g^*,
\end{equation}
of $\d$ on $G$, which, by \eqref{eq-delta-d-g-g}, is a left Poisson action of the Lie bialgebra $(\d, \delta_\d)$ on $(G, \piG)$. \hfill $\diamond$
\eex

Recall from $\S$\ref{subsec-twists} the notion of twists of Lie bialgebras. 
%{\bf need to find a reference.}

\ble{le-twisting} {\rm {\bf (Twists of Poisson actions.)}} 
Assume that the Lie bialgebra $(\g, \delta_{\g, t})$ is the twist of a Lie bialgebra $(\g, \delta_{\g})$ by a twisting element $t \in \wedge^2 \g$ of $(\g, \delta_\g)$.

1) If $(Y, \piY, \lam)$ is a left $(\g, \delta_\g)$-Poisson space, then $(Y, \piY +\lam(t), \lam)$ is a left $(\g, \delta_{\g,t})$-Poisson space;

2) If $(Y, \piY, \rho)$ is a right $(\g, \delta_\g)$-Poisson space, then $(Y, \piY - \rho(t), \rho)$ is a right $(\g, \delta_{\g,t})$-Poisson space.
\ele

\begin{proof}
 Assume that $(Y, \piY, \lam)$ is a left $(\g, \delta_\g)$-Poisson space. Then
\[
[\piY + \lam(t), \; \piY + \lam(t)] = 2[\lam(t), \,\piY] + [\lam(t), \,\lam(t)] = 
-\lam\left(2\delta_{\g}(t) +[t, \,t]\right) = 0.
\]
Thus $\piY+\lam(t)$ is a Poisson structure on $Y$. It follows from \eqref{eq-g1g1} that $\lam$ is also a left Poisson action of the Lie bialgebra $(\g, \delta_{\g,t})$ on the Poisson manifold $(Y, \piY+\lam(t))$. 2) is proved similarly.
\end{proof}

\subsection{Poisson structures defined by quasitriangular $r$-matrices}\lb{subsec-Poi-r}
Let $\g$ be any Lie algebra and let $\sigma: \g \to \V^1(Y)$ be a left or right action of $\g$ on a manifold $Y$. For $r = \sum_i x_i \ot y_i \in \g \ot \g$, one has the $2$-tensor field $\sigma(r)$ on $Y$ given by 
\[
\sigma(r) = \sum_i \sigma(x_i) \ot \sigma(y_i).
\]
Let $s$ be the symmetric part of $r$. Clearly, $\sigma(r)$ is skew-symmetric, i.e., $\sigma(r)$ is a bi-vector field on $Y$, if and only if $\sigma(s) = 0$. A subspace $\c$ of $\g$ is said to be {\it coisotropic with respect to $s$} if $s^\#(\c^0) \subset \c$, where recall that $\c^0 = \{\xi \in \g^*: \la \xi, \c\ra = 0\} \subset \g^*$.

\ble{le-admi-r-equi}
For $r \in \g \ot \g$, the $2$-tensor field $\sigma(r)$ on $Y$ is skew-symmetric if and only if the stabilizer subalgebra of $\g$ at every $y \in Y$ is coisotropic with respect to the symmetric part $s$ of $r$.
\ele

\begin{proof}
Let $y \in Y$ and let $\sigma_y: \g \to T_y, \, \sigma_y(x) = \sigma(x)(y)$ for $x \in \g$. Then $\ker \sigma_y \subset \g$ is the stabilizer subalgebra of $\g$ at $y$, and $(\ker \sigma_y)^0 = {\rm Im} (\sigma_y^*: T^*Y \to \g^*)$. It follows from the definitions that $(\sigma(s)(y))^\# = \sigma_y \circ s^\# \circ  \sigma_y^*: T_y^*Y \to T_yY$. Thus $\sigma(s)(y) = 0$ if and only if $\ker \sigma_y$ is coisotropic with respect to $s$.
\end{proof}

Let $r \in \g \ot \g$ be a quasitriangular $r$-matrix on a Lie algebra $\g$ 
so that $(\g, r)$ is a quasitriangular Lie bialgebra. Let $\lam: \g \to \V^1(Y)$ be a left Lie algebra action of $\g$ on a manifold $Y$. 

\bpr{pr-admi-Poi-r}
If the two tensor field $\lam(r)$ on $Y$ is skew-symmetric, it is Poisson,
and $(Y, -\lam(r), \lam)$ is a left $(\g, r)$-Poisson space (see \deref{de-Poisson-space}).
\epr

\begin{proof}
Let $r = \Lam + s$, where $\Lam \in \wedge^2 \g$ and $s \in (S^2(\g))^\g$, and assume that $\lam(s) = 0$, so that $\lam(r)$ is a bi-vector field on $Y$. By \eqref{eq-CYB-Lam} - \eqref{eq-CYB-Lam-s}, the Schouten bracket of $\lam(r)$ with itself is given by
\[
[\lam(r), \, \lam(r)] = [\lam(\Lam), \, \lam(\Lambda)] = -\lam([\Lam, \Lam]) = 2\lam({\rm CYB}(s)).
\]
Let $y \in Y$, and let the notation be as in the proof of \leref{le-admi-r-equi}. Let $\alpha, \beta, \gamma \in T_y^*Y$. Then by \eqref{eq-CYB-s} and by the assumption that $\ker \lam_y$ is coisotropic with respect to $s$, one has
\[
\lam({\rm CYB}(s))(y) (\alpha, \beta, \gamma) = \la \lam_y^*(\alpha), \; [s^\#(\lam_y^*(\beta), \; 
s^\#(\lam_y^*(\gamma))]\ra =0.
\]
Since
$[\lam(x),  -\lam(\Lam)] =-[\lam(x),  \lam(\Lam)] = \lam([x, \Lam]) = \lam(\delta_r(x))$ for all $x \in \g$,
where $\delta_r: \g \to \wedge^2 \g$ is defined in \eqref{eq-delta-r},
$(Y, -\lam(r), \lam)$ is a left $(\g, r)$-Poisson space.
\end{proof}

\bde{de-admi-Poi-r} 
For a quasitriangular $r$-matrix $r$ on a Lie algebra $\g$ and a left Lie algebra action $\lam: \g \to V^1(Y)$, if $\lam(r)$ is skew-symmetric (so it is a Poisson structure on $Y$ by \prref{pr-admi-Poi-r}), we call $-\lam(r)$  {\it the Poisson structure on $Y$ defined by $\lam$ and $r$}.
\ede

\bre{re-Poi-r}
1) When $(\g, \lara_\g)$ is a quadratic Lie algebra and the quasitriangular $r$-matrix on $\g$ is defined by a Lagrangian 
splitting of $(\g, \lara_\g)$, the construction in \prref{pr-admi-Poi-r} has been given in \cite[$\S$2]{lu-yakimov}, and 
an interpretation of the construction in the framework of Courant algebroids is given in \cite{Li-Bland-Mein}.

2) An observation analogous to \prref{pr-admi-Poi-r} for $\g$-quasi-Poisson spaces is made in
\cite[$\S$2.1]{David-Severa:quasi-Hamiltonian-groupoids}.
The $(\g, r)$-Poisson space $(Y, -\lam(r), \lam)$ in \prref{pr-admi-Poi-r} corresponds to the $\g$-quasi-Poisson space $(Y, 0, \lam)$ under an equivalence via twisting between the category of $(\g, r)$-Poisson spaces and that of $\g$-quasi-Poisson spaces, which will be explained in $\S$\ref{subsec-mixed-fusion} and \exref{ex-Qy-0}. \hfill $\diamond$
\ere

\sectionnew{Mixed twists of direct product Lie bialgebras and mixed product Poisson structures}\lb{sec-mixed-twists}
\subsection{Mixed twists of direct product Lie bialgebras}\lb{subsec-mixed-twists}
Let $(\g_1, \delta_1), \ldots, (\g_n, \delta_n)$, $n \geq 1$,  be
Lie bialgebras, and let $(\g, \delta_\g)$ be the direct product Lie bialgebra, i.e., 
$\g = \g_1 \oplus \cdots \oplus \g_n$ with the direct product Lie bracket, and
\[
\delta_\g(x_1, \ldots, x_n) = (\delta_{\g_1}(x_1), \; \ldots, \delta_{\g_n}(x_n)), \hs (x_1, \ldots, x_n) \in \g,
\]
(see notation in $\S$\ref{subsec-intro-notation}). 
For  $1 \leq j \leq n$, let $p_j: \g \to \g_j$ be the $j$'th projection. Recall from $\S$\ref{subsec-twists} the notion of twists of Lie bialgebras. In particular, if $t \in \wedge^2 \g$ is a twisting element for $(\g, \delta_\g)$,
one has the twisted Lie bialgebra $(\g, \delta_{\g, t})$, where
$\delta_{\g, t}(x) = \delta_\g(x) - [x, t]$ for $x \in \g$.

\bde{de-mixed-twists}
A twisting element $t \in \wedge^2 \g$ of the direct product Lie bialgebra $(\g, \delta_\g)$ is said to be {\it mixed} if
$p_j(t) = 0$ for each $1 \leq j \leq n$, and in this case the twisted Lie bialgebra $(\g, \delta_{\g, t})$ 
(see \eqref{eq-g1g1}) is 
called a {\it mixed twist} of the direct product Lie bialgebra $(\g, \delta_\g)$. 
Note that when $n = 1$, only the zero element is a mixed twisting element of $(\g, \delta_\g)$.
\ede

For $J = \{j_1, \ldots, j_k\} \subset \{1, \ldots, n\}$ with $1 \leq j_1 < \cdots < j_k \leq n$, let 
$\g_{{}_\sJ}  =\g_{j_1} \oplus \cdots \oplus \g_{j_k}$, and let
\begin{equation}\lb{eq-p-sJ-0}
p_\sJ: \;\; \g \lrw \g_{{}_\sJ}, \;\; (x_1, \ldots, x_n) \longmapsto (x_{j_1}, \ldots, x_{j_k}).
\end{equation}
Denote by $(\g_{{}_\sJ}, \, \delta_{\g_{{}_\sJ}})$ the direct product Lie bialgebra
of $(\g_{j_1}, \delta_{\g_{j_1}}), \ldots, (\g_{j_k}, \delta_{\g_{j_k}})$.  

\ble{le-p-sJ}
If $t \in \wedge^2 \g$ be a mixed twisting element of the direct product Lie bialgebra $(\g, \delta_\g)$, then $p_\sJ(t) \in \wedge^2 \g_{{}_\sJ}$ is a mixed twisting element of the direct product Lie bialgebra
$(\g_{{}_\sJ}, \delta_{\g_{{}_\sJ}})$ for
any $J \subset \{1, \ldots, n\}$, and 
\begin{equation}\lb{eq-p-sJ}
p_\sJ: \;\; (\g, \, \delta_{\g, t}) \lrw (\g_{{}_\sJ}, \; \delta_{\g_{{}_\sJ}, p_\sJ(t)})
\end{equation}
is a Lie bialgebra homomorphism. In particular, $p_j: (\g, \delta_{\g, t}) \to (\g_j, \delta_{\g_j})$ is
a Lie bialgebra homomorphism for each $1 \leq j \leq n$.
\ele

\begin{proof}
Let $t_\sJ = p_\sJ(t)$. As $p_\sJ: (\g, \delta_\g) \to (\g_{\sJ}, \delta_{\g_{{}_\sJ}})$ is a Lie bialgebra homomorphism, it follows from 
$2\delta_{\g}(t)+[t, t]=0$ that $2\delta_{\g_{{}_\sJ}}(t_\sJ) +[t_\sJ, t_\sJ]=0$, so $t_\sJ$ is a twisting
element, and in fact a mixed twisting element, of $(\g_{{}_\sJ}, \delta_{\g_{{}_\sJ}})$. It is also clear from the definitions of $\delta_{\g, t}$ and
$\delta_{\g_{{}_\sJ}, t_\sJ}$ that $p_\sJ$ in \eqref{eq-p-sJ} is a Lie bialgebra homomorphism.
\end{proof}

If $(\g_1, r_1), \ldots, (\g_n, r_n)$ are quasitriangular Lie bialgebras, one has the 
direct product quasitriangular Lie bialgebra $(\g, r)$, where again 
$\g = \g_1 \oplus \cdots \oplus \g_n$ is the direct product Lie algebra  and $r = (r_1, \ldots, r_n) \in \g \otimes \g$.
Let $\delta_{\g_j} = \delta_{r_j}$ 
(see \eqref{eq-delta-r}) for $1 \leq j \leq n$.  By 
\leref{le-t-r}, $t \in \wedge^2 \g$ is a twisting element for the direct
product Lie algebra $(\g, \delta_\g)$ of $(\g_1, \delta_{\g_1}), \ldots, (\g_n, \delta_{\g_n})$ if
and only if $r - t$ is 
quasitriangular $r$-matrix on $\g$.

\bde{de-mixed-t-r}
If $t \in \wedge^2 \g$ is such that $r - t$, where $r = (r_1, \ldots, r_n)$, is 
quasitriangular $r$-matrix on $\g$ and $p_j(t) = 0$ for each $1 \leq j \leq n$, 
we call $r - t$ a {\it mixed twist of $r$}, and 
we also call the quasitriangular Lie bialgebra
$(\g, r - t)$ a {\it mixed twist of the direct product quasitriangular Lie bialgebra $(\g, r)$}.
\ede

\subsection{Mixed twists and mixed product Poisson structures}\lb{subsec-mixed-mixed}
Let $(\g_1, \delta_{\g_1}), \ldots, (\g_n, \delta_{\g_n})$, $n \geq 2$, be Lie bialgebras and let
$(\g, \delta_\g)$ be the direct product Lie bialgebra as in $\S$\ref{subsec-mixed-twists}. Assume that 
$t \in \wedge^2\g$ is a mixed twisting element of $(\g, \delta_\g)$, and let $(\g, \delta_{\g, t})$ be 
the corresponding twisted Lie bialgebra. 
As a direct application of \leref{le-twisting} and an immediate consequence of $t$ being a mixed twisting element,
we have

\bpr{pr-mixed-mixed}
1) Let  $(Y_j, \pi_i, \lam_j)$, $1 \leq j \leq n$,
be left $(\g_j, \delta_{\g_j})$-Poisson spaces, and let $Y = Y_1 \times \cdots \times Y_n$ be the 
product manifold, equipped with the product Lie algebra action $\lam = (\lam_1, \ldots, \lam_n)$ of $\g$, i.e.,
\[
\lam: \;\; \g \lrw \V^1(Y), \;\; \lam(x_1, \ldots, x_n) = (\lam_1(x_1), \; \ldots, \; \lam_n(x_n)), \hs x_j \in \g_j.
\]
Then $\piY = (\pi_1,\; \ldots, \; \pi_n) +\lam(t)$ is a mixed product Poisson structure on $Y$, and
$(Y, \piY, \lam)$ is a left $(\g, \delta_{\g, t})$-Poisson spaces;

2) Let   $(Y_j, \pi_i, \rho_j)$, $1 \leq j \leq n$, 
be right $(\g_j, \delta_{\g_j})$-Poisson space, and let 
\[
\rho = (\rho_1, \ldots, \rho_n): \;\; \g \lrw \V^1(Y), \;\; \rho(x_1, \ldots, x_n) = (\rho_1(x_1), \; \ldots, \; \rho_n(x_n)), \hs x_j \in \g_j,
\]
where $Y = Y_1 \times \cdots \times Y_n$,
Then $\piY = (\pi_1,\; \ldots, \; \pi_n) -\rho(t)$ is a mixed product Poisson structure on $Y$, and
$(Y, \piY, \rho)$ is a right $(\g, \delta_{\g, t})$-Poisson spaces;

3) If $r_j \in \g_j \otimes \g_j$ is a quasitriangular structure for $(\g_j, \delta_{\g_j})$ and $\pi_j = -\lam_j(r_j)$
in 1) (resp. $\pi_j = \rho_j(r_j)$ in 2)) for $1 \leq j \leq n$, then $\piY = -\lam(r-t)$ in 1) (resp. $\piY = \rho(r-t)$ in 2)), where
$r = (r_1, \ldots, r_n)$.
\epr

\prref{pr-mixed-mixed} is the key to all the constructions in this paper. Indeed, 
the mixed product Poisson structures constructed in the following
$\S$\ref{sec-mixed}, $\S$\ref{sec-polyuble}, and $\S$\ref{sec-mixed-quasi} are all incidences of 
mixed twists of direct product Lie bialgebras and of their direct product Poisson actions in the spirit of
\prref{pr-mixed-mixed}.

\sectionnew{Two-fold mixed product Poisson structures associated to Lie bialgebras}\lb{sec-mixed}
Although a special case of what to be covered in $\S$\ref{sec-polyuble}, the
{\it two-fold} mixed product Poisson structure construction in this section is more basic, so we discuss it
first. The construction also appeared in \cite{Victor:thesis}.

\subsection{The construction}\lb{subsec-mixed-product}
Let $(\g, \delta_\g)$ be a Lie bialgebra, with dual Lie bialgebra $(\g^*, \delta_{\g^*})$ and double Lie bialgebra
$(\d, \delta_\d)$, and let $(\d', \delta_{\d'})$ be the dual Lie bialgebra of $(\d, \delta_\d)$. 
Recall from $\S$\ref{subsec-double} that as a Lie algebra $\d' = \g^* \oplus \g$ with the direct product Lie bracket. 
Define the direct product Lie co-bracket $\delta_{\d'}^{\prime}$ on $\d'$ by
\begin{equation}\lb{eq-delta-0}
\delta_{\d'}^\prime(\xi, x) = (\delta_{\g^*}(\xi), \;-\delta_\g(x)), \hs \hs \xi \in \g^*, \;x \in \g.
\end{equation}
Let $\{x_i\}_{i=1}^{m}$ be any basis of $\g$ and let $\{\xi_i\}_{i=1}^{m}$ be the dual basis for $\g^*$. 
Let 
\begin{equation}\lb{eq-t-1}
t = \sum_{i=1}^m (\xi_i, 0) \wedge (0, x_i) \in \wedge^2 \d'.
\end{equation}

\ble{le-twist-d-prime}
The element $-t \in \wedge^2 \d'$ is a mixed twisting element for the direct product Lie bialgebra
$(\d', \delta_{\d'}^\prime)$, and the twist of $(\d', \delta_{\d'}^\prime)$ by $-t$ is the Lie bialgebra $(\d', \delta_{\d'})$.
\ele

\begin{proof}
It is straightforward to check that $\delta_{\d'}^\prime$ is defined (see \deref{de-lag-splitting}) by the Lagrangian splitting $\d^2 = \d' + (\g \oplus \g^*)$  of the quadratic Lie algebra $(\d^2, \lara_{\d^2})$, and that the twist by $-t$ of $\d^2 = \d' + (\g \oplus \g^*)$ 
is the Lagrangian splitting $\d^2 = \d' + \d_{\rm diag}$ which defines the Lie bialgebra $(\d', \delta_{\d'})$.
\end{proof}

Let $(X, \piX, \rho)$ be a right $(\g^*, \delta_{\g^*})$-Poisson space and $(Y, \piY, \lam)$ a
left $(\g, \delta_\g)$-Poisson space. Let again $\{x_i\}_{i=1}^{m}$ be any basis of $\g$ and $\{\xi_i\}_{i=1}^{m}$ the dual basis for $\g^*$.  Define the bi-vector field $\piX \times_{(\rho, \lam)} \piY$ on the product manifold $X \times Y$ by 
\begin{equation}\lb{eq-pi-mixed-1}
\piX \times_{(\rho, \lam)} \piY  
=(\piX, \; \piY) -\sum_{i=1}^m (\rho(\xi_i), \, 0) \wedge (0, \; \lam(x_i)).
\end{equation}
The following \prref{pr-mixed-is-poisson} is now a direct consequence of 
\reref{re-quotient-action-0} and \prref{pr-mixed-mixed}.

\bpr{pr-mixed-is-poisson}
The bi-vector field $\piX \times_{(\rho, \lam)} \piY$ on $X \times Y$ is Poisson, and the map
\[
\rho_0: \;\; \d' \longrightarrow \V^1(X \times Y), \;\;\; \rho_0(\xi, x) = (\rho(\xi), \, 0) + (0, \, -\lam(x)),
\hs \xi \in \g^*, \, x \in \g
\]
is a right Poisson action of the Lie bialgebra $(\d', \delta_{\d'})$ on $(X \times Y, \,   \piX \times_{(\rho, \lam)} \piY)$.
\epr
 
\bde{de-mixed}
We call $\piX \times_{(\rho, \lam)} \piY$ the {\it mixed product Poisson structure} on $X \times Y$  {\it associated to the pair $(\rho, \lam)$} of Poisson actions. 
\ede

\bre{re-mixed-15} ({\bf Reduction to sub-Lie bialgebras, I.})
In the setting of \prref{pr-mixed-is-poisson}, assume that $(\p, \delta_\p)$ is a sub-Lie bialgebra of $(\g, \delta_\g)$ such that $\rho(\xi) = 0$ for all $\xi \in \p^0\subset \g^*$. Let $\lam_\p: \p \to \V^1(Y)$ be the restriction of $\lam$ to $\p$ and $\rho_{\p^*}: \p^* \cong \g^*/\p^0 \to \V^1(X)$  the Lie algebra action induced by $\rho$. Then $(Y, \piY, \lam_\p)$ is a left $(\p, \delta_\p)$-Poisson space and $(X, \piX, \rho_{\p^*})$ a right $(\p^*, \delta_{\p^*})$-Poisson space. Choosing a basis of $\p$ and extending it to a basis of $\g$, one sees  that 
$\piX \times_{(\rho, \lam)} \piY =\piX \times_{(\rho_{\p^*}, \lam_\p)} \piY$. \hfill $\diamond$

\ere

\bre{re-mixed-term}
Given a right $(\g^*, \delta_{\g^*})$-Poisson space $(X, \piX, \rho)$ and a left $(\g, \delta_\g)$-Poisson space $(Y, \piY, \lam)$,  let $\hat{\rho}: T^*X \to \g$ and $\hat{\lam}: T^*Y \to \g^*$ be respectively given by
\begin{align}\lb{eq-rho-hat}
\la \hat{\rho}(\alpha_p), \, \xi\ra& = \la \alpha_p, \, \rho(\xi)(p)\ra, \hs p \in X, \, \alpha_p \in T_p^*X,\, \xi \in \g^*,\\
\lb{eq-lam-hat}
\la \hat{\lam}(\beta_q),\,  x\ra& = \la \beta_q, \, \lam(x)(q)\ra, \hs q \in Y, \, \beta_q \in T_q^*Y, \, x \in \g,
\end{align}
and let $p_\sX: X \times Y \to X$ and $p_\sY: X \times Y \to Y$ be the two projections. Then  a Poisson structure $\pi$ on $X \times Y$ is equal to $\piX \times_{(\rho, \lam)} \piY$ if and only if the two projections $p_\sX: (X \times Y, \, \pi) \rightarrow (X, \piX)$ and  $p_\sY: (X \times Y, \, \pi) \rightarrow (Y, \piY)$ are Poisson and 
\begin{equation}\lb{eq-pi-mixed-2}
(\piX \times_{(\rho, \lam)} \piY) (p_\sX^*\alpha, \, p_\sY^*\beta) 
=-\la \hat{\rho}(\alpha), \, \hat{\lam}(\beta)\ra.
\end{equation}
for any $1$-form $\alpha$ on $X$ and $1$-form $\beta$ on $Y$. \hfill $\diamond$
\ere

The next \leref{le-mixed-poi-sub} follows directly from the definition of the Poisson structure $\piX \times_{(\rho, \lam)} \piY$.

\ble{le-mixed-poi-sub}
Let $\piX \times_{(\rho, \lam)} \piY$ be the mixed product Poisson structure on $X \times Y$ defined by the pair $(\rho, \lam)$.  If $X_1 \subset X$ is a $\g^*$-invariant Poisson submanifold of $(X, \piX)$ and $Y_1 \subset Y$  a $\g$-invariant Poisson submanifold of $(Y, \piY)$, then $X_1 \times Y_1$ is a Poisson submanifold of $(X \times Y, \, \piX \times_{(\rho, \lam)} \piY)$.
\ele

\bex{ex-semi-product} 
Let $(G, \piG)$ be a Poisson Lie group with Lie bialgebra $(\g, \delta_\g)$ and dual Lie bialgebra $(\g^*, \delta_{\g^*})$, and let $\lam$ be the left Poisson action of $(\g, \delta_\g)$ on $(G, \piG)$ given by $\lam(x) = x^R$ for $x \in \g$ (see \exref{ex-dressing}). For a right $(\g^*, \delta_{\g^*})$-Poisson space $(X, \piX, \rho)$, the mixed product Poisson structure $\piX \times_{(\rho, \lam)} \piG$ on $X \times G$ is called a {\it semi-direct product Poisson structure} on $X \times G$ in \cite{lu:thesis}.  In this case, the right Poisson action of $(\g^*, \delta_{\g^*})$ on $(X\times G, \, \piX \times_{(\rho, \lam)} \piG)$ generated (see \exref{ex-dressing}) by the Poisson map $\kappa_\sG:  (X\times G,  \piX \times_{(\rho, \lam)} \piG) \to (G, \piG), (x, g) \to g$, is the diagonal action 
\[
\g^* \lrw \V^1(X \times G), \;\; \xi \longmapsto (\rho(\xi), 0) + (0, \vr(\xi)), \hs \xi \in \g^*,
\]
where $\vr$ is the right dressing action of $(\g^*, \delta_{\g^*})$ on $(G, \piG)$. \hfill $\diamond$
\eex

\bex{dual of double mixed prod}
Let $(G, \piG)$ be a Poisson Lie group with Lie algebra $(\g, \delta_\g)$ and $(G^*, \pi_{\sG^*})$ any dual Poisson Lie group of $(G, \piG)$. Consider the dual pair of Poisson Lie groups
$$
(H, \pi_H) = (G^{\rm op}, \piG) \times (G, \piG) \hs \text{and} \hs (H^*, \pi_{\sH^*}) = (G^*, -\pi_{\sG^*}) \times (G^*, \pi_{\sG^*}),
$$
where $G^{\rm op}$ denotes $G$ with the opposite group structure, and the right and left Poisson actions 
\begin{align*} 
\rho:&\;\;\; (G^*, \pi_{\sG^*}) \times (H^*, \pi_{\sH^*})  \longrightarrow (G^*, \pi_{\sG^*}),\;\;\;(g^*, \, (h_1, h_2))\longmapsto  h_1^{-1}g^*h_2, \hs g^*,  h_1, h_2 \in G^*,       \\
\lambda: & \;\;\;(H, \pi_\sH) \times (G, -\piG) \longrightarrow (G, -\piG),\;\;\; ((g_1, g_2), g) \longmapsto g_1^{-1} g g_2^{-1}, \hs g,  g_1, g_2 \in G.
\end{align*}
Then $(D', \pi_{\sD'}) = (G^* \times G, \pi_{\sG^*} \times_{(\rho, \lambda)} (-\piG))$ is a Poisson Lie group with Lie bialgebra $(\d', \delta_{\d'})$, where recall from $\S$\ref{subsec-double} that $(\d', \delta_{\d'})$ is the dual Lie bialgebra of the double Lie bialgebra $(\d, \delta_\d)$ of $(\g, \delta_\g)$. Indeed, identify the Lie algebra of $D'$ with $\d' = \g^* \oplus \g$. If $\{x_i\}_{i=1}^m$ is a basis of $\g$ and $\{\xi_i\}_{i=1}^m$  the dual basis in $\g^*$, the mixed term of $\pi_{\sG^*} \times_{(\rho, \lambda)} (-\piG)$ is
\[
\mu_{(\rho, \lambda)}  = -\sum_{i=1}^m (\xi_i^R, 0) \wedge(0, x_i^R) - (\xi_i^L, 0) \wedge (0, x_i^L)  
   = -t^L + t^R,
\]
where $t \in \wedge^2\d'$ is defined in \eqref{eq-t-1}. Thus $\pi_{\sD'} = (\pi_{\sG^*}, \, -\piG) - t^L + t^R$ is multiplicative. Moreover, since the Lie bialgebra $(\d', \delta_{\d'})$ is the twisting by $t$ of the Lie bialgebra $(\d', \delta_0)$ as in the proof of \prref{pr-mixed-is-poisson}, the Lie bialgebra of $(D', \pi_{\sD'})$ is $(\d', \delta_{\d'})$. \hfill $\diamond$
\eex

%Further examples of $2$-fold mixed product Poisson structures  will be given in $\S$\ref{subsec-mixed-actions}.

\subsection{Morphisms between mixed product Poisson structures} 
For $i = 1, 2$, let $(\g_i, \delta_{\g_i})$ be a Lie bialgebra and $(\g_i^*, \delta_{\g_i^*})$ its dual Lie bialgebra,  let $(Y_i, \pi_{\sY_i}, \lam_i)$ be a left $(\g_i, \delta_{\g_i})$-Poisson space and $(X_i, \pi_{\sX_i}, \rho_i)$ a right $(\g_i^*, \delta_{\g_i^*})$-Poisson space. Given a triple $(\Phi, \Psi, \phi)$, where 
\[
\Phi: \;\;\ (X_1, \pi_{\sX_1}) \longrightarrow  (X_2, \pi_{\sX_2}) \hs \mbox{and} \hs
\Psi: \;\; (Y_1, \pi_{\sY_1}) \longrightarrow  (Y_2, \pi_{\sY_2})
\]
are Poisson maps and $\phi: (\g_1, \delta_{\g_1}) \to (\g_2, \delta_{\g_2})$ is a Lie bialgebra homomorphism such that for every  $x_1 \in \g_1$, the vector fields $\lam_1(x_1) \in \V^1(Y_1)$ and $\lam_2(\phi(x_1))\in \V^1(Y_2)$ are $\Psi$-related and that for every $\xi_2 \in \g_2^*$, the vector fields $\rho_1(\phi^*(\xi_2)) \in \V^1(X_1)$ and $\rho_2(\xi_2)\in \V^1(X_2)$ are $\Phi$-related, let
\[
\Phi \times \Psi: \;\;\; X_1 \times Y_1 \lrw X_2 \times Y_2, \;\;\; (p_1, q_1) \longmapsto (\Phi(p_1), \, \Psi(q_1)),
\hs p_1 \in X_1, q_1 \in Y_1.
\]

\bpr{pr-mixed-morphism}
$\Phi \times \Psi: (X_1 \times Y_1, \; \pi_{\sX_1} \times_{(\rho_1, \lam_1)} \pi_{\sY_1}) \rightarrow (X_2 \times Y_2, \; \pi_{\sX_2} \times_{(\rho_2, \lam_2)} \pi_{\sY_2})$ is Poisson.
\epr

\begin{proof}
Let  $\mu_i = \pi_{\sX_i} \times_{(\rho_i, \lam_i)} \pi_{\sY_i} -(\pi_{\sX_i}, 0) - (0, \pi_{\sY_i})$ be the mixed terms of $\pi_{\sX_i} \times_{(\rho_i, \lam_i)} \pi_{\sY_i}$, $i = 1, 2$. Since $\Phi$ and $\Psi$ are Poisson, it suffices to show that $\mu_1$ and $\mu_2$ are $(\Phi \times \Psi)$-related.

Fix  $(p_1, q_1) \in X_1 \times Y_1$ and let $p_2 = \Phi(p_1)$ and $q_2 = \Psi(q_1)$. For $i = 1, 2$, let $\hat{\rho}_i:  T^*_{p_i} X_i \to \g_i$ and $\hat{\lam}_i: T^*_{q_i} Y_i \to \g_i^*$ be defined as in \reref{re-mixed-term}. Then the assumptions on $\Phi$ and $\Psi$ imply that
\[
\phi \hat{\rho}_1 \Phi^* = \hat{\rho}_2: \; T_{p_2}^* X_2 \lrw \g_2 \hs
\mbox{and} \hs \hat{\lam}_1 \Psi^* = \phi^* \hat{\lam}_2: \; T_{q_2}^*Y_2 \lrw \g_1^*.
\]
It follows that for any $(\alpha_2, \beta_2) \in T_{p_2}^*X_2 \times  T_{q_2}^* Y_2$, one has 
\begin{align*}
\mu_1(p_1, q_1)(\Phi^*(\alpha_2), \Psi^*(\beta_2))&=-\la \hat{\rho}_1(\Phi^*(\alpha_2)), \; \hat{\lam}_1(\Psi^*(\beta_2))\ra
= -\la \phi \hat{\rho}_1 \Phi^*(\alpha_2), \; \hat{\lam}_2(\beta_2)\ra\\
&=-\la \hat{\rho}_2 (\alpha_2), \; \hat{\lam}_2(\beta_2)\ra
=\mu_2(p_2, q_2)(\alpha_2, \beta_2).
\end{align*}
Thus $(\Phi \times \Psi)(\mu_1(p_1, q_1)) = \mu_2(p_2, q_2)$.
\end{proof}

We now look at a situation in which the mixed term of a mixed product Poisson structure vanishes under a Poisson morphism.
Let $((\g, \delta_\g), \,(\g^*, \delta_{\g^*}))$ be a dual pair of Lie bialgebras.

\ble{le-vanish}
Let $(X, \piX, \rho)$ be a right $(\g^*, \delta_{\g^*})$-Poisson space and $(Y, \piY, \lam)$ a
left $(\g, \delta_\g)$-Poisson space, and assume that $\Phi: (X, \piX) \to (X', \pi_{\sX'})$ and 
$\Psi: (Y, \piY) \to (Y', \pi_{\sY'})$ are Poisson morphisms. If there is a subspace $\a \subset \g$ such that
$\Psi(\lam(x)) = 0$ for all $x \in \a$ and $\Phi(\rho(\xi)) = 0$ for all $\xi \in \a^0\subset \g^*$,
then $(\Phi \times \Psi)(\piX \times_{(\rho, \lam)} \piY) = \pi_{\sX'} \times \pi_{\sY'}$.
\ele

\begin{proof}
Let $\{x_i\}_{i=1}^m$ be a basis of $\g$ such that $\{x_i\}_{i=1}^{k =\dim \a}$ is a basis of $\a$, and
let $\{\xi_i\}_{i=1}^m$ be the dual basis of $\g^*$. As $\Psi(\lam(x_i)) = 0$ for $1 \leq i \leq k$ and $\Phi(\rho(\xi_i)) =0$
for $k+1 \leq i \leq n$, the mixed term of $\piX \times_{(\rho, \lam)} \piY$ is mapped to zero by $\Phi \times \Psi$.
\end{proof}

\sectionnew{Polyubles of Lie bialgebras and mixed product Poisson structures}\lb{sec-polyuble}

\subsection{Polyubles of Lie bialgebras}\lb{subsec-polyuble} 
Let $(\g, \delta_\g)$ be a Lie bialgebra, with dual Lie bialgebra $(\g^*, \delta_{\g^*})$, double Lie algebra $(\d, \lara_\d)$, and double Lie bialgebra $(\d, \delta_\d)$.  For an integer $n \geq 1$, consider the quadratic Lie algebra $(\d^n, \lara_{\d^n})$, where $\d^n = \d \oplus \cdots \oplus \d$ has the direct product Lie bracket, its elements denoted as $(a_1, \ldots, a_n)$ with $a_j \in \d$, and the bilinear form $\lara_{\d^n}$ is given by
\begin{equation} \lb{lara_n}
\la (a_1, \ldots, a_n), \;(a'_1, \ldots, a'_n) \ra_{\d^n} = \sum_{j = 1}^n (-1)^{j+1} \la a_j, \;a'_j \ra_\d, \hs a_j, \;
a_j^\prime \in \d. 
\end{equation}
Let again $\d_{\rm diag} = \{(a,a) : a \in \d \} \subset \d \oplus \d$. For an integer $n \geq 1$, define
\begin{align} \lb{eq-g-odd}
&\g_{(2n-1)} = \overbrace{\d_{\rm diag} \oplus \cdots \oplus \d_{\rm diag}}^{n-1}  \oplus \;\g \subset \d^{2n-1}, \hs\hs
\g_{(2n-1)}^* = \g^* \oplus \overbrace{\d_{\rm diag} \oplus \cdots \oplus \d_{\rm diag}}^{n-1}\subset \d^{2n-1}, \\
\lb{eq-g-even}
&\g_{(2n)} = \overbrace{\d_{\rm diag} \oplus \cdots \oplus \d_{\rm diag}}^n \subset \d^{2n}, \hs \hs
\g_{(2n)}^* = \g^* \oplus \overbrace{\d_{\rm diag} \oplus \cdots \oplus \d_{\rm diag}}^{n-1} \oplus \;\g \subset \d^{2n}.
\end{align}
It is straightforward to check that for each $n \geq 1$,
\begin{equation}\lb{eq-n-splitting}
\d^n = \g_{(n)} + \g_{(n)}^*
\end{equation}
is a Lagrangian splitting of $(\d^n, \lara_{\d^n})$. Denote by $\left(\g_{(n)}, \delta_{\g_{(n)}}\right)$ and $\left(\g_{(n)}^*, \delta_{\g_{(n)}^*}\right)$ the dual pair of Lie bialgebras defined by the Lagrangian splitting in \eqref{eq-n-splitting}, and denote by $(\d^n, \delta_{\d^n})$ the corresponding double Lie bialgebra. Note that when $n = 1, 2$, the splitting in \eqref{eq-n-splitting} are respectively $\d = \g + \g^*$ and $\d^2 = \d_{\rm diag} + \d'$, where $\d' = \g^* \oplus \g$, and we thus recover the two dual pairs of Lie bialgebras $((\g, \delta_\g), (\g^*, \delta_{\g^*}))$ and $((\d, \delta_\d), (\d', \delta_{\d'}))$ (see $\S$\ref{subsec-double}).

\bde{de-uble}
The Lie bialgebra $(\g_{(n)}, \delta_{\g_{(n)}})$ is called the {\it $n$-uble} of $(\g, \delta_\g)$.
\ede

\bre{re-ubles} 
1) The $1$-uble and $2$-uble of a Lie bialgebra are respectively the Lie bialgebra itself and its Drinfeld double. For $n = 3$, {\it triples} of quasitriangular Lie bialgebras were constructed by Grabowski in \cite{grabowski}. In the factorizable case, it follows from \cite[Theorem 4.5]{grabowski} and our \prref{pr-uble-mixed-power} in $\S$\ref{sec-mixed-quasi}
that the triple of a Lie bialgebra is isomorphic to its $3$-uble.

2) {\it Polyubles} for factorizable Lie bialgebras and factorizable Poisson Lie groups \cite{RT:factorizable} were constructed in \cite{fock-rosly:polyubles, fock-rosly:r-matrix}, and we also refer to \cite{fock-rosly:polyubles, fock-rosly:r-matrix, Li-Bland-Severa} for interpretations of such polyubles as moduli spaces of flat connections on (quilted) surfaces. Due to limitation of space, we will only present in this section  properties of polyubles (for arbitrary Lie bialgebras) that are needed for the rest of the paper. 
\hfill $\diamond$
%A more detailed study on polyubles will be given elsewhere. 
\ere

We now show that the Lie bialgebras $\left(\g_{(n)}, \delta_{\g_{(n)}}\right)$ for $n \geq 3$ and 
$\left(\g_{(n)}^*, \delta_{\g_{(n)}^*}\right)$ for $n \geq 2$ are mixed twists of direct product Lie bialgebras
(see $\S$\ref{sec-mixed-twists}).
Note that the Lie brackets on $\g_{(n)}$ for $n \geq 3$ and on $\g_{(n)}^*$ for $n \geq 3$ are indeed direct product Lie brackets:
identify $\d_{\rm diag} \cong \d$ via $(a, a) \mapsto a$ for $a \in \d$ in \eqref{eq-g-odd} and \eqref{eq-g-even}, one has the identification of Lie algebras
\begin{equation}\lb{eq-iden-ggg}
\g_{(2n+2)} \cong \d^{n+1}, \hs 
\g_{(2n)}^* = \g^* \oplus \d^{n-1} \oplus \g, \hs 
\g_{(2n+1)} = \d^n  \oplus \g, \hs 
\g_{(2n+1)}^* = \g^* \oplus \d^n, \hs n \geq 1. 
\end{equation}

\bnota{nota-tn}
Let $n\geq 1$, and for $1 \leq j \leq n+1$ and $a \in \d$, let $(a)_j = (0, \ldots, 0, a, 0, \ldots, 0) \in \d^{n+1}$, 
where
$a$ is at the $j$'th entry of $(a)_j$. Let $\{x_i\}_{i=1}^m$ be any basis of $\g$ and $\{\xi_i\}_{i=1}^m$ the dual
basis of $\g^*$. Set
\begin{equation}\lb{eq-t-n}
t_{n+1} = \sum_{1 \leq j < k \leq n+1} \; \sum_{i=1}^m (\xi_i)_j \wedge (x_i)_k \in \wedge^2 (\d^{n+1}).
\end{equation}
Note that $t_{n+1}$ in fact lies in $\wedge^2 (\g^* \oplus \d^{n-1} \oplus \g) \subset \left(\wedge^2(\d^n \oplus \g)\right)
\cap \left(\wedge^2 (\g^* \oplus \d^n)\right) 
\subset \wedge^2 (\d^{n+1})$.   
Using the identifications in \eqref{eq-iden-ggg}, we can regard $t_{n+1}$ as an element in $\wedge^2 V$, where $V$
can be either of  
$\g_{(2n+2)}, \;$, $\g_{(2n)}^*, \;$
$\g_{(2n+1)}$, or $\g_{(2n+1)}^*$.
\enota

\ble{le-uble-mixed-1} Under the Lie algebra identification $\g_{(2n+2)}\cong \d^{n+1}$ in 
\eqref{eq-iden-ggg}, where 
$n \geq 1$, the Lie bialgebra
$\left(\g_{(2n+2)}, \;\delta_{\g_{(2n+2)}}\right)$ is the mixed twist of 
the direct product Lie bialgebra $\left(\g_{(2n+2)}, \;\delta_{\g_{(2n+2)}}^\prime\right)$ by the mixed twisting
element $t_{n+1}\in \wedge^2 \g_{(2n+2)}$, where 
\[
\delta_{\g_{(2n+2)}}^\prime(a_1, \,\ldots, \,a_{n+1}) =
(\delta_\d(a_1), \; \ldots, \;\delta_{\d}(a_{n+1})), \hs a_1, \, \ldots, a_{n+1} \in \d.
\]
\ele

\begin{proof}
Consider the two Lagrangian splittings of the quadratic Lie algebra $(\d^{2n+2}, \lara_{\d^{2n+2}})$:
\begin{equation}\lb{eq-two-splittings-n}
\d^{2n+2} = \g_{(2n+2)} + \overbrace{\d' \oplus \cdots \oplus \d'}^{n+1}
= \g_{(2n+2)} + \g_{(2n+2)}^*,
\end{equation}
where recall from $\S$\ref{subsec-double} that $\d^\prime = \g^* \oplus \g \subset \d^2$. 
Clearly,
the co-bracket on $\g_{(2n+2)}$ defined by the first splitting 
in \eqref{eq-two-splittings-n} is $\delta_{\g_{(2n+2)}}^\prime$. It is straightforward to check that 
second Lagrangian splitting in \eqref{eq-two-splittings-n} is the twist by $t_{n+1} \in \wedge^2 \g_{(2n+2)}$ of the first splitting.  It is also evident that the twisting element $t_{n+1}$ is mixed. 
\end{proof}

Let $r_\d \in \d \otimes \d$ be again the quasitriangular $r$-matrix on $\d$ defined by the Lagrangian splitting
$\g = \g + \g^*$ of the quadratic Lie algebra $(\d, \lara_\d)$. 
For $n \geq 1$, let $r_{\d^{n+1}} \in \d^{n+1} \otimes \d^{n+1}$ be the quasitriangular $r$-matrix on $\d^{n+1}$
defined by the Lagrangian splitting $\d^{n+1} = \g_{(n+1)} + \g_{(n+1)}^*$ of $(\d^{n+1}, \lara_{\d^{n+1}})$.
Define the $(n+1)$'st {\it alternating sum} of $r_\d$ by
\[
{\rm Alt}^{n+1}(r_\d) = \begin{cases} (r_\d, \; -r^{21}_\d, \; r_\d, \; \cdots, \; -r^{21}_\d, \; r_\d)\in \d^{n+1}
\otimes \d^{n+1}, & \;\;\; n \;
\mbox{even},\\
(r_\d, \; -r^{21}_\d,  \; \cdots, \; r_\d, \; -r^{21}_\d)\in \d^{n+1} \otimes \d^{n+1}, & \;\;\; n \; \mbox{odd}.\end{cases}
\]

\ble{le-double-uble-r}
For $n \geq 1$, the quasitriangular Lie bialgebra $(\d^{n+1}, \, r_{\d^{n+1}})$
is the mixed twist by $t_{n+1} \in \wedge^2 \d^{n+1}$ of the direct product quasitriangular Lie bialgebra
$(\d^{n+1}, \, {\rm Alt}^{n+1}(r_\d))$, i.e.,
\[
\displaystyle r_{\d^{n+1}} = {\rm Alt}^{n+1}(r_\d) - t_{n+1} \in \d^{n+1} \otimes \d^{n+1}.
\]
 Consequently, 
$J_{2n}: \, \g_{(2n)} \to \d^n: (a_1, a_1, a_2, a_2, \ldots, a_n, a_n) \mapsto (a_1, a_2, \ldots, a_n)$, where $a_1, \ldots, a_n \in \d$, 
is a Lie bialgebra isomorphism from $(\g_{(2n)}, \delta_{\g_{(2n)}})$ to $(\d^n, \delta_{\d^n})$.
\ele

\begin{proof}
For $a \in \d$, and $1 \leq j \leq k \leq n+1$, let
\begin{equation}\lb{eq-ajk}
(a)_{(j,k)} =  (0, \ldots,0,\underset{j^{\text{th}}}{a}, a, \ldots, \underset{k^{\text{th}}}{a},  0, \dots,0) \in \d^{n+1}.
\end{equation}
Assume first that $n = 2q-1$ is odd. If $\{x_i\}_{i=1}^m$ is a basis of $\g$ and $\{\xi_i\}_{i=1}^m$ is the dual basis of $\g^*$, then $\left\{ (x_i)_{(2r-1, 2r)}, \; (\xi_i)_{(2r-1, 2r)}: 1 \leq r \leq q\right\}_{i=1}^m$ is a basis for $\g_{(n+1)}$, with its dual basis of $\g_{(n+1)}^*$ given by $\left\{(\xi_i)_{(1, 2r-1)}, \; -(x_i)_{(2r, n+1)}: \; 1 \leq r \leq q\right\}_{i=1}^m$, and one proves that $r_{\d^{n+1}} = {\rm Alt}^{n+1}(r_\d) - t_{n+1}$ by using the definition of $r_{\d^{n+1}}$ given in \deref{de-r-on-d}. The case when $n$ is even is similarly proved. Combining with \leref{le-uble-mixed-1}, one sees that
the Lie algebra identification $J_{2n}$ is also a Lie bialgebra isomorphism from $(\g_{(2n)}, \delta_{\g_{(2n)}})$ to $(\d^n, \delta_{\d^n})$.
\end{proof}

Analogous to \leref{le-uble-mixed-1}, we have

\ble{le-uble-mixed-2}  Let $n \geq 1$ and let the Lie algebra identifications be as in  
\eqref{eq-iden-ggg}.

1) The Lie bialgebra
$\left(\g_{(2n+1)}, \;\delta_{\g_{(2n+1)}}\right)$ is the mixed twist by the element 
$t_{n+1}\in \wedge^2  \g_{(2n+1)}$ of 
the direct product Lie bialgebra $\left(\g_{(2n+1)}, \;\delta_{\g_{(2n+1)}}^\prime\right)$,  where 
\[
\delta_{\g_{(2n+1)}}^\prime(a_1, \,\ldots, \,a_{n}, \, x) =
(\delta_\d(a_1), \; \ldots, \;\delta_{\d}(a_{n}), \; \delta_\g(x)), \hs a_1, \, \ldots, \, a_n  \in \d, \, x \in \g;
\]

2) The Lie bialgebra
$\left(\g_{(2n)}^*, \;\delta_{\g_{(2n)}^*}\right)$ is the mixed twist by 
the element $-t_{n+1}\in \wedge^2  \g_{(2n)}^*$ of
the direct product Lie bialgebra $\left(\g_{(2n)}^*, \;\delta_{\g_{(2n)}^*}^\prime\right)$, where 
\[
\delta_{\g_{(2n)}^*}^\prime(\xi, \,a_1, \,\ldots, \,a_{n-1}, \,x) = (\delta_{\g^*}(\xi), \;
-\delta_\d(a_1), \; \ldots, -\delta_{\d}(a_{n-1}), \; -\delta_\g(x)), \hs \xi\in \g^*, \,a_1, \ldots, a_{n-1} \in \d, \, x \in \g;
\]

3) The Lie bialgebra
$\left(\g_{(2n+1)}^*, \;\delta_{\g_{(2n+1)}^*}\right)$ is the mixed twist by 
the element $-t_{n+1}\in \wedge^2  \g_{(2n+1)}^*$ of
the direct product Lie bialgebra $\left(\g_{(2n+1)}^*, \;\delta_{\g_{(2n+1)}^*}^\prime\right)$, where 
\[
\delta_{\g_{(2n+1)}^*}^\prime(\xi, \,a_1, \,\ldots, \,a_{n}) = (\delta_{\g^*}(\xi), \;
-\delta_\d(a_1), \; \ldots, \;-\delta_{\d}(a_{n})), \hs \xi\in \g^*, a_1, \, \ldots, \, a_n  \in \d.
\]
\ele

\begin{proof}
One checks directly that the embeddings of direct product Lie algebras
\begin{align*}
&\g_{(2n+1)} \lrw \g_{(2n+2)}, \;\;\; (a_1, \,a_1, \ldots, \,a_{n},\, a_{n}, \,x) \longmapsto
(a_1, \,a_1, \,\ldots, \,a_{n}, \,a_{n}, \,x, \,x),\\ 
&\g_{(2n)}^* \lrw \g_{(2n+2)},  \;\;\; (\xi,\, a_1, \,a_1, \,\ldots,\, a_{n-1}, \,a_{n-1},\, x) \longmapsto
(\xi, \,\,\xi, \,a_1, \,a_1, \,\ldots, \,a_{n-1}, \,a_{n-1}, \,x, \,x),\\
&\g_{(2n+1)}^* \lrw \g_{(2n+2)}, \;\;\; (\xi, \,a_1, \,a_1, \,\ldots, \,a_{n},\, a_{n}) \longmapsto
(\xi, \,\xi, \,a_1, \,a_1, \,\ldots, \,a_{n}, \,a_{n}),
\end{align*}
where $\xi \in \g^*, x \in \g$, and $a_j \in \d$ for $1 \leq j \leq n$, are also respective Lie bialgebra embeddings of 
$\displaystyle\left(\g_{(2n+1)}, \delta_{\g_{(2n+1)}}\!\right)$, 
$\displaystyle\left(\g_{(2n)}^*, -\delta_{\g_{(2n)}^*}\right)$, and 
$\displaystyle\left(\g_{(2n+1)}^*, -\delta_{\g_{(2n+1)}^*}\right)$ into $(\g_{(2n+2)}, \delta_{\g_{(2n+2)}})$.
As $t_{n+1} \in \wedge^2 \g_{(2n+2)}$ actually lies in the images of the above embeddings, \leref{le-uble-mixed-2}
follows directly from \leref{le-uble-mixed-1}.
\end{proof}

\subsection{Mixed product Poisson structures associated to polyubles}\lb{subsec-n-fold} 
Let again $(\g, \delta_\g)$ be a Lie bialgebra, with dual Lie algebra $(\g^*, \delta_{\g^*})$ and double Lie bialgebra $(\d, \delta_\d)$. Assume that $(Y_j, \pi_j, \sigma_j)$, for $1 \leq j \leq n$ and $n \geq 2$, are  left
$(\d, \delta_\d)$-Poisson spaces. Set
\[
\rho_j = -\sigma_j|_{\g^*}: \;\;\; \g^* \lrw \V^1(Y_j) \hs \mbox{and} \hs 
\lam_j = \sigma_j|_\g: \;\;\; \g \lrw \V^1(Y_j), \hs 1 \leq j \leq n,
\]
so that $\rho_j$ and $\lam_j$ are respectively right and left Poisson actions of $(\g^*, \delta_{\g^*})$ and $(\g, \delta_\g)$ on $(Y_j, \pi_j)$.
Let again $\{x_i\}_{i=1}^m$ be any basis of $\g$ and $\{\xi_i\}_{i=1}^m$ the dual basis of $\g^*$. Define the bi-vector field $\piY$ on the product manifold $Y = Y_1 \times Y_2 \times \cdots \times Y_n$ by
\begin{equation}\lb{eq-pi-n}
\piY=(\pi_1, \; \pi_2, \; \cdots, \; \pi_n)
- \sum_{1 \leq j < k \leq n}\; \sum_{i=1}^m (0, \ldots, \underset{j^{\text{th}}\,\text{entry}}{\rho_j(\xi_i)}, \ldots, 0) \wedge 
(0, \ldots, \underset{k^{\text{th}}\,\text{entry}}{\lam_{k}(x_i)}, \ldots, 0),
\end{equation}
where $(\pi_1, \, \pi_2, \, \cdots, \, \pi_n)$ is the direct product Poisson structure on $Y$ (see notation in $\S$\ref{subsec-intro-notation}). Let $\sigma: \d^n \to \V^1(Y)$ be the direct product Lie algebra action, i.e., 
\begin{equation}\lb{eq-sigma-0-0}
\sigma(a_1, \, a_2, \, \ldots, \, a_n) = (\sigma_1(a_1), \;\sigma_2(a_2), \; \ldots, \; \sigma_{n}(a_{n})),\hs a_1, \ldots, a_n \in \d.
\end{equation}
The following \prref{pr-n-mixed-product} is a direct consequence of \leref{le-double-uble-r} and \prref{pr-mixed-mixed}.

\bpr{pr-n-mixed-product}
The triple $(Y, \piY, \sigma)$ is a left $(\d^n, \delta_{\d^n})$-Poisson space. Moreover, for
$1 \leq j < k \leq n$, one has $p_{jk}(\piY) = \pi_j \times_{(\rho_j, \lam_k)} \pi_j$, 
where $p_{jk}: Y \to Y_j \times Y_k$ is the projection to the $j$'th and the $k$'th factor.
\epr

%\bde{de-n-mixed}
%We call $\piY$ in \eqref{eq-pi-n} the {\it mixed product Poisson structure on $Y_1 \times \cdots \times Y_n$ associated to the %Poisson actions $(\sigma_1, \sigma_2, \ldots, \sigma_n)$} and denote it by 
%\[
%(\pi_1 \times \cdots \times \pi_n)_{(\sigma_1, \,\sigma_2, \,\ldots, \,\sigma_{n})}.
%\]
%\ede

\bre{re-n-extended} If $n \geq 2$ and if we are given
$(Y_1, \pi_1, \rho_1)$ as a right $(\g^*, \delta_{\g^*})$-Poisson space, $(Y_n, \pi_n, \lam_n)$ as a left $(\g, \delta_\g)$-Poisson space, and for $2 \leq j \leq n-1$, $(Y_j, \pi_j, \sigma_j)$ as a left $(\d, \delta_\d)$-Poisson space,
one can still define the bi-vector field $\piY$ on $Y = Y_1 \times Y_2 \times \cdots \times Y_n$
as in \eqref{eq-pi-n}, and \leref{le-uble-mixed-2} and  \prref{pr-mixed-mixed} imply that 
$(Y, \piY, \rho)$ is a right $\left(\g_{(2n-2)}^*, \, \delta_{\g_{(2n-2)}^*}\right)$-Poisson space, where
$\rho$ is the right Lie algebra action 
$\rho: \g_{(2n-2)}^*  \to\V^1(Y)$ given by
\begin{equation}\lb{eq-rho-0-0}
\rho(\xi, a_2, a_2, \ldots, a_{n-1}, a_{n-1}, x) = (\rho_1(\xi), \, -\sigma_2(a_2), \; \ldots, \; -\sigma_{n-1}(a_{n-1}), \; 
-\lam_n(x)),
\end{equation}
generalizing the two-fold mixed product Poisson structure construction in \prref{pr-mixed-is-poisson}. Using \leref{le-uble-mixed-2}, one can similarly apply \prref{pr-mixed-mixed} to obtain 
mixed product Poisson structures with Poisson actions by the Lie bialgebras 
$\displaystyle\left(\g_{(2n-1)}^*, \, \delta_{\g_{(2n-1)}^*}\right)$ and 
$\displaystyle\left(\g_{(2n-1)}, \, \delta_{\g_{(2n-1)}}\right)$. 
\ere

\sectionnew{Mixed twists of product quasitriangular Lie bialgebras and mixed product Poisson structures}\lb{sec-mixed-quasi}
In $\S$\ref{sec-mixed-quasi},  
we give a construction of \ certain mixed twists of direct product quasitriangular Lie bialgebras,
generalizing that of polyubles of Lie bialgebras in $\S$\ref{sec-polyuble}.
Throughout $\S$\ref{sec-mixed-quasi}, let $r$ be a quasitriangular $r$-matrix on a Lie algebra $\g$, and let 
$\delta_r: \g \to \wedge^2 \g, \delta_r(x) = \ad_x r$ be the co-bracket on $\g$ (see \eqref{eq-delta-r}). 
The dual map of $\delta_r$ is the Lie bracket on $\g^*$ given by
\begin{equation}\lb{eq-xi-eta-1}
[\xi, \eta]  = \ad_{r_-(\xi)}^* \eta - \ad_{r_+(\eta)}^* \xi =  \ad_{r_+(\xi)}^* \eta -\ad_{r_-(\eta)}^* \xi.
\end{equation}
where $r_+, r_-: \g^* \to \g$ are given in \eqref{eq-r-pm-0}. 
Recall that a (left or right) Poisson space
of the Lie bialgebra $(\g, \delta_r)$ is also called a Poisson space of the quasitriangular Lie bialgebra $(\g, r)$.

%Recall also  from \leref{le-lpm-dual} the 
%dual pair of Lie bialgebras 
%$(\f_-, \delta_r|_{\f_-})$ and $(\f_+, -\delta_r|_{\f_+})$, where $\f_+ = {\rm Im} (r_+) \subset \g$ and
%$\f_- = {\rm Im}(\f_-) \subset \g$. 

\subsection{The quasitriangular $r$-matrix $r^{(n)}$ on $\g^n$}\lb{subsec-mixed-powers}
Let $n \geq 1$ and consider the direct product Lie algebra $\g^n = \g \oplus \cdots \oplus \g$ ($n$-copies). For $X \in \g^{\ot k}$, $k \geq 1$, and $1 \leq j \leq n$, let $(X)_j \in (\g^n)^{\ot k}$ be the image of $X$ under the embedding of $\g$ into $\g^n$ as the $j$'th summand. If $r =\sum_i x_i \ot y_i \in \g \ot \g$, define ${\rm Alt}^n(r) \in \g^n \ot \g^n$ to be the alternating sum
\begin{equation}\lb{eq-Alt-r}
{\rm Alt}^n(r) = (r, \, -r^{21}, \, r, \, -r^{21}, \, \cdots) =\sum_{1 \leq j \leq n,  \; j\; {\rm is \; odd}} (r)_j \;\;+ \sum_{ 1 \leq j \leq n,  \;j \;{\rm is \; even}} (-r^{21})_j,
\end{equation}
where recall that $r^{21} = \sum_i y_i \ot x_i$, and define ${\rm Mix}^n(r) \in \wedge^2 (\g^n)$
\begin{equation}\lb{eq-r-n-def-0}
{\rm Mix}^n(r) = \sum_{1 \leq j < k \leq n} \left({\rm Mix}^n(r)\right)_{j,k},
\end{equation}
where $\left({\rm Mix}^n(r)\right)_{j,k} =\sum_i (y_i)_j \wedge (x_i)_k$ for $1 \leq j < k \leq n$.
The elements $\left({\rm Mix}^n(r)\right)_{j, k} \in \wedge^2 (\g^n)$, where $1 \leq j < k \leq n$, are indeed well-defined:
the linear map $\left({\rm Mix}^n(r)\right)_{j, k}^\#: (\g^*)^n \cong (\g^n)^* \to \g^n$ is given by
\begin{equation}\lb{eq-Mix-r-sharp}
\left({\rm Mix}^n(r)\right)_{j, k}^\#(\xi_1, \ldots, \xi_n) = 
(-r_+(\xi_k))_j + (-r_-(\xi_j))_k, \hs (\xi_1, \ldots, \xi_n) \in (\g^*)^n \cong (\g^n)^*.
\end{equation}
Define now $r^{(n)} \in \g^n \otimes \g^n$ by
\begin{equation}\lb{eq-r-n-def}
r^{(n)} = {\rm Alt}^n(r) - {\rm Mix}^n(r) \in \g^n \otimes \g^n.
\end{equation}
The next \thref{th-uble-r} says that $r^{(n)}$ is a quasitriangular $r$-matrix on $\g^n$. For $n \geq 2$, 
$r^{(n)}$ is thus a mixed twist (see \deref{de-mixed-t-r}) of the quasitriangular $r$-matrix
${\rm Alt}^n(r)$ on $\g^n$.

\bre{re-rd-n}
Consider the quasitriangular $r$-matrix $r_\d$ on the double
Lie algebra $(\d, \lara_\d)$ of an arbitrary Lie bialgebra $(\g, \delta_\g)$ defined by
the Lagrangian splitting $\d = \g + \g^*$.  
By \leref{le-double-uble-r}, $r_\d^{(n)}$ coincides with the quasitriangular $r$-matrix on $\d^n$ defined by
the Lagrangian splitting $\d^n = \g_{(n)} + \g_{(n)}^*$ in \eqref{eq-n-splitting}.
Note also that if $\phi: (\g, r) \to (\g, r')$ is a homomorphism of quasitriangular Lie bialgebras, then
$\phi(r^{(n)}) = (r')^{(n)}$.
\hfill $\diamond$
\ere

\bth{th-uble-r}
For any integer $n \geq 1$, $r^{(n)}$ is a quasitriangular $r$-matrix on the direct product Lie algebra $\g^n$, and the dual of $\delta_{r^{(n)}}: \g^n \to \wedge^2 \g^n$ is the Lie bracket $[\; , \; ]_{r^{(n)}}$ on $(\g^*)^n \cong (\g^n)^*$ given by 
\[
[(\xi_1, \ldots, \xi_n), \;(\eta_1, \ldots, \eta_n)]_{r^{(n)}} = (\zeta_1, \ldots, \zeta_n), 
\]
where $(\xi_1, \ldots, \xi_n), \;(\eta_1, \ldots, \eta_n) \in (\g^*)^n$, and for each $1 \leq j \leq n$,
\begin{equation}\lb{eq-zetaj}
\zeta_j = [\xi_j, \, \eta_j] + \ad_{r_-(\xi_1 + \cdots + \xi_{j-1}) + r_+(\xi_{j+1} + \cdots + \xi_n)}^* \eta_j - 
\ad_{r_-(\eta_1 + \cdots + \eta_{j-1}) + r_+(\eta_{j+1} + \cdots + \eta_n)}^* \xi_j.
\end{equation}
Moreover, for any $1 \leq k \leq m \leq n$, the map $\phi_{m, k}: \g^m \to \g^n$ given by
\begin{equation}\lb{eq-phi-mk}
\phi_{m, k}(x_1, x_2, \ldots, x_m) = (x_1, \, \ldots, x_{k-1}, \, 
\overbrace{x_k, \; x_k, \; \ldots,\; x_k}^{(n-m+1)\, \mbox{terms}}, \, x_{k+1}, \, \ldots,
\, x_m), \hs x_j \in \g,
\end{equation}
is a Lie bialgebra homomorphism from $(\g^m, \delta_{r^{(m)}})$ to $(\g^n, \delta_{r^{(n)}})$.
\eth

\begin{proof} Note that the symmetric part of $r^{(n)}$ is the alternating sum $(s, -s, s, \cdots, (-1)^ns)$, where $s$ is the symmetric part of $r$. To show that $r^{(n)}$ is a quasitriangular $r$-matrix on $\g^n$, we need to show that ${\rm CYB}(r^{(n)}) = 0$. Let $(\d, \delta_\d)$ be the double Lie bialgebra of $(\g, \delta_r)$, and consider
the Lie algebra homomorphism $p_+: \d \to \g$ given by $p_+(x + \xi) = x + r_+(\xi)$ for $x \in \g$ and $\xi \in \g^*$.
By \eqref{eq-r-ppm}, $p_+(r_\d) = r$. Thus $r^{(n)} =p_+\left(r_{\d}^{(n)}\right)$, where 
$p_+(a_1, \ldots, a_n) =(p_+(a_1), \ldots, p_+(a_n)) \in \g^n$ for $(a_1, \ldots, a_n) \in \d^n$. 
As $r_\d^{(n)}$ is a quasitriangular $r$-matrix on $\d^n$ by \reref{re-rd-n}, 
\[
{\rm CYB}(r^{(n)}) = p_+^n\left({\rm CYB}(r_{\d}^{(n)})\right) = 0.
\]

To compute the Lie bracket $[\;\, \; ]_{r^{(n)}}$ on $(\g^n)^* \cong (\g^*)^n$, consider first the case of $n = 2$. 
It is easy to check that 
the two linear maps $(r^{(2)})_\pm: \g^* \oplus \g^* \to \g \oplus \g$ are respectively given by
\begin{align}\lb{eq-r-2-p}
(r^{(2)})_+(\xi_1, \, \xi_2)& = (r_+(\xi_1+\xi_2), \; \,r_-(\xi_1 + \xi_2)), \\
\lb{eq-r-2-m}
(r^{(2)})_-(\xi_1, \, \xi_2) &= (r_-(\xi_1)+r_+(\xi_2),\, \; r_-(\xi_1) + r_+(\xi_2)),
\hs \xi_1, \xi_2 \in \g^*.
\end{align}
It follows that for $\xi_1, \xi_2, \eta_1, \eta_2 \in \g^*$, one has
\[
[(\xi_1, \xi_2), \; (\eta_1, \eta_2)]_{r^{(2)}} = 
\left([\xi_1, \eta_1] + \ad_{r_+(\xi_2)}^*\eta_1 -\ad_{r_+(\eta_2)}^*\xi_1, \;\, 
[\xi_2, \eta_2] + \ad_{r_-(\xi_1)}^*\eta_2 -\ad_{r_-(\eta_1)}^*\xi_2\right).
\]
The formula for  $[\;\, \; ]_{r^{(n)}}$ now follows from the observation that,
as $r^{(n)}$ is a mixed twist of ${\rm Alt}^n(r)$, for each pair $1 \leq j < k \leq n$,
the projection $p_{jk}: \g^n \to \g \oplus \g$ given by
$(x_1, \ldots, x_n) \mapsto (x_j, x_k)$ is a Lie bialgebra
homomorphism from $(\g^n, \delta_{r^{(n)}})$ to $(\g^2, \delta_{r^{(2)}})$. 

Let $1 \leq k \leq m \leq n$. Using \eqref{eq-zetaj}, one checks by direct computations (we omit the details) that the dual map $\phi_{m, k}^*$ of $\phi_{m, k}$, which is given by
\[
\phi_{m, k}^*: \;\; (\g^*)^n \lrw (\g^*)^m, \;\;\;
\phi_{m, k}^*(\xi_1, \ldots, \xi_n) = (\xi_1, \;  \ldots, \; \xi_{k-1}, \;  \xi_k + \cdots + \xi_l, \; \xi_{l+1}, \;\ldots, \;\xi_n), 
\]
where $l = n-m+k$, is a Lie algebra homomorphism from $((\g^*)^n, \, [\; , \; ]_{r^{(n)}})$ to $((\g^*)^m, \, [\; , \; ]_{r^{(m)}})$. 
\end{proof}

\bre{re-EK-III} Starting from the quasitriangular Lie bialgebra $(\g, r)$,  P. Etingof and D. Kazhdan constructed in \cite{EK:III} a Lie bracket on $(\g^*)^n  =\g^* \oplus \cdots \oplus \g^*$, $n \geq 2$, which, together with the co-bracket dual to the direct product Lie bracket on $\g^n$, form a {\it locally factored Lie bialgebra with equal components}, i.e., for each $1 \leq j \leq n$, the $j$'th summand $(\g^*)_j$ is a sub-Lie bialgebra, isomorphic to the dual Lie bialgebra of $(\g, r)$, and $[(\g^*)_i, \, (\g^*)_j] \subset (\g^*)_j \oplus (\g^*)_j$ for each pair $i \neq j$. One checks directly that the Lie bracket $[\; , \; ]_{r^{(n)}}$ given in \thref{th-uble-r} is precisely the Lie bracket on $(\g^*)^n$ in \cite[Proposition 1.9]{EK:III}. %In particular, \thref{th-uble-r} implies that the Lie bialgebra of Etingof and Kazhdan in \cite[Proposition 1.9]{EK:III} is %co-quasitriangular. \hfill $\diamond$
\hfill $\diamond$
\ere 

\bre{re-tau-r}
One may obtain other quasitriangular $r$-matrices on $\g^n$ from $r^{(n)}$. Indeed, let $S_n$ be the permutation group of $\{1, 2, \ldots, n\}$, and for $\tau \in S_n$, define $\phi_\tau \in {\rm Aut}(\g^n)$ by $\phi_\tau(x_1, \ldots,\, x_n) = (x_{\tau^{-1}(1)},  \ldots,  x_{\tau^{-1}(n)})$ for $x_1, \ldots, x_n \in \g$. Then $\phi_\tau(r^{(n)})$ is a quasitriangular $r$-matrix on $\g^n$ for every $\tau \in S_n$. Write $r = \Lam + s$, where $\Lam \in \wedge^2 \g$ and $s \in (S^2(\g))^\g$, so that
\[
r^{(n)} = (\Lam, \Lam, \ldots, \Lam) - {\rm Mix}^n(\Lam) - {\rm Mix}^n(s) + (s, -s, s, \ldots, (-1)^{n-1}s),
\]
where ${\rm Mix}^n(\Lam) \in \wedge^2 \g^n$ and ${\rm Mix}^n(s) \in \wedge^2 \g^n$ are defined using \eqref{eq-r-n-def-0} by replacing $r$ by $\Lam$ and $s$. Then for every $\tau \in S_n$,
\[
\phi_\tau\left(r^{(n)}\right) = (\Lam, \Lam, \ldots, \Lam) - 
\phi_\tau\left({\rm Mix}^n(\Lam)\right) - \phi_\tau\left({\rm Mix}^n(s)\right) + (\varepsilon(1) s, \,
\varepsilon(2)s, \ldots, \,\varepsilon(n) s)
\]
for some sign function $\varepsilon: \{1, 2 \ldots,n\} \to \{1, -1\}$. It is easy to show that $\phi_\tau\left({\rm Mix}^n(\Lam)\right) = {\rm Mix}^n(\Lam)$, while 
\[
\phi_\tau\left({\rm Mix}^n(s)\right) = 
\sum_{1 \leq j < k \leq n,\;\, 
\tau^{-1}(j) < \tau^{-1}(k)} \left({\rm Mix}^n(s)\right)_{j,k} 
-\sum_{1 \leq j < k \leq n, \;\, \tau^{-1}(j) >\tau^{-1}(k)} \left({\rm Mix}^n(s)\right)_{j,k}.
\]
For an arbitrary function $\varepsilon: \{1, 2 \ldots,n\} \to \{1, -1\}$ and any $\tau \in S_n$, define
\begin{equation}\lb{eq-r-ep-tau-n}
r^{(\varepsilon, \tau, n)} = (\varepsilon(1) s, \,
\varepsilon(2)s, \ldots, \,\varepsilon(n) s) + (\Lam, \Lam, \ldots, \Lam) - 
{\rm Mix}^n(\Lam) - \phi_\tau\left({\rm Mix}^n(s)\right)  \in \g^n \ot \g^n.
\end{equation}
By \eqref{eq-CYB-Lam} - \eqref{eq-CYB-Lam-s}, $r^{(\varepsilon, \tau, n)}$ is a quasitriangular $r$-matrix on $\g^n$.
By \deref{de-mixed-twists}, the corresponding Lie bialgebra $(\g^n, \delta_{r^{(\varepsilon, \tau, n)}})$ is a mixed twist of 
the $n$-fold direct product Lie bialgebra $(\g, \delta_r)$ with itself.
As a special case, for any $J = \{j_1, \ldots, j_k\} \subset \{1, \ldots, n\}$ with $j_1 < \cdots < j_k$ and
$p_\sJ: \g^n \to \g^k$ the projection given by $(x_1, \ldots, x_n) \mapsto (x_{j_1}, \ldots, x_{j_k})$, one has
\[
p_\sJ(r^{(n)}) = r^{(\varepsilon, e, k)} = (\varepsilon(1) s, \,
\varepsilon(2)s, \ldots, \,\varepsilon(k) s) + (\Lam, \Lam, \ldots, \Lam) - 
{\rm Mix}^k(r),  
\]
where $e$ is the identity element of $S_k$, and  $\varepsilon(i) = (-1)^{j_i-1}$ for $i = 1, \ldots, k$.
\hfill $\diamond$
\ere

For $n \geq 1$, we now compare the $n$-uble $(\g_{(n)}, \delta_{\g_{(n)}})$ of the Lie bialgebra $(\g, \delta_r)$ 
with the Lie bialgebra $(\g^n, \delta_{r^{(n)}})$.
Consider first the $(2n)$-uble 
$(\g_{(2n)}, \delta_{\g_{(2n)}})$ of $(\g, \delta_r)$, where $n \geq 1$.  
Recall from \leref{le-double-uble-r}  that one has the Lie bialgebra isomorphism 
\[
 \left(\g_{(2n)}, \, \delta_{\g_{(2n)}}\right) \lrw \left(\d^n, \;\delta_{r_\d^{(n)}}\right): \;\;\;
(a_1, \,a_1, \,\ldots, \,a_n, \,a_n) \longmapsto (a_1, \,\ldots, \,a_n), \hs a_j \in \d.
\]
 Let $p_{2n}: \d^n \to \g^{2n}$ be the Lie algebra homomorphism given by
\[
p_{2n}(a_1,\, a_2, \,\ldots, \,a_n) = (p_+(a_1), \, p_+(a_2), \, \ldots, \, p_+(a_n), \, p_-(a_n), \, p_-(a_{n-1}), \, \ldots, \, p_-(a_1)),      \hs a_j \in \d,
\]
where $p_\pm: \d \to \g$ are the Lie algebra homomorphisms given in \eqref{eq-ppm}. 

\ble{le-dn-g2n}
One has $p_{2n}\left(r_\d^{(n)}\right) = r^{(2n)}$. Consequently,  
$J_{2n}: (\g_{(2n)}, \, \delta_{\g_{(2n)}}) \to \left(\g^{2n}, \, \delta_{r^{(2n)}}\right)$ given by 
\[
J_{2n}(a_1, \, a_1, \, \ldots, \, a_n, \, a_n) = (p_+(a_1), \, p_+(a_2), \, \ldots, \, p_+(a_n), \, p_-(a_n), \, p_-(a_{n-1}), \, \ldots, \, p_-(a_1)),      \hs a_j \in \d,
\]
is a Lie bialgebra homomorphism.
\ele

\begin{proof} 
Using \eqref{eq-r-ppm}, one proves that $p_{2n}\left(r_\d^{(n)}\right) = r^{(2n)}$ by a straightforward calculation. 
\end{proof}

For the $(2n-1)$-uble $(\g_{(2n-1)}, \delta_{\g_{(2n-1)}})$ of $(\g, \delta_r)$, where $n \geq 1$, define $J_{2n-1}: \g_{(2n-1)} \to \g^{2n-1}$ by
\[
J_{2n-1}(a_1, \, a_1, \, \ldots, \, a_{n-1}, \, a_{n-1}, \, x) = (p_+(a_1), \, \ldots, \, p_+(a_{n-1}), \, x, \, p_-(a_{n-1}), \, \ldots, \, p_-(a_1)),
\] 
where $a_1, \ldots, a_{n-1} \in \d$ and $x \in \g$.

\ble{le-dn-g2n-1}
The map $J_{2n-1}: (\g_{(2n-1)},\, \delta_{\g_{(2n-1)}}) \to (\g^{2n-1}, \, \delta_{r^{(2n-1)}})$
is a Lie bialgebra homomorphism. 
\ele

\begin{proof} One has $J_{2n-1} = p_{(2n)} \circ J_{2n} \circ I_{2n-1}$, where 
$I_{2n-1}: (\g_{(2n-1)}, \delta_{\g_{(2n-1)}}) \to (\g_{(2n)}, \delta_{\g_{(2n)}})$
is the Lie bialgebra embedding given by
\[
I_{2n-1}(a_1, \,a_1, \ldots, \,a_{n-1},\, a_{n-1}, \,x) \longmapsto
(a_1, \,a_1, \,\ldots, \,a_{n-1}, \,a_{n-1}, \,x, \,x), \hs a_j \in \d, \, x \in \g,
\]
and $p_{(2n)}: (\g^{2n}, \delta_{r^{(2n)}}) \to (\g^{2n-1}, \,\delta_{r^{(2n-1)}})$ is the  
surjective Lie bialgebra homomorphism (see \reref{re-tau-r})
\[
p_{(2n)}(x_1, \, x_2,\, \ldots, \,x_{2n}) = (x_1, \, x_2, \, \ldots, \, x_{n}, \, x_{n+2}, \, \ldots, x_{2n}),\hs
x_j \in \g.
\]
\end{proof}

\bpr{pr-uble-mixed-power}
When $r$ factorizable, the Lie bialgebra homomorphisms $J_n: (\g_{(n)}, \delta_{\g_{(n)}}) \to (\g^n, \delta_{r^{(n)}})$
in \leref{le-dn-g2n} and \leref{le-dn-g2n-1} are Lie bialgebra isomorphisms.
\epr

\begin{proof}
When $r$ is factorizable, the map $\d \to \g \oplus \g, a \mapsto (p_+(a), \, p_-(a))$, $a \in \d$, is 
a Lie algebra isomorphism. Thus $J_n: \g_{(n)} \to \g^n$ is a Lie algebra isomorphism for any $n \geq 1$.
\end{proof}

\subsection{Mixed product Poisson structures associated to quasitriangular $r$-matrices}\lb{subsec-Poi-r-mixed}

Assume that $(Y_j, \pi_j, \lam_j)$, $1 \leq j \leq n$,  are left $(\g, r)$-Poisson spaces. 
Let $Y = Y_1 \times \cdots \times Y_n$ be the product manifold, and let $\lam = (\lam_1, \ldots, \lam_n): \g^n \to \V^1(Y)$ be the 
direct product left Lie algebra action of $\g^n$ on $Y$ given by
\begin{equation}\lb{eq-sigma-gn}
\lam(x_1, \; \ldots, \; x_n) = (\lam_1(x_1),\;  \ldots, \; \lam_n(x_n)), \hs x_j \in \g.
\end{equation}
If $r = \sum_i x_i \otimes y_i \in \g \otimes \g$, define the bi-vector field $\piY$ on $Y$ by  
\begin{align}\lb{eq-piY-Mix-r}
\piY &= (\pi_1, \ldots, \pi_n) +\lam \left({\rm Mix}^n(r)\right)\\
\nonumber
& =(\pi_1, \; \pi_2, \; \cdots, \; \pi_n)
+ \sum_{1 \leq j < k \leq n}\; \sum_i(0, \ldots, \underset{j^{\text{th}}\,\text{entry}}{\lam_j(y_i)}, \ldots, 0) \wedge 
(0, \ldots, \underset{k^{\text{th}}\,\text{entry}}{\lam_{k}(x_i)}, \ldots, 0),
\end{align}
where again $(\pi_1, \; \pi_2, \; \cdots, \; \pi_n)$ is the direct product Poisson structure on $Y = Y_1 \times \cdots \times Y_n$.

\bth{th-piY-Mix-r}
The triple $(Y, \piY, \lam)$ is a left $(\g^n, r^{(n)})$-Poisson space, and $(Y, \, \piY, \lam_{\rm diag})$ is a left $(\g, r)$-Poisson space, where $\lam_{\rm diag}$ is the 
diagonal action
\begin{equation}\lb{eq-lam-diagonal}
\lam_{\rm diag}: \;\; \g \lrw \V^1(Y), \;\;\; \lam_{\rm diag}(x) = (\lam_1(x), \;\ldots, \;\lam_n(x)), \hs x \in \g,
\end{equation}
of $\g$ on $Y$. 
If, in addition, $\pi_j = -\lam_j(r)$ for each $1 \leq j \leq n$ (see \prref{pr-admi-Poi-r}), then $\piY = -\lam(r^{(n)})$.
\eth

\begin{proof}
The first and the third statements of \thref{th-piY-Mix-r} follow directly from \prref{pr-mixed-mixed}.
By \thref{th-uble-r}, for any $1 \leq k \leq m \leq n$, $(Y, \piY, \lam \circ \phi_{m, k})$ is a left $(\g^m, r^{(m)})$-Poisson space, where $\phi_{m,k}: \g^m \to \g^n$ is given in \eqref{eq-phi-mk}. In particular, 
the triple $(Y, \, \piY, \lam_{\rm diag})$ is a left $(\g, r)$-Poisson space.
\end{proof}

\bde{de-mixed-product-r}
Given $(\g, r)$-Poisson spaces $(Y_j, \pi_j, \lam_j)$ for $1 \leq j \leq n$ with $n \geq 2$, we also denote the 
Poisson structure $\piY$ on $Y = Y_1 \times \cdots \times Y_n$ given in \eqref{eq-piY-Mix-r} by
\[
\piY = (\pi_1 \times \cdots \times \pi_n)_{(\lam_1, \ldots, \lam_n)},
\]
and we will refer to the left $(\g, r)$-Poisson space $(Y, \, \piY, \lam_{\rm diag})$,
where $\lam_{\rm diag}$ is the diagonal action of $\g$ on $Y$ given in  \eqref{eq-lam-diagonal}, as the {\it fusion product}
of $(Y_j, \pi_j, \lam_j)$, $1 \leq j \leq n$.
\ede

\bre{re-step-by-step}
Recall from \leref{le-lpm-dual} the dual pair of Lie bialgebras $(\f_-, \delta_r|_{\f_-})$ and $(\f_+, -\delta_r|_{\f_+})$. 
Let $\{x_i\}_{i=1}^m$ be a basis of $\g$ such that $\{x_i\}_{i=1}^l$ is a basis of $\f_-$, and let $\{\xi_i\}_{i=1}^m$ be the dual basis of $\g^*$. Then $\{r_+(\xi_i)\}_{i=1}^l$ is a basis of $\f_+$, dual to the basis $\{x_i\}_{i=1}^l$ of $\f_-$ under the pairing $\lara_{(\f_-, \f_+)}$  in \eqref{eq-pairing-ll}, and
\[
r = \sum_{i=1}^m x_i \otimes r_+(\xi_i) = \sum_{i=1}^l x_i \otimes r_+(\xi_i) \in \f_- \otimes \f_+ \subset \g \otimes \g
\]
(see \reref{re-r-fpm}).
In the context of \thref{th-piY-Mix-r}, it then follows that for any $1 \leq j < k \leq n$, 
the projection  of $\piY$ to $Y_j \times Y_k$ is precisely the mixed product Poisson structure $\pi_j \times_{(-\lam_j|_{\f_+}, \lam_k|_{\f_-})} \pi_k$ formed out of the  right $(\f_+, -\delta_r|_{\f_+})$-Poisson space $(Y_j, \pi_j, -\lam_j|_{\f_+})$
and the
left $(\f_-, \delta_r|_{\f_-})$-Poisson space $(Y_k, \pi_k, \lam_k|_{\f_-})$.
Furthermore, let $1 \leq j \leq n$, and let $Y_{(j)} = Y_1 \times \cdots \times Y_{j}$ and 
$Y_{[j]} = Y_{j+1} \times \cdots \times Y_n$. It then follows from 
\begin{align*}
{\rm Mix}^n(r) &= \left({\rm Mix}^j(r), \; 0\right) + \left(0, \; {\rm Mix}^{n-j}(r)\right) \\
& \hs + \sum_{i=1}^l
(\underbrace{r_+(\xi_i), \;r_+(\xi_i), \;\ldots,\; r_+(\xi_i)}_{j \, \text{terms}},\; 0, \;0, \;\ldots,\; 0) 
\wedge (0, \;0,\; \ldots, \;0, \;\underbrace{x_i, \;x_i, \;\ldots, \;x_i}_{n-j \, \text{terms}})
\end{align*}
that under the obvious identification $Y = Y_{(j)}\times Y_{[j]}$, $\piY$ is also a two-fold mixed product, namely
\[
\piY = \pi_{(j)} \times_{(-\lam_{(j)}|_{\f_+}, \; \lam_{[j]}|_{\f_-})} \pi_{[j]},
\]
where $\pi_{(j)} = (\pi_1 \times \cdots \times \pi_j)_{(\lam_1, \ldots, \lam_j)}, \; $ 
$\pi_{[j]} = (\pi_{j+1} \times \cdots \times \pi_n)_{(\lam_{j+1}, \ldots, \lam_n)}$, and 
\begin{align*}
&\lam_{(j)}: \; \g \lrw V^1(Y_{(j)}), \;\;
\lam_{(j)}(x) = (\lam_1(x), \; \ldots, \; \lam_j(x)), \hs x \in \g, \\ 
&\lam_{[j]}: \; \;\g \lrw V^1(Y_{[j]}), \;
\lam_{[j]}(x) = (\lam_{j+1}(x), \; \ldots, \; \lam_{n}(x)), \hs x \in \g.
\end{align*}
Set, for notational simplicity,  
\begin{align*}
&\lam_{j}^+ = \lam_{j}|_{\f_+}: \; \f_+ \lrw \V^1(Y_j), \hs \;\;
\lam_{(j)}^+ = \lam_{(j)}|_{\f_+}: \; \f_+ \lrw \V^1(Y_{(j)}),\\
&\lam_{j}^- = \lam_j|_{\f_-}: \;
\f_- \lrw \V^1(Y_{j}),\hs \;\; \lam_{[j]}^- = \lam_{[j]}|_{\f_-}: \;
\f_- \lrw \V^1(Y_{([j])}.
\end{align*}
Then $\piY$ can also be realized as the successive $2$-fold mixed products 
\begin{align*}
\piY &= \left(\left(\pi_1\times_{\left(-\lam_1^+, \lam_2^-\right)} \pi_2\right)
\times_{\left(-\lam_{(2)}^+, \lam_3^-\right)} \cdots 
\times_{\left(-\lam_{(n-2)}^+, \lam_{n-1}^-\right)}\pi_{n-1}\right) 
\times_{\left(-\lam_{(n-1)}^+, \lam_n^-\right)} \pi_n\\
& = \pi_1\times_{\left(-\lam_1^+, \lam_{[2]}^-\right)} \left(\pi_2 \times_{\left(-\lam_2^+, \lam_{[3]}^-\right)} \cdots 
\times_{\left(-\lam_{n-2}^+, \lam_{[n-1]}^-\right)} \left(\pi_{n-1} \times_{\left(-\lam_{n-1}^+, \lam_n^-\right)} \pi_n\right)\right).
\end{align*}
  \hfill $\diamond$ 
\ere

We now give a typical example of the construction in \thref{th-piY-Mix-r}.

\bex{ex-class-examples}
Suppose that $G$ is a connected Lie group with Lie algebra $\g$ and let $s \in (S^2\g)^\g$ be the symmetric part of $r$. Suppose that $Q_1, \ldots, Q_n$ are closed Lie subgroups of $G$ such that  the Lie algebra $\q_j$ of $Q_j$ for each $j$ is coisotropic with respect to $s$ (see $\S$\ref{subsec-Poi-r}). Let $\lam_j$ be the left action of $\g$ on $G/Q_j$ induced by the left action of $G$ on $G/Q$ by left translation. By \prref{pr-admi-Poi-r},  $\lam_j(r)$ is a Poisson structure on $G/Q_j$ for each $j$, and by \thref{th-piY-Mix-r}, $\lam(r^{(n)})$ is a  mixed product Poisson structure on the product manifold $G/Q_1 \times \cdots \times G/Q_n$, where $\lam = (\lam_1, \ldots, \lam_n)$ is the direct product action of $\g^n$ on $G/Q_1 \times \cdots \times G/Q_n$.
Note that any twist $r'$ of $r$, having the same symmetric part as $r$, also gives rise to the Poisson
structure $\lam\left((r')^{(n)}\right)$ on $G/Q_1 \times \cdots \times G/Q_n$.  We will return to these examples in $\S$\ref{subsec-flags} for the case when $G$ is a complex semi-simple Lie group.    \hfill $\diamond$ 
\eex

Let $G$ be again any connected Lie group with Lie algebra $\g$. We end $\S$\ref{subsec-Poi-r-mixed} by 
a discussion on the Poisson Lie group
$\left(G^n, \pi_\sG^{(n)}\right)$, where $n \geq 1$, and 
\begin{equation}\lb{eq-piG-n}
\pi_\sG^{(n)} = \left(r^{(n)}\right)^L - \left(r^{(n)}\right)^R,
\end{equation}
with $\left(r^{(n)}\right)^L$ and $\left(r^{(n)}\right)^R$  respectively denoting the left and right invariant $2$-tensor 
fields on $G^n$ with value $r^{(n)}$ at the identity element. Here $r^{(1)} = r$ by convention, and we also write
$\pi_\sG^{(1)} = \piG$. 
By \prref{pr-uble-mixed-power}, when $r$ is factorizable,  $\displaystyle \left(G \times G, \, \pi_\sG^{(2)}\right)$ is  a
Drinfeld double of the Poisson Lie group $(G, \piG)$.

Consider again the pair $((\f_+, -\delta_r|_{\f_+}), \, (\f_-, \delta_r|_{\f_-}))$ of dual Lie bialgebras, where 
$\f_\pm = {\rm Im}(r_\pm) \subset \g$. 
Let $F_+$ and $F_-$ be the connected subgroups of $G$ with Lie algebras $\f_-$ and $\f_+$ respectively, and denote by the same symbol the restrictions of $\piG$ to $F_+$ and $F_-$. Then the Poisson Lie groups $(F_+, -\piG)$ and $(F_-, \piG)$ form a dual pair, and so do $(F_+, \piG)$ and $(F_-^{\rm op}, \piG)$, where $F_-^{\rm op}$  denotes the manifold $F_-$ with the opposite group structure. Consider now the dual pair of Poisson Lie groups 
\begin{equation}\lb{eq-LLs}
(F, \pi_{\sF}) =(F_-, \piG) \times (F_-^{\rm op}, \piG) \hs \mbox{and} \hs 
(F^*, \pi_{\sF^*})= (F_+, -\piG) \times (F_+, \piG).
\end{equation}
By the multiplicativity of the Poisson structure $\piG$, one has the right and left Poisson actions
\begin{align}\lb{eq-rho-F}
\rho:&\;(G, \piG) \times (F^*, \; \pi_{\sF^*})  \longrightarrow (G, \piG),\;\;(g, \; (f_1, f_2))\longmapsto  f_1^{-1}gf_2, 
\;\; g \in G, \, f_1, f_2 \in F_+,\\
\lb{eq-lam-F}
\lambda: & \;(F, \; \pi_{\sF})\times (G, \piG) \longrightarrow (G, \piG),\;\; 
((f_{-1}, f_{-2}), g) \longmapsto f_{-1} g f_{-2}, \;\; g \in G, \, f_{-1}, f_{-2} \in F_-.
\end{align}
For $1 \leq j < k \leq n$, let $p_{jk}: G^n \to G \times G$ be the projection to the $j$'th and the $k$'th factor.

\bpr{pr-mixed-Gn} For any $n \geq 2$, $\pi_\sG^{(n)}$ is a mixed product Poisson structure on $G^n$,
and $p_{jk}\left(\pi_\sG^{(n)}\right) = \pi_\sG^{(2)}$ for every $1 \leq j < k \leq n$;
Moreover, $\pi_\sG^{(2)} = \piG \times_{(\rho, \lam)} \piG$.
\epr

\begin{proof} The first statement  is clear from the definition of $r^{(n)}$. To prove the second one, let $\{x_i\}_{i=1}^m$ be a basis of $\g$ such that $\{x_i\}_{i=1}^k$ is a basis of $\f_-$, and let $\{\xi_i\}_{i=1}^m$ be the dual basis of $\g^*$. By \eqref{eq-r-fpm}, 
one has
\[
\pi_\sG^{(2)} = (\piG, \, \piG) -\sum_{i=1}^k \left((r_+(\xi_i)^L, 0) \wedge (0, x_i^L)+
(r_+(\xi_i)^R, 0) \wedge (0, x_i^R)\right) = \piG \times_{(\rho, \lam)} \piG.
\]
\end{proof}

\subsection{Fusion products of Poisson spaces and quasi-Poisson spaces}\lb{subsec-mixed-fusion}
Write $r = s + \Lam$, where  $s \in (S^2\g)^\g$ and $\Lam \in \wedge^2 \g$, and let $\phi_s = -2{\rm CYB}(s) \in (\wedge^3 \g)^\g$. By \eqref{eq-CYB-Lam} - \eqref{eq-CYB-Lam-s}, $[\Lam, \Lam] = \phi_s$, and  
\[
\phi_s(\xi, \eta, \zeta) = -2
\la \xi,  \;[s^\#(\eta), \, s^\#(\zeta)]\ra, \hs \xi, \eta, \zeta \in \g^*.
\]
Recall from \cite{AA} that a {\it $(\g, \phi_s)$-quasi-Poisson space} is a triple $(Y, Q_\sY, \lam)$, where $Y$ is a manifold, $\lam: \g \to \V^1(Y)$ a left Lie algebra action, and $Q_\sY$ a $\g-$invariant bi-vector field on $Y$  such that
\[
[Q_\sY, Q_\sY] = \lam(\phi_s). 
\]

Let ${\mathcal P}{(\g,\delta_r)}$ be the category of left $(\g, \delta_r)$-Poisson spaces with $\g$-equivariant Poisson morphisms as morphisms, and let ${\mathcal{QP}}{(\g, \phi_s)}$ be the category of  $(\g, \phi_s)$-quasi-Poisson spaces whose morphisms are the $\g$-equivariant morphisms respecting the quasi-Poisson bi-vector fields. The following \leref{le-quasi-twisting} is proved in \cite{Anton-Yvette, David-Severa:quasi-Hamiltonian-groupoids}. Due to our sign conventions, we include a proof for the convenience of the reader.

\ble{le-quasi-twisting} \cite{Anton-Yvette, David-Severa:quasi-Hamiltonian-groupoids} 
One has the equivalence of categories
\begin{equation}\lb{eq-Qm-pim}
{\mathcal P}{(\g,\delta_r)} \; \longleftrightarrow \; {\mathcal{QP}}{(\g, \phi_s)}, \;\;\; (Y, \; \piY, \, \lam) \; 
\longleftrightarrow \; (Y, \; Q_\sY, \;\lam), 
\end{equation}
where $Q_\sY = \piY +\lam(\Lam)$, and the map on morphisms is the identity map.
\ele

\begin{proof} 
Let $Y$ be a manifold with a left action $\lam$ of $\g$, and let $\piY$ and $Q_\sY$ be bi-vector fields on $Y$ related by $Q_\sY = \piY + \lam(\Lam)$. Assume first that $(Y, Q_\sY, \lam)$ is $(\g, \phi_s)$-quasi-Poisson.  Then
\[
[\pi_\sY, \, \pi_\sY] = [Q_\sY, \, Q_\sY] - \lam([\Lam, \Lam]) = \lam(\phi_s - [\Lam, \Lam]) = 0,
\]
and for $x \in \g$,  $[\lam(x), \pi_\sY] = -[\lam(x), \lam(\Lam)] = \lam([x, \Lam]) = \lam(\delta_r(x))$. Thus 
$(Y, \piY, \lam) \in {\mathcal P}{(\g,\delta_r)}$. Conversely, assume that 
$(Y, \, \pi_\sY, \, \lam)\in {\mathcal P}{(\g,\delta_r)}$. Then 
$[\lam(x), \, Q_\sY] = [\lam(x), \, \pi_\sY] -\lam([x, \Lam]) = [\lam(x), \, \pi_\sY] - \lam(\delta_r(x)) = 0$ 
for all $x \in \g$, and 
\[
[Q_\sY, \, Q_\sY] = 2[\lam(\Lam), \, \pi_\sY] -\lam([\Lam, \Lam]) = \lam(2[\Lam, \Lam]-[\Lam, \Lam]) = \lam([\Lam, \Lam]) =
\lam(\phi_s).
\]
Thus $(Y, Q_\sY, \lam) \in {\mathcal{QP}}{(\g, \phi_s)}$. Clearly the two functors in \eqref{eq-Qm-pim} are inverse functors of each other.
\end{proof}

\bex{ex-Qy-0}
If $\lam: \g \to \V^1(Y)$ is a left action of $\g$ on a manifold $Y$ such that $\lam(s) = 0$ (see $\S$\ref{subsec-Poi-r}), then $\lam(\phi_s)=0$ by the proof of \prref{pr-admi-Poi-r},  so $(Y, Q_Y = 0, \lam)$ is a left $(\g, \phi_s)$-quasi-Poisson space, which, under the equivalence of categories in \leref{le-quasi-twisting}, corresponds to the left $(\g, \delta_r)$-Poisson space $(Y, -\lam(r), \lam)$ in \prref{pr-admi-Poi-r}.    \hfill $\diamond$
\eex

\bre{re-two-versions}
Note that if $(Y, Q_\sY, \lam)$ is a $(\g, \phi_s)$-quasi-Poisson space, so is $(Y, -Q_\sY, \lam)$. If $(Y, Q_\sY, \lam)$ corresponds to the $(\g, \delta_r)$-Poisson space $(Y, \piY, \lam)$ as in  \leref{le-quasi-twisting}, then $(Y, -Q_\sY, \lam)$ corresponds to $(Y, \piY-2Q_\sY, \lam)$.  \hfill $\diamond$
\ere

\bre{re-r-rprime}
The category ${\mathcal{QP}}{(\g, \phi_s)}$ depends only on $s \in (S^2\g)^\g$, while 
${\mathcal P}{(\g,\delta_r)}$ depends only on $\Lam \in \wedge^2 \g$ subject to the condition $[\Lam, \Lam] = \phi_s$.
If $r' = r - t$ is a twist of $r$ (see \deref{de-twist-r}), the composition functor 
\[
{\mathcal P}{(\g,\delta_r)} \lrw {\mathcal{QP}}{(\g, \phi_s)} \lrw {\mathcal P}{(\g,\delta_{r'})}: \;\;\; (Y, \, \piY, \, \lam)
\longmapsto (Y, \, \piY + \lam(t), \, \lam)
\]
is the equivalence of categories given in \leref{le-twisting}.
\ere

Let $n \geq 2$ be any integer and let $\phi_s^n = (\phi_s, \phi_s, \ldots, \phi_s) \in \wedge^3(\g^n)$.  Applying \leref{le-quasi-twisting} to the quasitriangular $r$-matrix $r^{(n)}$ on $\g^n$ defined in $\S$\ref{subsec-mixed-powers}, one has the equivalence of categories
\begin{align}\lb{eq-Qn-pin}
{\mathcal P}{(\g^n, \delta_{r^{(n)}})}& \; \longleftrightarrow \; {\mathcal{QP}}{(\g^n, \phi_s^n)},\\
\lb{eq-piY-QY-Gn}
(Y, \; \piY, \, \lam) \; &
\longleftrightarrow \; (Y, \; Q_\sY, \;\lam), \hs \mbox{where} \;\;\; Q_\sY = \piY +\lam(\Lam, \; \ldots, \; \Lam)
-\lam({\rm Mix}^n(r)).
\end{align}
By \thref{th-uble-r}, one has the restriction functor 
\begin{equation}\lb{eq-Res}
{\mathcal P}(\g^n, \delta_{r^{(n)}}) \; \stackrel{{\rm Res}}{\lrw}\; {\mathcal P}(\g, \delta_r), \;\;\; (Y, \;\piY,\; \lam) 
\longmapsto (Y, \; \piY, \;\lam \circ ({\rm diag})_n),
\end{equation}
where $({\rm diag})_n:  \g \lrw \g^n, ({\rm diag})_n(x) = (x, x, \ldots, x)$ for $x \in \g$.
%i.e.,$({\rm diag})_n = \phi_{1, 1}$ in the notation of \eqref{eq-phi-mk}. 
On the other hand, by \cite{Anton-Yvette, AA} (see also \cite{Li-Bland-Severa}), one has the fusion functor 
\begin{equation}\lb{eq-Qn-fus}
{\mathcal{QP}}{(\g^n, \phi_s^n)} \; \stackrel{{\rm Fus}}{\longrightarrow} \; {\mathcal{QP}}{(\g, \phi_s)},\;\;\;
(Y, \; Q_\sY, \; \lam) \longmapsto (Y, \; Q_\sY + \lam({\rm Mix}^n(s)), \; \lam \circ ({\rm diag})_n).
\end{equation}

\bpr{pr-functors-commute}
One has the commutative diagram
\begin{align*}\lb{eq-commute-functors}
{\mathcal P}(\g^n, \;&\delta_{r^{(n)}}) \; \longleftrightarrow \; {\mathcal{QP}}{(\g^n, \phi_s^n)}\\
{\rm Res} \downarrow &\hs \hs \hs \hs \hs \downarrow {\rm Fus}\\
{\mathcal P}(\g, \; &\delta_r) \hs\longleftrightarrow  \hs{\mathcal{QP}}{(\g, \phi_s)},
\end{align*}
where the functors represented by the top and bottom horizontal arrows are respectively given in \eqref{eq-Qn-pin} and \eqref{eq-Qm-pim}.
\epr

\begin{proof}
Using the definitions of the functors, it is enough to show that
\begin{equation}\lb{eq-diagonal-Lam}
({\rm diag})_n(\Lam) = (\Lam, \, \ldots, \Lam) - {\rm Mix}^n(\Lam),
\end{equation}
which is straightforward to check.
\end{proof}

\bco{co-mixed-fusion-0}
Under the equivalence in \eqref{eq-Qm-pim}, the fusion product 
in ${\mathcal P}(\g, \delta_r)$ defined in \deref{de-mixed-product-r} corresponds to the fusion product in ${\mathcal{QP}}(\g, \phi_s)$. 
\eco

\bre{re-fusion-order} 
For each permutation $\tau \in S_n$, the element 
\[
\phi_{\tau}\left({\rm Mix}^n(r)\right) =  {\rm Mix}^n(\Lam) + \phi_{\tau}\left({\rm Mix}^n(s)\right) \in \wedge^2 (\g^n)
\]
is also a mixed twisting element of the $n$-fold direct product of the Lie bialgebra $(\g, \delta_r)$ with itself. Given $(\g, \delta_r)$-Poisson spaces $(Y_j, \pi_j, \lam_j)$ for $1 \leq j \leq n$ with $n \geq 2$, we then have  the left $(\g, \delta_r)$-Poisson space $(Y, \, \pi_{\sY}^\tau, \lam_{\rm diag})$, called the {\it $\tau$-fusion product}, of $\{(Y_j, \pi_j, \lam_j)\}_j$, where $\pi_{\sY}^\tau$ is now the mixed product Poisson structure on $Y = Y_1 \times \cdots \times Y_n$ defined by
\[
\pi_{\sY}^\tau = (\pi_1, \ldots, \pi_n) +\lam \left(\phi_\tau\left({\rm Mix}^n(r)\right)\right).
\]
Replacing the fusion functor ${\rm Fus}$ in \eqref{eq-Qn-fus} by 
\begin{equation}\lb{eq-Qn-fus-tau}
{\mathcal{QP}}{(\g^n, \phi_s^n)} \; \stackrel{{\rm Fus}^\tau}{\;\;\longrightarrow\;\;} \; {\mathcal{QP}}{(\g, \phi_s)},\;\;\;
(Y, \; Q_\sY, \; \lam) \longmapsto (Y, \; Q_\sY + \lam(\phi_\tau\left({\rm Mix}^n(s)\right)), \; \lam \circ ({\rm diag})_n),
\end{equation}
it follows from \eqref{eq-diagonal-Lam} again that under the functor in \leref{le-quasi-twisting}, $\tau$-fusion products of $(\g, \delta_r)$-Poisson spaces correspond to $\tau$-fusion products of $(\g, \phi_s)$-quasi-Poisson spaces. See \cite{Li-Bland-Severa} for examples of $\tau$-fusion products of $(\g, \phi_s)$-quasi-Poisson spaces for non-trivial $\tau \in S_n$.  \hfill $\diamond$
\ere

\sectionnew{Quotient Poisson structures as mixed products}\lb{sec-mixed-quotient}

In $\S$\ref{sec-mixed-quotient}, we present examples that have motivated our introduction to the notion of mixed 
product Poisson structures. Namely, we identify  a class of quotient Poisson manifolds of product Poisson Lie groups
as mixed product Poisson manifolds (see \deref{de-mixed-general}). 

\subsection{Quotient Poisson structures}\lb{subsec-quo-Poi}
Let $(G, \piG)$ be a Poisson Lie group. Recall that a {\it coisotropic subgroup} of $(G, \piG)$ is, by definition, a Lie subgroup $Q$ of $G$ which is also a coisotropic submanifold with respect to $\piG$, i.e., $\piG(q) \in T_qQ \ot T_qG + T_qG \ot T_qQ$ for every $q \in Q$.

\ble{le-coiso-quo}\cite{STS2}
If $(Y, \piY, \rho)$ is a right $(G, \piG)$-Poisson space and $Q$ a closed coisotropic subgroup of $(G, \piG)$ such that $Y/Q$ is a smooth manifold, then the Poisson structure $\piY$ projects to a well-defined Poisson structure on $Y/Q$.
\ele

We will refer to the projection of $\piY$ to $Y/Q$ as the {\it quotient Poisson structure} on $Y/Q$.
As an example, assume that $(Y, \piY)$ is a Poisson manifold with a left Poisson action of a closed Poisson subgroup $(Q, \pi_\sQ=\piG|_\sQ)$ of $(G, \piG)$. Then one has the right Poisson action of the product Poisson Lie group $(Q \times Q, \, \pi_\sQ \times (-\pi_\sQ))$ on the product Poisson manifold $(G \times Y, \, \piG \times \piY)$  by
\[
(g, \,y) \cdot (q_1, \, q_2) = (g q_1, \, q_2^{-1}y), \hs g \in G, \, y \in Y, \, q_1, q_2 \in Q.
\]
As $Q_{\rm diag} = \{(q, q): q \in Q\}$ is a coisotropic subgroup of $(Q \times Q, \, \pi_{\sQ} \times (-\pi_{\sQ}))$, the
Poisson structure $\piG \times \piY$ on $G \times Y$ projects to a well-defined Poisson structure
on the quotient of $G \times Y$ by $Q_{\rm diag}$, denoted by $G\times_Q Y$.  Repeating the construction, one has the following \ldref{ld-GGG-QQQ}.

\bld{ld-GGG-QQQ} 
If $Q_1, \ldots, Q_n$ are closed Poisson subgroups of $(G, \piG)$, the quotient manifold
\[
Z = G \times_{Q_1} \cdots \times_{Q_{n-1}} G/Q_n
\]
in \eqref{eq-Z} has the well-defined Poisson structure $\pi_\sZ \!= \!\varpi_\sZ(\pi_\sG^n)$, where $\varpi_\sZ: \!G^n \!\to Z$ is the projection. We will refer to $(Z, \pi_\sZ)$ as a quotient Poisson manifold of the product Poisson Lie group
$(G^n, \pi_\sG^n)$. 
\eld

We now look at the case when $n = 1$, leaving the case of an arbitrary $n$ to $\S$\ref{subsec-mixed-actions}. Assume 
thus $Q$ is a closed Poisson Lie subgroup of $(G, \piG)$ and consider the Poisson structure
\[
\pi_{\sGQ} := \varpi_{{\scriptscriptstyle G/Q}} (\piG)
\]
on $G/Q$, where $\varpi_{{\scriptscriptstyle G/Q}}: G \to G/Q$ is the projection. Let $(\g, \delta_\g)$ be the Lie bialgebra of $(G, \piG)$, and let $(\g^*, \delta_{\g^*})$ and $(\d, \delta_\d)$ be respectively the dual and double Lie bialgebras of $(\g, \delta_\g)$. Recall from \eqref{eq-varsig-G} the left Poisson action $\varsigma$ of $(\d, \delta_\d)$ on $(G, \piG)$, i.e.,
\[
\varsigma: \;\; \d \lrw \V^1(G), \;\; \varsigma(x+\xi) = x^R -\pi_\sG^\#(\xi^R), \hs x \in \g, \, \xi \in \g^*.
\]
As $Q$ is a Poisson subgroup, it follows from \eqref{eq-multi-vr} that
\begin{equation}\lb{eq-vrGQ}
\vr_\sGQ: \; \g^* \lrw \V^1(G/Q), \;\;\; \vr_\sGQ(\xi) \,\stackrel{{\rm def}}{=}\, \varpi_{\sGQ}(\pi_\sG^\#(\xi^R)), \hs \xi \in \g^*,
\end{equation}
is a well-defined Poisson action, called the {\it right dressing action},  of $(\g^*, \delta_{\g^*})$ on $(G/Q, \pi_\sGQ)$. As
\begin{equation}\lb{eq-lamGQ}
\lam_\sGQ: \; \g \lrw \V^1(G/Q), \;\;\; \lam_\sGQ(x) \,\stackrel{{\rm def}}{=}\, \varpi_{\sGQ}(x^R), \hs x \in \g,
\end{equation}
is also well-defined, one has the well-defined left Lie algebra action
\begin{equation}\lb{eq-sigma-GQ}
\varsigma_\sGQ = \varpi_{\sGQ} \circ \varsigma: \;\; \d \lrw \V^1(G/Q), \;\;\; \varsigma_\sGQ (x+\xi) = 
\varpi_{\sGQ}(x^R -\pi_\sG^\#(\xi^R)), \ \hs
x \in \g, \, \xi \in \g^*,
\end{equation}
of $\d$ on $G/Q$, making $(G/Q, \pi_\sGQ, \varsigma_\sGQ)$ into a left $(\d, \delta_\d)$-Poisson space. Recall 
from \eqref{eq-double-r} the quasitriangular $r$-matrix $r_\d$ on $\d$ associated to the Lagrangian splitting $\d = \g + \g^*$. 

\ble{le-GQ-d}
For any closed Poisson Lie subgroup $Q$ of $(G, \piG)$, one has $\pi_\sGQ =- \varsigma_\sGQ(r_\d)$.
\ele

\begin{proof}  
Let $\{x_i\}_{i=1}^m$ be any basis of $\g$ and $\{\xi_i\}_{i=1}^m$ the dual basis of $\g^*$. By the definition of the bundle map $\pi_\sG^\#: T^*G \to TG$, one has
\[
\piG = \sum_{i=1}^m x_i^R \ot \pi_\sG^\#(\xi_i^R) = -\sum_{i=1}^m \varsigma(x_i) \ot \varsigma(\xi_i) = -\varsigma(r_\d).
\]
Composing with $\varpi_{\sGQ}: G \to G/Q$ gives
\begin{equation}\lb{eq-piSGQ-var}
\pi_{\sGQ} = \sum_{i=1}^m \varpi_{\sGQ}(x_i^R) \ot \varpi_\sGQ(\pi_\sG^\#(\xi_i^R))= -\varsigma_{\sGQ} (r_\d).
\end{equation}
\end{proof}

For the remainder of $\S$\ref{subsec-quo-Poi}, let $(\g, r)$ be  quasitriangular Lie bialgebra, and let $G$ be any connected Lie group with Lie algebra $\g$. One then has the Poisson Lie group $(G, \piG)$, where $\piG = r^L -r^R$. We now consider some particular examples of Poisson Lie subgroups of $(G, \piG)$. We first compute the right dressing action $\vr$ of $(\g^*, \delta_{\g^*})$ on $(G, \piG)$. Recall from $\S$\ref{subsec-r-matrices}
that $r_+ = r^\#$ and $r_- = -(r^{21})^\# = -r_+^*$ are  Lie bialgebra homomorphisms from $(\g^*, -\delta_{\g^*})$ to $(\g, \delta_r)$. For $g \in G$, let $L_g$ and $R_g$ denote again  the respective left and right translations on $G$ by $g$.

\ble{le-dressing-quasi}
For $\xi \in \g^*$ and $g \in G$, one has
\begin{equation}\lb{eq-dressing-quasi}
\pi_\sG^\#(\xi^R)(g) =L_g r_+(\Ad_g^*\xi) - R_g r_+(\xi) =  L_g r_-(\Ad_g^*\xi) - R_g r_-(\xi).
\end{equation}
\ele

\begin{proof}
Writing $r =\sum_j x_i \ot y_j$ where $x_j, y_j \in \g$, one has
\begin{align*}
\pi_\sG^\#(\xi^R)(g) &= L_g\left(\sum_j \la R_{g^{-1}}^* \xi, \, L_g x_j\ra y_j\right) 
- R_g\left(\sum_j \la R_{g^{-1}}^* \xi, \, R_g x_j\ra y_j\right)\\
& = L_g r_+(\Ad_g^*\xi) - R_g r_+(\xi).
\end{align*}
The second identity in \eqref{eq-dressing-quasi} is proved using $\piG(g) = L_g (-r^{21}) - R_g (-r^{21})$.
\end{proof}

Consider again (see $\S$\ref{subsec-r-matrices}) the sub-Lie bialgebras $(\f_-, \delta_r|_{\f_-})$ and $(\f_+, \delta_r|_{\f_+})$ of $(\g, \delta_r)$, where $\f_\pm = {\rm Im} (r_\pm)$.

\bpr{pr-l-pm-q}
Let $Q$ be a Lie subgroup of $G$ with Lie algebra $\q \subset \g$ such that $\f_+ \subset \q$ or $\f_- \subset \q$. Then $Q$ is a Poisson Lie subgroup of $(G, \piG)$. Moreover, if $Q \subset G$ is closed, then
\begin{align}\lb{eq-vr-GQ-1} 
\vr_\sGQ(\xi) &= -\lam_\sGQ (r_+(\xi)), \hs \xi \in \g^*, \hs \mbox{if} \; \; \f_+ \subset \q;\\
\lb{eq-vr-GQ-2}
\vr_\sGQ(\xi) &= -\lam_\sGQ (r_-(\xi)), \hs \xi \in \g^*, \hs \mbox{if} \; \f_- \subset \q.
\end{align}
In both cases, one has $\pi_\sGQ = -\lam_\sGQ(r)$.
%If $Q \subset G$ is closed and $\f_- \subset \q$, then
%\begin{equation}
%\vr_\sGQ(\xi) = -\lam_\sGQ (r_-(\xi)), \hs \xi \in \g^*;
%\end{equation}
\epr

\begin{proof}
As $\q\supset \f_+$ or $\q \supset \f_-$, it follows from \leref{le-dressing-quasi} that $\pi_\sG^\#(\xi^R)(g)\in T_gQ$ for every $\xi \in \g^*$ and $g \in Q$. Thus $Q$ is a Poisson Lie subgroup of $(G, \piG)$.

Assume that $Q$ is closed in $G$. Then \eqref{eq-vr-GQ-1} and \eqref{eq-vr-GQ-2} follow from \eqref{eq-dressing-quasi}. Assume that $\f_+ \subset \q$. Let $\{x_i\}_{i=1}^m$  be basis of $\g$ and $\{\xi_i\}_{i = 1}^m$ the dual basis of $\g^*$ such that $\{x_i\}_{i=1}^k$ is a basis of $\f_-$. Then ${\rm Span}\{\xi_{k+1}, \ldots, \xi_m\} = \f_-^0 = \ker r_+$. By \eqref{eq-piSGQ-var} and \eqref{eq-vr-GQ-1}, 
\[
\pi_\sGQ = \sum_{i=1}^k \lam_{\sGQ}(x_i) \ot \vr_\sGQ(\xi_i) =-\lam_\sGQ 
\left(\sum_{i=1}^k x_i \ot r_+(\xi_i)\right) = -\lam_\sGQ(r).
\]
If $\f_- \subset \q$, one has $\pi_\sGQ = -\lam_\sGQ(-r^{21})= -\lam_\sGQ(r)$.
\end{proof}

%Studying this class of quotient Poisson structures is 
%one of the main motivations for this paper.

\subsection{Quotient Poisson structures as mixed products}\lb{subsec-mixed-actions}
In $\S$\ref{subsec-mixed-actions}, fix a connected Poisson Lie group $(G, \piG)$ and let $(\g, \delta_\g)$ be its Lie bialgebra. Let $(\g^*, \delta_{\g^*})$ and $(\d, \delta_\d)$ be respectively the dual and the double Lie bialgebra of $(\g, \delta_\g)$. Recall the right dressing action $\vr$ of $(\g^*, \delta_{\g^*})$ on $(G, \piG)$ given in \eqref{eq-vr-p}.

Let $(Y, \piY, \lam)$ be a left $(G, \piG)$-Poisson space. Let $J$ be the diffeomorphism
\[
J: \;\; G \times Y \lrw G \times Y, \;\; (g, y) \longmapsto (g, \; gy), \hs g \in G,\; y \in Y.
\]
Recall from $\S$\ref{subsec-mixed-product} that one has the mixed product Poisson structure 
$\piG \times_{(\varrho, \lam)} \piY$ on $G \times Y$.

\ble{le-GGY}
One has $J(\piG \times \piY) = \piG \times_{(\varrho, \lam)} \piY$. 
\ele

\begin{proof} 
Let $\pi = J(\piG \times \piY)$. With $\kappa_\sG: G \times Y \to G$ and $\kappa_\sY: G \times Y \to Y$ denoting the two projections, it is clear that $\kappa_\sG: (G \times Y, \pi) \to (G, \piG)$ and $\kappa_\sY: (G \times Y, \pi) \to (Y, \piY)$ are Poisson maps. Let $\mu$ be the mixed term of $\piG \times_{(\varrho, \lam)} \piY$ and let $g \in G$, $y \in Y$. It remains to show that
\[
\pi(g, y)(\kappa_\sG^* \theta_g, \; \kappa_\sY^* \zeta_y) = \mu(g, y)(\kappa_\sG^* \theta_g, \; \kappa_\sY^* \zeta_y), 
\hs \theta_g \in T_g^*G, \; \zeta_y \in T^*_yY.
\] 
Let $y_1 = g^{-1}y \in Y$, and let $\lam_{y_1}: G \to Y, \lam_{y_1}(h) = hy_1$ for $h \in G$. Let $\iota_g: Y \to G \times Y$ and $\iota_{y_1}: G \to G \times Y$ be respectively given by $\iota_g(y') =(g, y')$ and $\iota_{y_1}(g') = (g', y_1)$. 
Then
\begin{align*}
\pi(g, y)(\kappa_\sG^* \theta_g, \; \kappa_\sY^* \zeta_y)
&=\piG(g)(\iota_{y_1}^*J^*\kappa_\sG^* \theta_g, \; \iota_{y_1}^*J^*\kappa_\sY^* \zeta_y)
+\piY(y_1)(\iota_{g}^*J^*\kappa_\sG^* \theta_g, \; \iota_{g}^*J^*\kappa_\sY^* \zeta_y)\\
& = \piG(g)(\theta_g, \; \lam_{y_1}^* \zeta_y)-\la \theta_g, \; \varrho(r_g^* \lam_{y_1}^* \zeta_y)\ra 
=\mu(g, y)(\kappa_\sG^* \theta_g, \; \kappa_\sY^* \zeta_y).
\end{align*}
This proves that $\pi = \piG \times_{(\varrho, \lam)} \piY$.
\end{proof}

Assume now that $Q$ is a closed Poisson Lie subgroup of $(G, \piG)$, and consider the quotient manifold $G \times_\sQ Y$ with the quotient Poisson structure $\pi=\varpi(\piG \times \piY)$, where $\varpi: G \times Y \to G \times_\sQ Y$ is the  projection. Let $J_\sQ$ be the diffeomorphism
\[
J_\sQ: \;\; G \times_\sQ Y \lrw (G/Q) \times Y, \;\; J_\sQ(\varpi(g, y)) = (gQ, \; gy), \hs g \in G, \;y \in Y.
\]
Recall the right Poisson action $\vr_\sGQ$  of the Lie bialgebra $(\g^*, \delta_{\g^*})$ on $(G/Q, \pi_\sGQ)$ given in \eqref{eq-vrGQ}.

\ble{le-mixed-GQY} 
As Poisson structures on $(G/Q) \times Y$, one has 
\begin{equation}\lb{eq-pi-GQY}
J_\sQ(\pi) =  \pi_{{\scriptscriptstyle G/Q}} 
\times_{(\varrho_{\scriptscriptstyle G/Q}, \; \lam)} \piY.
\end{equation}
\ele

\begin{proof} 
\leref{le-mixed-GQY} follows from \leref{le-GGY} and the commutative diagram
\begin{align*}
&G \times Y \;\;\stackrel{J}{\longrightarrow} \;\; G\times Y\\
&\;\; \varpi \downarrow   \hs \hs \hs \hs\downarrow \varpi_{\scriptscriptstyle G/Q} \times {\rm Id}_\sY\\
&G \times_\sQ Y \;\; \stackrel{J_\sQ}{\longrightarrow} \;\; G/Q\times Y.
\end{align*}
\end{proof}

Let now $Q_1, \ldots, Q_n$ be closed Poisson Lie subgroups of $(G, \piG)$ and consider the quotient 
\[
Z = G \times_{Q_1} G\times_{Q_2} \cdots \times_{Q_{n-1}} G/Q_n
\]
defined in \eqref{eq-Z} with the quotient Poisson structure $\pi_{\sZ} = \varpi_\sZ(\pi_\sG^n)$, where $\varpi_\sZ: G^n \to Z$ is the projection. Set $[g_1, \ldots, g_n]_\sZ = \varpi_\sZ(g_1, \ldots, g_n) \in Z$ for $(g_1, \ldots, g_n) \in G^n$
and recall the diffeomorphism $J_\sZ: Z \to G/Q_1  \times \cdots \times G/Q_n$ given in \eqref{eq-IZ}, i.e.,
\[
J_{\sZ}[g_1, g_2, \ldots, g_n]_\sZ = (g_1Q_1, \; g_1g_2Q_2, \;\ldots, \;g_1g_2\cdots g_nQ_n), \hs g_j \in G.
\]
By \leref{le-GQ-d}, for each $1 \leq j \leq n$, one has the $(\d, \delta_\d)$-Poisson space $(G/Q_j, \pi_{{\sGQ}_j}, \varsigma_j)$, where $\varsigma_j = \varsigma_{{\scriptscriptstyle G/Q_j}}$ (see \eqref{eq-sigma-GQ}), and  $\pi_{{\sGQ}_j} = -\varsigma_j(r_\d)$, where $r_\d \in \d \ot \d$ is the quasitriangular $r$-matrix on $\d$ defined by the Lagrangian splitting $\d = \g + \g^*$. Consider the direct product left Lie algebra action $\sigma = (\varsigma_1, \; \ldots, \; \varsigma_n)$ of $\d^n$ on the product manifold $G/Q_1 \times \cdots \times G/Q_n$.

\bth{th-GQ-n-mixed}
One has (see notation in \deref{de-mixed-product-r})
\begin{equation}\lb{eq-IZ-0}
J_{\sZ}(\pi_\sZ)= (\pi_{{\scriptscriptstyle G/Q_1}} \times \cdots \times \pi_{{\scriptscriptstyle G/Q_n}})_{(\varsigma_1, \; \ldots, \; \varsigma_n)}= -\sigma\left(r_\d^{(n)}\right).
\end{equation} 
\eth

\begin{proof}
When $n = 1$, 
\eqref{eq-IZ-0} holds by \leref{le-GQ-d}. Assume that $n > 1$ and that \eqref{eq-IZ-0} holds for 
$n-1$. By taking $Q = Q_1$ and $Y = G \times_{Q_2} \cdots \times_{Q_{n-1}} G/Q_n$ in \leref{le-mixed-GQY} and by noting that
${\rm Im}((r_\d)_+) = \g^* \subset \d$ and ${\rm Im}((r_\d)_-) = \g\subset \d$, \eqref{eq-IZ-0} holds 
by \leref{le-mixed-GQY} and \reref{re-step-by-step}.
\end{proof}

\bre{re-special-pairs-one} ({\bf Reduction to sub-Lie bialgebras, II.}) 
For a Poisson Lie group $(G, \piG)$ with Lie bialgebra $(\g, \delta_\g)$,  a pair $((\p, \delta_\p), \,  Q)$, where $(\p, \delta_\p = \delta_\g|_\p)$ is a sub-Lie bialgebra of $(\g, \delta_\g)$ and $Q$ is a closed Poisson Lie subgroup of $(G, \piG)$, is said to be {\it special} if $\vr_{\sGQ}(\xi) = 0$ for all $\xi \in \p^0$, where $\vr_\sGQ$ is given in \eqref{eq-vrGQ}. In this case, the left Poisson action $\varsigma_{\sGQ}$ of $(\d, \delta_\d)$ on $(G/Q, \pi_\sGQ)$ induces a left action, denoted by $\varsigma_{\p}$, of the double Lie bialgebra of $(\p, \delta_\p)$ on $(G/Q, \pi_\sGQ)$ by $\varsigma_\p (x+\xi + \p^0) = \varsigma_\sGQ(x+\xi)$,  $x \in \p, \xi \in \g^*$ (see \exref{ex-double-of-sub} and \reref{re-mixed-15}), and in the setting of \thref{th-GQ-n-mixed} with $Q_j = Q$ for each $j$, one checks directly that the Poisson structure $J_\sZ(\pi_\sZ)$ on $G/Q \times \cdots \times G/Q$ is also given by
\[
J_\sZ(\pi_\sZ) = (\pi_{{\scriptscriptstyle G/Q}} \times \cdots \times 
\pi_{{\scriptscriptstyle G/Q}})_{(\varsigma_{\p}, \ldots, \varsigma_{\p})}  = -\sigma_{\p}
\left(r_{\d_\p}^{(n)}\right),
\]
where $\sigma_{\p} = (\varsigma_{\p}, \ldots, \varsigma_{\p})$ is the direct product action of $\d_\p$ on $Z$ and $r_{\d_\p}$ is the quasitriangular $r$-matrix on $\d_\p$ associated to the Lagrangian splitting $\d_\p = \p + \p^*$.
\hfill$\diamond$
\ere

\sectionnew{Primary examples}\lb{sec-primary-examples}
In $\S$\ref{sec-primary-examples}, 
we identify a specific class of mixed Poisson
structures, which are the main motivation for the work in the present paper. 
Examples include those defined on products of flag varieties and related spaces that are mentioned in the Introduction.
 
\subsection{The set-up}\lb{subsec-setup} Throughout $\S$\ref{sec-primary-examples}, we fix a quasitriangular Lie bialgebra
$(\g, r)$, and let $G$ be any connected Poisson Lie group with Lie algebra $\g$. One then has the Poisson Lie groups
$(G, \piG)$ and $\displaystyle \left(G \!\times \!G, \, \pi_\sG^{(2)}\right)$, where $\piG = r^L - r^R$, and
$\pi_\sG^{(2)} = (r^{(2)})^L - (r^{(2)})^R$ (see $\S$\ref{subsec-Poi-r-mixed}).

Recall the two Lie subalgebras $\f_\pm = {\rm Im} (r_\pm)$ of $\g$.
Let $(Q_+, Q_-)$ be a pair of connected and closed Lie subgroups of $G$ such that their respective Lie algebras
$\q_+$ and $\q_-$ satisfy 
\begin{equation}\lb{eq-ff-qq}
\f_+ \subset \q_+ \hs \mbox{and} \hs \f_- \subset \q_-.
\end{equation}  
By \prref{pr-l-pm-q}, both $Q_+$ and $Q_-$ are Poisson Lie subgroups of $(G, \piG)$. 
By \eqref{eq-r-2-p}, 
$\q_+ \oplus \q_- \supset {\rm Im}\left((r^{(2)})_+\right)$, so  
$Q_+ \times Q_-$ is a Poisson Lie subgroup of $\displaystyle \left(\!G \!\times\! G, \, \pi_\sG^{(2)}\right)$
again by \prref{pr-l-pm-q}. One thus has the following six series of 
quotient Poisson manifolds (see \ldref{ld-GGG-QQQ})
\begin{equation}\lb{eq-series-1}
(X_n, \, \pi_{\sX_n}), \hs (X_{-n}, \, \pi_{\sX_{-n}}), 
\hs 
(Y_n, \, \pi_{\sY_n}), \hs (Y_{-n}, \, \pi_{\sY_{-n}}), \hs  ({\mathbb X}_n, \; \pi_{\scriptscriptstyle{\mathbb{X}_n}}), \hs ({\mathbb Y}_n, \; \pi_{\scriptscriptstyle{\mathbb{Y}_n}})
\end{equation}
of the product Poisson Lie groups $(G^n, \, \pi_\sG^n)$ and  $((G \times G)^n, \,(\pi_\sG^{(2)})^n)$,  
where $n \geq 2$ is any integer,   and 
\begin{align}\lb{eq-D-nn-0}
&X_n = G \times_{Q_+} \cdots \times_{Q_+} G, \hs\hs X_{-n} = G \times_{Q_-} \cdots \times_{Q_-} G,\\
\lb{eq-X-nn-0}
&Y_n = G \times_{Q_+} \cdots \times_{Q_+} G/Q_+, \hs \hs Y_{-n} = G \times_{Q_-} \cdots \times_{Q_-} G/Q_-,\\
\lb{eq-DD-n-0}
&{\mathbb X}_n = (G \times G)\times_{(Q_+ \times Q_-)} \cdots \times_{(Q_+ \times Q_-)}  (G \times G),\\
\lb{eq-XX-n-0}
&{\mathbb Y}_n  = (G \times G)\times_{(Q_+ \times Q_-)} \cdots \times_{(Q_+ \times Q_-)}  (G \times G)/(Q_+ \times Q_-).
\end{align}
\thref{th-GQ-n-mixed} applies to any $(Z, \pi_\sZ)$ in the series in \eqref{eq-series-1} and identifies
the Poisson structures $\pi_\sZ$ as a mixed product Poisson structure defined by  a quasitriangular
$r$-matrix. We show in  $\S$\ref{sec-primary-examples} that the assumptions on $\q_+$ and $\q_-$ in \eqref{eq-ff-qq} allow
one to 
``reduce" the quasitriangular $r$-matrix in \thref{th-GQ-n-mixed} in the same spirit as the reduction of the
$r$-matrix in \leref{le-GQ-d} to that in \prref{pr-l-pm-q} 
for the case of $n = 1$.

\subsection{The Poisson manifolds $(X_n, \, \pi_{\sX_n})$, 
$(X_{-n}, \, \pi_{\sX_{-n}})$}\lb{subsec-DD} 
Let $n \geq 1$ and consider first the Poisson manifold $(X_n, \, \pi_{\sX_n})$.
Define
the direct sum quasitriangular $r$-matrix $r^{\la n+1\ra}$ on $\g^{n+1} = \g^n \oplus \g$ by
\begin{equation}\lb{eq-r-modify}
r^{\la n+1\ra} = \begin{cases} (r^{(n)}, \, 0) + (0, \, -r), & \;\; \mbox{if} \;\; n = 2m+1 \;  \mbox{is odd},\\
(r^{(n)}, \, 0) + (0, \, r^{21}), & \;\; \mbox{if} \;\; n = 2m \; \mbox{is even},\end{cases}
\end{equation}
and let $\lam$ be the left Lie algebra action of $\g^{n+1}$ on  the product manifold $(G/Q_+)^{n-1} \times G$ by 
\[
\lam(x_1, \ldots, x_n, x_{n+1}) = (\lam_\sGQp(x_1), \; \ldots, \;\lam_{\sGQp}(x_{n-1}), \; x_n^R-x_{n+1}^L), \hs x_j \in \g.
\]
Let $J_{\sX_n}: X_n \to (G/Q_+)^{n-1} \times G$ be the diffeomorphism given by
\[
J_{\sX_n}([g_1, g_2, \cdots, g_n]_{\sX_n}) = \left(g_1Q_+, \; g_1g_2Q_+, \, \ldots, g_1g_2 \cdots g_{n-1}Q_+, \; g_1g_2 \cdots g_n\right),
\hs g_j \in G.
\]

\bpr{pr-QQQG}
As Poisson structures on $(G/Q_+)^{n-1} \times G$, one has
$J_{\sX_n}(\pi_{\sX_n}) = -\lam\left(r^{\la n+1 \ra}\right)$.
\epr

\begin{proof} We apply \thref{th-GQ-n-mixed} by taking $Q_1 = \cdots = Q_{n-1} = Q_+$ and $Q_n = \{e\}$. 
Let the notation be as in \thref{th-GQ-n-mixed}. 
Let $\{x_i\}_{i=1}^m$ be a basis of $\g$ such that $\{x_i\}_{i=1}^k$
is a basis of $\f_-$, and let $\{\xi_i\}_{i=1}^m$ be the dual basis of $\g^*$. Then $r_\d = \sum_{i=1}^m x_i \otimes \xi_i$.
By \thref{th-GQ-n-mixed}, 
\[
J_{\sX_n}(\pi_{\sX_n}) = (\pi_{\sGQp}, \; \ldots, \; \pi_{\sGQp}, \; \piG) + \sigma \left({\rm Mix}^n\left(r_\d\right)\right).
\]
Let $1 \leq j < k \leq n$ and consider (see \eqref{eq-r-n-def-0})
\[
\mu_{jk} =\sigma \left(\left({\rm Mix}^n\left(r_\d\right)\right)_{jk}\right) = \sum_{i=1}^m 
(0, \ldots, \underset{j^{\text{th}}\,\text{entry}}{\sigma_{\sGQp}(\xi_i)}, \ldots, 0) \wedge 
(0, \ldots, \underset{k^{\text{th}}\,\text{entry}}{\sigma_{\sGQp}(x_i)}, \ldots, 0).
\]
Using the fact that $\f_+ \subset \q_+$ and by \prref{pr-l-pm-q}, one has
\[
\mu_{jk} = \sum_{i=1}^k 
(0, \ldots, \underset{j^{\text{th}}\,\text{entry}}{\lam_{\sGQp}(r_+(\xi_i))}, \ldots, 0) \wedge 
(0, \ldots, \underset{k^{\text{th}}\,\text{entry}}{\lam_{\sGQp}(x_i)}, \ldots, 0)
= \lam\left(\left({\rm Mix}^n(r)\right)_{jk}\right).
\]
Thus
\[
J_{\sX_n}(\pi_{\sX_n}) = (\pi_{\sGQp}, \; \ldots, \; \pi_{\sGQp}, \; \piG) + 
\lam \left(\left({\rm Mix}^n\left(r\right), \; 0\right)\right).
\]
On the other hand, again by \prref{pr-l-pm-q}, one has
\[
-\lam\left(\left(r^{(n)}, 0\right)\right) = \begin{cases} 
(\pi_{\sGQp}, \; \ldots, \; \pi_{\sGQp}, \; -r^R) + \lam \left(\left({\rm Mix}^n\left(r\right), \; 0\right)\right),& \;\; 
\mbox{if} \;\; n \;\; \mbox{is odd},\\
(\pi_{\sGQp}, \; \ldots, \; \pi_{\sGQp}, \; (r^{21})^R) +\lam \left(\left({\rm Mix}^n\left(r\right), \; 0\right)\right),& \;\; 
\mbox{if} \;\; n \;\; \mbox{is even},\end{cases}
\]
Let $r' = -r$ if $n$ is odd and $r' = r^{21}$ if $n$ is even. Then
$-\lam (0, \, r') =  r^L$ if $n$ is odd and 
$-\lam (0, \, r') =-(r^{21})^L$ if $n$ is even. 
Since $\piG = r^L - r^R = (r^{21})^R - (r^{21})^L$, one has
\[
-\lam\left(r^{\la n+1 \ra}\right) = 
-\lam\left(\left(r^{(n)}, 0\right)\right)-\lam (0, \, r') = J_{\sX_n}(\pi_{\sX_n}).
\] 
\end{proof}

Replacing $r$ by $-r^{21}$, we have
\[
J_{\sX_{-n}}(\pi_{\sX_{-n}}) = -\lam\left((-r^{21})^{\la n+1 \ra}\right),
\]
where $J_{\sX_{-n}}$ is the diffeomorphism $X_{-n} \to (G/Q_-)^{n-1} \times G$ given by
\[
J_{\sX_{-n}}([g_1, g_2, \cdots, g_n]_{\sX_{-n}}) = \left(g_1Q_-, \; g_1g_2Q_-, \, \ldots, g_1g_2 \cdots g_{n-1}Q_-, \; g_1g_2 \cdots g_n\right),
\hs g_j \in G.
\]

\subsection{The Poisson manifolds $(Y_{n}, \, \pi_{\sY_{\scriptscriptstyle{n}}})$ and $(Y_{-n}, \, \pi_{\sY_{\scriptscriptstyle{-n}}})$}\lb{subsec-XX}
Let $n \geq 1$, and let
$J_{\sY_n}: Y_n \to (G/Q_+)^n$ and $J_{\sY_{-n}}: Y_{-n} \to (G/Q_-)^n$ be the diffeomorphisms, respectively given by
\begin{align*}
&J_{\sY_{n}}([g_1, g_2, \cdots, g_n]_{\sY_{n}}) = \left(g_1Q_+, \; g_1g_2Q_+, \, \ldots,  \; g_1g_2 \cdots g_nQ_+\right),
\hs g_j \in G, \\
&J_{\sY_{-n}}([g_1, g_2, \cdots, g_n]_{\sY_{-n}}) = \left(g_1Q_-, \; g_1g_2Q_-, \, \ldots,  \; g_1g_2 \cdots g_nQ_-\right),
\hs g_j \in G.
\end{align*}
Using the surjective Poisson morphisms
\begin{align}\lb{eq-XD-1}
&(X_n, \,\pi_{\sX_n}) \lrw (Y_n, \,\pi_{\sY_n}), \;\; [g_1, g_2, \ldots, g_n]_{\sX_n} \longmapsto [g_1, g_2, \ldots, g_n]_{\sY_n},
\hs g_j \in G,\\
\lb{eq-XD-2}&(X_{-n}, \,\pi_{\sX_{-n}}) \lrw (Y_{-n}, \,\pi_{\sY_{-n}}), \;\; [g_1, g_2, \ldots, g_n]_{\sX_{-n}} \longmapsto [g_1, g_2, \ldots, g_n]_{\sY_{-n}},
\hs g_j \in G,
\end{align}
one has the following direct consequence of \prref{pr-QQQG}.

\bth{th-main-GGG-QQQ}
As Poisson structures on $(G/Q_+)^n$ and $(G/Q_-)^n$ respectively, one has
\[
J_{\sY_n}(\pi_{\sY_n}) = -\lam\left(r^{(n)}\right) \hs \mbox{and} \hs
J_{\sY_{-n}}(\pi_{\sY_{-n}}) = -\lam\left((-r^{21})^{(n)}\right),
\]
where the action  $\lam$  of $\g^n $ on $(G/Q_\pm)^n$ is the direct product of $\lam_{\scriptscriptstyle{G/Q_{\pm}}}$
of $\g$ 
on each factor.
\eth

\subsection{The Poisson structures $\pi_{\scriptscriptstyle{\mathbb{X}_n}}$ and 
$\pi_{\scriptscriptstyle{\mathbb{Y}_n}}$ as two-fold mixed products}\lb{subsec-DD_XX}
Let  $n \geq 1$, and consider the diffeomorphism
$S_{\sXbb_n}: {\mathbb{X}}_n \to X_n \times X_{-n}$ given by
\begin{equation}\lb{eq-S-n}
S_{\sXbb_n}([g_1, k_2, g_2, k_2, \ldots, g_n, k_n]_{\scriptscriptstyle{\mathbb{X}_n}}) =
\left([g_1, g_2, \ldots, g_n]_{\scriptscriptstyle{{X}_n}}, \; [k_1, k_2, \ldots, k_n]_{\scriptscriptstyle{{X}_{-n}}}\right),
\hs g_j, k_j \in G.
\end{equation}
We now express the Poisson structure $S_{\sXbb_n}(\pi_{\scriptscriptstyle{\mathbb{X}_n}})$ on $X_n \times X_{-n}$ as a mixed product.
To this end, recall from \eqref{eq-LLs} the pair of dual Poisson Lie groups 
\[
(F, \pi_{\sF}) =(F_-, \piG) \times (F_-^{\rm op}, \piG) \hs \mbox{and} \hs 
(F^*, \pi_{\sF^*})= (F_+, -\piG) \times (F_+, \piG),
\]
where $F_+$ and $F_-$ are the connected Lie subgroup of $G$ with Lie algebras $\f_+$ and $\f_-$ respectively.
Consider the right Poisson action $\rho_{\sX_n}$ of $(F^*, \; \pi_{\sF^*})$ on $(X_n, \pi_{\scriptscriptstyle{{X}_n}})$
and the left Poisson action $\lam_{\sX_{-n}}$ of $(F, \; \pi_{\sF})$ on 
$(X_{-n}, \pi_{\scriptscriptstyle{{X}_{-n}}})$, respectively given by
\begin{align*}
&\rho_{\sX_n}([g_1, g_2, \ldots, g_n]_{\scriptscriptstyle{{X}_n}}, \; (f_1, f_2))=  [f_1^{-1}g_1, \,g_2, \,\ldots,\, g_nf_2]_{\scriptscriptstyle{{X}_n}}, 
\hs g_j \in G, \; f_1, f_2 \in F_+,\\
&\lam_{\sX_{-n}}((f_{-1}, f_{-2}), \; [g_1, g_2, \ldots, g_n]_{\scriptscriptstyle{{X}_{-n}}}) =
[f_{-1} g_1, \, g_2, \, \ldots,\, g_n f_{-2}]_{\scriptscriptstyle{{X}_{-n}}}, \hs g_j \in G, \; f_{-1}, f_{-2} \in F_-.
\end{align*}

\bpr{pr-Dnnn}
One has $S_{\sXbb_n}(\pi_{\scriptscriptstyle{\mathbb{X}_n}}) = \pi_{\scriptscriptstyle{{X}_n}} \times_{(\rho_{\sX_n}, \lam_{\sX_{-n}})}
\pi_{\scriptscriptstyle{{X}_{-n}}}$.
\epr

\begin{proof}
Let $S$ be the diffeomorphism $(G \times G)^n \lrw G^n \times G^n$ given by
\[
S(g_1, k_2, g_2, k_2, \ldots, g_n, k_n) =
((g_1, g_2, \ldots, g_n), \, (k_1, k_2,\ldots, k_n)), \hs g_j, k_j \in G.
\]
By \prref{pr-mixed-Gn}, $\displaystyle S\left(\left(\pi_\sG^{(2)}\right)^n\right) = 
\pi_\sG^n \times_{(\rho^n, \, \lam^n)} \pi_\sG^n$, where $\rho$ is the right Poisson action of
$(F^*, \; \pi_{\sF^*})$ on $(G, \piG)$ given in \eqref{eq-rho-F}, $\lam$ is the left Poisson action of 
$(F, \pi_\sF)$ on $(G, \piG)$ given in \eqref{eq-lam-F}, and $\rho^n$ and $\lam^n$ denote respectively the
direct product actions of the direct product Poisson Lie groups $((F^*)^n, \; \pi_{\sF^*}^n)$ and
$(F^n, \pi_\sF^n)$ on $(G^n, \pi_\sG^n)$. Let $\varpi_n$ and $\varpi_{-n}$ be 
respectively the projections from $G^n$ to $X_n$ and $X_{-n}$. Then  
\[
S_{\sXbb_n}(\pi_{\scriptscriptstyle{\mathbb{X}_n}}) = (\varpi_n \times\varpi_{-n})
\left(S\left(\left(\pi_\sG^{(2)}\right)^n\right)\right)= 
(\varpi_n \times\varpi_{-n})\left(\pi_\sG^n \times_{(\rho^n, \, \lam^n)} \pi_\sG^n
\right).
\]
Consider the
product Poisson Lie group $(M, \pi_{\sM}) = (F_-^{\rm op}, \piG) \times (F_-, \piG)$ and its dual Poisson
Lie group $(M^*, \pi_{\sM^*}) = (F_+,\piG) \times (F_+, -\piG)$, and note the direct product decompositions
\begin{align}\lb{eq-LM-1}
(F^n, \; \pi_{\sF}^n) &= (F_-, \piG) \times (M^{n-1}, \; \pi_{\sM}^{n-1}) \times (F_-^{\rm op}, \piG),\\
\lb{eq-LM-2}
((F^*)^n, \; \pi_{\sF^*}^n) &= (F_+, -\piG) \times ((M^*)^{n-1}, \; \pi_{\sM^*}^{n-1}) \times (F_+, \piG).
\end{align}
The action $\lam^n$ of $(F^n, \; \pi_{\sF}^n)$ on $(G^n, \pi_\sG^n)$ restricts to left Poisson
actions by the Poisson subgroups 
\begin{align*}
(M^{n-1}, \; \pi_{\sM}^{n-1})& \cong \{e\} \times (M^{n-1}, \; \pi_{\sM}^{n-1}) \times \{e\} \subset (F^n, \; \pi_{\sF}^n),\\
(F, \pi_{\sF}) &\cong (F_-, \piG) \times \{e\} \times (F_-^{\rm op}, \piG) \subset (F^n, \; \pi_{\sF}^n),
\end{align*} 
which will be respectively denoted by $\lam^{(n-1)}$ and $\lam^\prime$. Similarly, $\rho^n$ restricts to right Poisson
actions $\rho^{(n-1)}$ of $((M^*)^{n-1}, \; \pi_{\sM^*}^{n-1})$ on $(G^n, \pi_\sG^n)$ 
and $\rho^\prime$ of $(F^*, \pi_{\sF^*})= (F_+, -\piG) \times (F_+, \piG)$
on $(G^n, \pi_\sG^n)$ via \eqref{eq-LM-2}. We then have
$\pi_\sG^n \times_{(\rho^n, \lam^n)} \pi_\sG^n= \pi_\sG^n \times_{(\rho^\prime, \lam^\prime)} \pi_\sG^n
+ \mu$, where $\mu$ is the mixed term of the mixed product Poisson structure $\pi_\sG^n \times_{(\rho^{(n-1)}, \lam^{(n-1)})} \pi_\sG^n$. Thus
\[
S_{\sXbb_n}(\pi_{\scriptscriptstyle{\mathbb{X}_n}}) = (\varpi_n \times\varpi_{-n})
\left(\pi_\sG^n \times_{(\rho^\prime, \lam^\prime)} \pi_\sG^n\right) + (\varpi_n \times\varpi_{-n})(\mu)
=\pi_{\scriptscriptstyle{{X}_n}} \times_{(\rho_{\sX_n}, \lam_{\sX_{-n}})}
\pi_{\scriptscriptstyle{{X}_{-n}}} + (\varpi_n \times\varpi_{-n})(\mu).
\]
It remains to show that $(\varpi_n \times\varpi_{-n})(\mu)=0$.

Consider the coisotropic subgroup $A = \{(f_-, f_-^{-1}): f_- \in F_-\}$ of $(M, \pi_{\sM})$ and the 
coisotropic subgroup $B = \{(f, f): f \in F_+\}$ of 
$(M^*, \pi_{\sM^*})$. It is easy to see that the annihilator of the Lie algebra of $A$ in the Lie algebra of $M^*$
is precisely the Lie algebra of $B$. Thus the annihilator of the Lie algebra of $A^{n-1}$ in the Lie algebra
of $(M^*)^{n-1}$ is precisely the Lie algebra of $B^{n-1}$. Since $X_n = G^n/B^{n-1}$, where $B^{n-1}$ acts on 
$G^n$ by $\rho^{(n-1)}$ via $B^{n-1} \subset (M^*)^{n-1}$, and  $X_{-n} = A^{n-1} \backslash G^n$, where $A^{n-1}$ acts on $G^n$
by $\lam^{(n-1)}$ via $A^{n-1} \subset M^{n-1}$, it follows from \leref{le-vanish} that $(\varpi_n \times\varpi_{-n})(\mu)=0$.
\end{proof}

Similarly, define the diffeomorphism $S_{\sYbb_n}: {\mathbb{Y}}_n \to Y_n \times Y_{-n}$ by
\begin{equation}\lb{eq-X-n}
S_{\sYbb_n}([g_1, k_2, g_2, k_2, \ldots, g_n, k_n]_{\scriptscriptstyle{\mathbb{Y}_n}}) =
\left([g_1, g_2, \ldots, g_n]_{\scriptscriptstyle{{Y}_n}}, \; [k_1, k_2, \ldots, k_n]_{\scriptscriptstyle{{Y}_{-n}}}\right),
\hs g_j, k_j \in G.
\end{equation}
Define the  right Poisson action $\rho_{\sY_n}$ of $(F_+, -\piG)$ on $(Y_n, \pi_{\scriptscriptstyle{{Y}_n}})$
and the left Poisson action $\lam_{\sY_{-n}}$ of $(F_-,  \piG)$ on 
$(Y_{-n}, \pi_{\scriptscriptstyle{{Y}_{-n}}})$, respectively given by
\begin{align*}
&\rho_{\sY_n}([g_1, g_2, \ldots, g_n]_{\scriptscriptstyle{{Y}_n}}, \; f)=  [f^{-1}g_1, \,g_2, \,\ldots,\, g_n]_{\scriptscriptstyle{{Y}_n}}, 
\hs g_j \in G, \; f \in F_+,\\
&\lam_{\sY_{-n}}(f_-, \; [g_1, g_2, \ldots, g_n]_{\scriptscriptstyle{{Y}_{-n}}}) =
[f_- g_1, \, g_2, \, \ldots,\, g_n]_{\scriptscriptstyle{{Y}_{-n}}}, \hs g_j \in G, \; f_- \in F_-.
\end{align*}
Using again the Poisson morphisms in \eqref{eq-XD-1} - \eqref{eq-XD-2}, one has the 
following consequence of \prref{pr-Dnnn}.

\bpr{pr-Xnnn}
One has $S_{\sYbb_n}(\pi_{\scriptscriptstyle{\mathbb{Y}_n}}) = \pi_{\scriptscriptstyle{{Y}_n}} \times_{(\rho_{\sY_n}, \lam_{\sY_{-n}})}
\pi_{\scriptscriptstyle{{Y}_{-n}}}$.
\epr

\subsection{A weak Poisson dual pair}\lb{subsec-poi-dual-pair}

\bde{de-weak-pair}
A {\it weak Poisson dual pair} is a pair of surjective Poisson submersions
\begin{equation}\lb{eq-rho-YZ}
\Phi_\sX: \;(Z, \, \piZ) \lrw (X, \,\piX) \hs \mbox{and} \hs  \; \Phi_\sY: \; (Z, \,\piZ) \lrw (Y, \,\piY),
\end{equation}
between Poisson manifolds such that 
the  map
\[
(\Phi_\sX, \; \Phi_\sY): \;\; (Z, \,\piZ) \lrw (X \times Y, \; \pi_\sX \times \piY), \;\;\; z \longmapsto
(\Phi_\sX(z), \; \Phi_\sY(z)), \hs z \in Z,
\]
is Poisson, where $\piX \times \piY$ is the product Poisson structure on $X \times Y$.
\ede

\bre{re-dualpairs}
When $(Z, \piZ)$ is symplectic and when the tangent spaces to the
fibers of $\Phi_{\sX}$ and $\Phi_{\sY}$ are the symplectic orthogonals of each other, the pair 
$(\Phi_{\sX}, \, \Phi_{\sY})$ is called a 
{\it symplectic dual pair}. 
Our notion of {\it weak Poisson dual pairs} thus generalizes that of symplectic dual pairs.
\hfill $\diamond$
\ere

Let the setting, in particular, the pair of subgroups $(Q_+, Q_-)$ of $G$, be as in $\S$\ref{subsec-setup}.
Similar to the Poisson manifold $(Y_{-n}, \pi_{\sY_{-n}})$, where $n \geq 1$, introduce the left action of $Q_-^n$ on $G^n$ by
\[
(q_{-1}, \, q_{-2}, \, \ldots, \, q_{-n}) \cdot (g_1, g_2, \ldots, g_n) = (q_{-1} g_1 q_{-2}^{-1}, \;q_{-2} g_2 q_{-3}^{-1}, \; \ldots, \; q_{-n} g_n), \hs q_{-j} \in Q_-, \, g_j \in G,
\]
and denote by 
\begin{equation}\lb{eq-minus-X}
Y_{-n}^\prime = Q_- \backslash G \times_{Q_-} G \times_{Q_-} \cdots \times_{Q_-} G
\end{equation}
the corresponding quotient of $G^n$ by $Q_-^n$, and by 
\[
\Phi_-: \;\; G^n \lrw Y_{-n}^\prime,\;\;\; (g_1, g_2, \ldots, g_n) \longmapsto [g_1, \, g_2, \ldots, g_n]_{\sY_{-n}^\prime}, \hs g_j \in G,
\]
the corresponding quotient map. Then 
\[
\pi_{\sY_{-n}^\prime} : =\Phi_-(\pi_{\sG}^n)
\]
is a well-defined Poisson structure on $Y_{-n}^\prime$.
Let $\Phi_+: G^n \to X_n$ be the projection from $G^n$ to $Y_n$, i.e., 
$\Phi_+(g_1, \ldots, g_n) = [g_1, \ldots, g_n]_{\sY_n}$ for $(g_1, \ldots, g_n) \in G^n$.

\bpr{pr-poi-dual-pair}
For any $n \geq 1$, the pair 
\[
\Phi_+: \;\; (G^n, \, \pi_{\sG}^n) \lrw (Y_n, \, \pi_{\sY_n}) \hs \mbox{and} \hs
\Phi_-: \;\; (G^n, \, \pi_{\sG}^n) \lrw \left(Y_{-n}^\prime, \, \pi_{\sY_{-n}^\prime}\right) 
\]
is a weak Poisson dual pair.
\epr 

\begin{proof} By definition, both $\Phi_+$ and $\Phi_-$ are surjective Poisson submersions.
Denote the Poisson submersion in \eqref{eq-XD-1} by $\phi_+: 
(X_n, \pi_{\sX_n}) \to (Y_n, \, \pi_{\sY_n})$, and consider the Poisson submersion
\[
\phi_-: \;\; (X_{-n}, \, \pi_{\sX_{-n}}) \lrw (Y_{-n}^\prime, \, \pi_{\sY_{-n}^\prime}), \;\; [g_1, \ldots, g_n]_{\sX_{-n}}
\longmapsto [g_1, \ldots, g_n]_{\sY_{-n}^\prime}, \hs g_j \in G.
\]
By \prref{pr-Dnnn}, one has the Poisson morphism
\[
S_{\sXbb_n}: \;\; (\Xbb_n, \, \pi_{\sXbb_n}) \lrw 
(X_n \times X_{-n}, \; \,\pi_{\scriptscriptstyle{{X}_n}} \times_{(\rho_{\sX_n}, \lam_{\sX_{-n}})}
\pi_{\scriptscriptstyle{{X}_{-n}}}).
\]
By the definition the the projections $\phi_+$ and $\phi_-$, the mixed term in the mixed product 
Poisson structure $\pi_{\scriptscriptstyle{{X}_n}} \times_{(\rho_{\sX_n}, \lam_{\sX_{-n}})}
\pi_{\scriptscriptstyle{{X}_{-n}}}$ on $X_n \times X_{-n}$ vanishes under $\phi_+ \times  \phi_-:
X_n \times X_{-n} \to Y_n \times Y_{-n}^\prime$. Thus one has the Poisson morphism
\[
(\phi_+ \times \phi_-) \circ S_{\sXbb_n}: \;\;  (\Xbb_n, \, \pi_{\sXbb_n}) \lrw
\left(Y_n \times Y_{-n}^\prime, \;\, \pi_{\sY_n} \times \pi_{\sY_{-n}^\prime}\right).
\]
As $\displaystyle (G, \piG) \to \left(G \times G, \; \pi_\sG^{(2)}\right), \,g \mapsto (g, g)$ is Poisson,
one has the Poisson  embedding
\[
I_n: \;\; (G, \, \pi_\sG^n) \lrw \left((G \times G)^n, \, \left(\pi_\sG^{(2)}\right)^n\right),\;\;
(g_1, \,g_2, \,\ldots,\,g_n) \longmapsto (g_1, \,g_1,\, g_2, \,g_2, \,\ldots, \,g_n, \,g_n),
\hs g_j \in G.
\]
With $\Phi: \left((G \times G)^n, \, \left(\pi_\sG^{(2)}\right)^n\right) \to (\Xbb_n, \pi_{\sXbb_n})$ denoting the projection,
one sees that 
\[
(\Phi_+, \Phi_-) = (\phi_+ \times \phi_-)\circ S_{\sXbb_n} \circ \Phi\circ I_n: \;\; (G^n, \pi_\sG^n) \lrw
\left(Y_n \times Y_{-n}^\prime, \; \,\pi_{\sY_n} \times \pi_{\sY_{-n}^\prime}\right)
\]
is a Poisson morphism.
\end{proof}

\subsection{Holomorphic Poisson structures on flag varieties}\lb{subsec-flags}
Let $G$ be a connected complex semi-simple Lie group with Lie algebra $\g$. Recall that a flag variety of $G$ is a (left) $G$-homogeneous space with parabolic stabilizer subalgebras. Let $\F_1, \ldots, \F_n$ be any finite collection of flag varieties of $G$, and for each $1 \leq j \leq n$, let $\lam_j$ be the corresponding (left) Lie algebra action of $\g$ on $\F_j$. 

\bpr{pr-flags}
For any quasitriangular $r$-matrix $r$ on $\g$, $\lam(r^{(n)})$ is a mixed product Poisson structure on the product manifold $\F= \F_1 \times \cdots \times \F_n$,  where $\lam =(\lam_1, \ldots, \lam_n)$ is the direct product action of $\g^n $ on ${\mathcal F}$.
\epr

\begin{proof}
Let $s \in (S^2\g)^\g$ be the symmetric part of $r$. If $\g = \g_1 \oplus \cdots \oplus \g_k$ is the decomposition of $\g$ into simple factors, then \cite[Lemma 1]{Delorme:Manin}, $s = (s_1, \ldots, s_k)$, where, $s_i \in S^2(\g_i)^{\g_i}$, $1 \leq i \leq k$, is a scalar multiple of the element in $S^2(\g_i)^{\g_i}$ corresponding to the Killing form of $\g_i$. Thus every parabolic subalgebra of $\g$, being a direct sum of parabolic subalgebras of the $\g_i$'s, is coisotropic with respect to  $s$. \prref{pr-flags} now follows from  \thref{th-piY-Mix-r} (see also \exref{ex-class-examples}).
\end{proof}

We remark that all factorizable quasitriangular $r$-matrices on a complex simple Lie algebra have been classified by
Belavin-Drinfeld in \cite{BD}. Assume again that $\g$ is semi-simple.
To give more concrete examples of the Poisson structures on flag varieties and related spaces, 
fix a symmetric non-degenerate invariant
bilinear form $\lara_\g$ on $\g$, and fix also  a pair $(\b, \b_-)$ of opposite Borel subalgebras of $\g$.  Let
$\h = \b \cap \b_-$, a Cartan subalgebra of $\g$. Let $\Delta$ and $\Delta_+\subset\Delta$ be respectively the set of roots 
for the pair $(\g, \h)$ and $(\b, \h)$, and let 
$\displaystyle \g = \h + \sum_{\alpha \in \Delta_+} \g_{\alpha} + \sum_{\alpha \in \Delta_+} \g_{-\alpha}$
be the corresponding root decomposition. 
Let $\{h_i\}_{i=1}^r$ be a basis of $\h$ such that $\la h_i, h_j\ra_\g = \delta_{ij}$
for $1 \leq i, j \leq \dim \h$, and let $E_\alpha \in \g_\alpha$ and $E_{-\alpha} \in \g_{-\alpha}$ be root vectors for 
$\alpha \in \Delta_+$ such that $\la E_\alpha, E_{-\alpha}\ra_\g = 1$. The {\it standard quasitriangular $r$-matrix}
associated to the choice of the triple $(\b, \b_-, \lara_\g)$ is the element $r_{\rm st} \in \g \otimes \g$
given by
\begin{equation}\lb{eq-r-st-0}
r_{\rm st} = \frac{1}{2} \sum_{i=1}^{\dim \h} h_i \ot h_i + \sum_{\alpha \in \Delta_+} E_{-\alpha} \otimes E_\alpha.
\end{equation}
The Poisson Lie group $(G, \pist)$, where $\pist = r_{\rm st}^L - r_{\rm st}^R$, is called
a {\it standard complex semi-simple Lie group}, and it has $(G \times G, \, \Pist)$ as  
a Drinfeld double Poisson Lie group, where  $\Pist$ is the Poisson structure on 
the product group $G \times G$ given by $\Pist = (r_{\rm st}^{(2)})^L - 
(r_{\rm st}^{(2)})^R$. 

It is clear from the definition in \eqref{eq-de-l-pm} that  Lie subalgebras $\f_-$
and $\f_+$ of $\g$ associated to $r_{\rm st}$ are respectively given by
$\f_- = \b_-$ and $\f_+ = \b$. Let $B$ and $B_-$ be the subgroups of $G$ with Lie algebras $\b$ and $\b_-$ respectively.
Taking $Q_+ = B$ and $Q_- = B_-$ and applying the constructions in 
$\S$\ref{subsec-setup}, one arrives at the four series of complex manifolds in $\S$\ref{subsec-moti-intro}, respectively 
denoted as
\[
(F_n, \; \,\pi_{\sF_n} \!\!= \pi_n), \hs (\Fbb_n, \;\, \pi_{\sFbb_n}\! \!= \Pi_n), 
\hs (\wF_n, \;\, \pi_{\swF_n} \!\!= \tilde{\pi}_n),
\hs (\wFbb_n, \;\, \pi_{\swFbb_n}\!\! = \widetilde{\Pi}_n).
\]

Using results in $\S$\ref{subsec-DD} and $\S$\ref{subsec-XX}, we can thus identify the four series of
Poisson manifolds with mixed product Poisson manifolds, but more importantly, their mixed product Poisson structures are now
defined by quasitriangular $r$-matrices.
In \cite{Lu-Victor:flags}, a sequel to the current paper, we
first develop a general theory on torus orbits of symplectic leaves for Poisson structures 
defined by quasitriangular $r$-matrices and then apply the general theory to obtain descriptions of 
the $T$- orbits of symplectic leaves for $(F_n, \, \pi_n), \, (\Fbb_n, \, \Pi_n), \, (\wF_n, \, \tilde{\pi}_n)$ and 
$(\wFbb_n, \, \widetilde{\Pi}_n)$, where $T = B \cap B_-$, a maximal torus of $G$. See \cite{Lu-Victor:flags} for detail.
 
For $n \geq 1$, introduce also the quotient manifold 
\[
F_{-n}^\prime = B_- \backslash G \times_{B_-} G \times_{B_-}\cdots \times_{B_-} G
\]
of $G^n$ as in \eqref{eq-minus-X} and denote by $\pi_{-n}^\prime$ projection of $\pi_{\rm st}^n$ from $G^n$ to 
$F_{-n}^\prime$. Let $\Phi_+$ and $\Phi_-$ be the projections
\begin{align*}
\Phi_+: &\;\; G^n \lrw F_n, \;\; (g_1, g_2, \ldots, g_n) \longmapsto [g_1, g_2, \ldots, g_n]_{\sF_n},
\\
\Phi_-: &\;\; G^n \lrw F_{-n}^\prime, \;\; (g_1, g_2, \ldots, g_n) \longmapsto [g_1, g_2, \ldots, g_n]_{\sF_{-n}^\prime},
\hs g_j \in G.
\end{align*}
The following \prref{pr-BGB} is a special case of \prref{pr-poi-dual-pair} and is used in \cite{Balazs-Lu:BS}.

\bpr{pr-BGB}
For any $n \geq 1$, the pair
\[
\Phi_+: \;\; (G^n, \, \pi_{\rm st}^n) \lrw (F_n, \, \pi_n) \hs \mbox{and} \hs
\Phi_-: \;\; (G^n, \, \pi_{\rm st}^n) \lrw \left(F_{-n}^\prime, \; \pi_{-n}^\prime\right) 
\]
is a weak Poisson dual pair.
\epr 
%%%%%%%%%%%%Bibliography%%%%%%%%%%%%%%%%%%%%%%%

\end{document}